\documentclass[11pt,reqno, a4paper]{amsart}
\usepackage{marginnote}

\usepackage{textcomp,microtype,lmodern,mathrsfs,fbb,scalerel,stackengine}
\usepackage{todonotes}
\usepackage{amsmath}
\usepackage{amsfonts}  
\usepackage{amssymb}
\usepackage{amsthm}
\usepackage{url}
\usepackage{graphicx}
\usepackage[noadjust]{cite}
\usepackage{caption}
\usepackage{subfig}
\usepackage{dutchcal}
\usepackage{bbm}
\usepackage{tikz}
\usetikzlibrary{decorations.pathreplacing}
\usepackage[T1]{fontenc}
\usepackage[utf8]{inputenc}
\usepackage{csquotes}
\usepackage{mathtools,apptools}
\usepackage{pdfpages}
\usepackage[top=.8in, bottom=.8in, left=1.1in, right=1.1in]{geometry}
\usepackage{fancyhdr}
\usepackage{nicefrac}
\usepackage{tikz-cd}
\usepackage{esint,dsfont}
\usepackage[]{hyperref}
\hypersetup{breaklinks=true,hypertexnames=false,colorlinks=true,urlcolor=green!50!blue,citecolor=green!50!blue,linkcolor=green!50!blue,bookmarksnumbered,bookmarksopen}
\usepackage{orcidlink}

\DeclareMathOperator{\pv}{p.v.}
\DeclareFontFamily{U}{mathx}{}
\DeclareFontShape{U}{mathx}{m}{n}{ <-> mathx10 }{}
\DeclareSymbolFont{mathx}{U}{mathx}{m}{n}
\DeclareFontSubstitution{U}{mathx}{m}{n}
\DeclareMathAccent{\widecheck}{0}{mathx}{"71}
\DeclareMathOperator{\supp}{supp}

\theoremstyle{plain}
\newtheorem{theorem}{Theorem}[section]
\newtheorem{Definition}[theorem]{Definition}
\newtheorem{Corollary}[theorem]{Corollary}
\newtheorem{Lemma}[theorem]{Lemma}

\newtheorem{Remark}[theorem]{Remark}

\newtheorem{Proposition}[theorem]{Proposition}

\newcommand{\Proof}{\noindent{\it \textbf{Proof.}\hspace{2mm}}}

\numberwithin{equation}{section}

\newcommand{\Z}{\mathbb Z}
\newcommand{\R}{\mathbb R}
\newcommand{\C}{\mathbb C}
\newcommand{\Rn}{\R^N}
\newcommand{\N}{{\mathbb N}}
\newcommand{\Hs}{\dot{H}^s(\R^N)}
\newcommand{\hs}{\dot{H}^s}
\newcommand{\fhs}{(\dot{H}^s)'}
\newcommand{\Fhs}{(\Hs)'}

\newcommand{\lst}{L^{2_s^{\star}(t)}_{|x|^{-t}}}
\newcommand{\twostar}{2_s^{\star}(t)}
\def\cc{C_c^\infty(\Rn)}
\def\dx{\,{\rm d}x}
\def\dy{\,{\rm d}y}
\def\dz{\,{\rm d}z}
\def\Flap{\left(-\Delta\right)^s}
\def\flap{\left(-\Delta\right)^{\frac s2}}
\def\alap{\left(-\Delta\right)^\alpha}
\newcommand{\id}{\mathrm{Id}}
\stackMath
\newcommand\Widehat[1]{%
\savestack{\tmpbox}{\stretchto{%
  \scaleto{%
    \scalerel*[\widthof{\ensuremath{#1}}]{\kern-.6pt\bigwedge\kern-.6pt}%
    {\rule[-\textheight/2]{1ex}{\textheight}}
  }{\textheight}%
}{0.5ex}}%
\stackon[1pt]{#1}{\tmpbox}%
}
\allowdisplaybreaks

\begin{document}
\begin{abstract}
We establish sharp quantitative multi-bubble stability for nonnegative critical points of the fractional Hardy--Sobolev inequality in the low-dimensional regime $2s<N<6s-2t$. For functions whose energy is close to that of a finite superposition of bubbles, we prove that the Euler--Lagrange deficit controls linearly the distance, in the homogeneous fractional Sobolev norm, to the multi-bubble set, and we recover the precise bubble configuration. This yields quantitative rigidity under arbitrary finite weak interactions. The proof combines a localization scheme adapted to the Hardy weight, Kato--Ponce type weighted commutator estimates, a bubble-wise spectral gap inequality, and a sharp bubble-bubble interaction analysis. We also show that the linear rate is optimal by constructing a matching counterexample.
\end{abstract}

\title[On the sharp multi-bubble stability for fractional Hardy--Sobolev equations: a quantitative approach in low dimensions]{On the sharp multi-bubble stability for fractional Hardy--Sobolev equations: a quantitative approach in low dimensions}
	
\author{Souptik Chakraborty {\tiny and} Utsab Sarkar}
	
\address{\orcidlink{0009-0004-1867-0560}\,\textbf{Souptik Chakraborty} Gandhi Institute of Technology and Management (GITAM), Department of Mathematics and Statistics, Rudraram, Patancheru Mandal, Hyderabad, Telangana, 502329, India. \normalfont\emph{Email address}: \texttt{schakrab8@gitam.edu, soupchak9492@gmail.com}}

\address{\orcidlink{0000-0003-4343-0724}\,\textbf{Utsab Sarkar} Department of Mathematics, Indian Institute of Technology Bombay, Powai, Mumbai 400076, Maharashtra, India. \normalfont\emph{Email address}: \texttt{utsab@math.iitb.ac.in, reachutsab@gmail.com}}
	
	

\maketitle
\tableofcontents

\section{Introduction} Hardy--Sobolev inequalities, both local and fractional, play a central role in the analysis of critical elliptic equations, capturing the sharp interaction between nonlocal diffusion, singular weights, and critical exponents; see~\cite{L83, CKN84, CLO06, FS08, GS15, GRSZ18}. A major theme in recent works have been the development of \emph{stability theories}: understanding how the Euler--Lagrange deficit associated with a critical functional governs the distance to certain set of extremal functions in a quantitative fashion. In the classical Sobolev case, a comprehensive picture has emerged, from the early variational and topological methods of~\cite{BE91,S84,Bahri-Coron-88,AB89} to the \emph{quantitative multi-bubble theory} pioneered by Figalli--Glaudo and subsequent developments~\cite{CFM18,FG20,DSW21,WW22,CKW26}.

In the fractional and weighted setting the landscape is more delicate. The nonlocal operator obstructs standard localization techniques, the Hardy weight modifies the linearized spectrum, and bubble interactions occur across multiple scales. In our earlier work \cite{CS25}, we established sharp quantitative stability for the fractional Hardy--Sobolev inequality and the associated Euler--Lagrange equation in the \emph{single-bubble} regime and identified the structural mechanisms underlying rigidity. The present work is a natural continuation of that program.

Our goal here is to establish the {multi-bubble stability theory} for the fractional Hardy--Sobolev equation~\eqref{maineq} in the \emph{low-dimensional} regime $2s < N < 6s - 2t$. We prove that any function that nearly solves the equation~\eqref{maineq} must lie close, in the natural fractional energy topology, to a finite superposition of weakly-interacting bubbles, and that its distance to the set of \emph{sums of bubbles} is quantitatively controlled by the Euler--Lagrange deficit~\eqref{def-ga}. This characterizes the full structure of near-extremal configurations for the fractional Hardy--Sobolev inequality~\eqref{fhs0} and yields a weighted, genuinely nonlocal analogue of the multi-bubble rigidity phenomena previously available only in the classical or unweighted fractional Sobolev settings. In particular, it extends the quantitative multi-bubble theory for critical Sobolev and Hardy--Sobolev inequalities and their associated Euler--Lagrange equations on $\Rn$ to the fractional framework; see~\cite{CFW13,DNK23JFA,DNK23CalVar,A23,CKW25,CS25}.

\subsection{Background and hurdles}
We begin by recalling the fractional Hardy--Sobolev inequality on $\Rn$,\begin{equation}\label{fhs0}
	\mu_{s,t}(\Rn)
	\left(\int_{\Rn}\frac{|u|^{\twostar}}{|x|^t}\,\mathrm{d}x\right)^{\frac{2}{\twostar}}
	\le
	\int_{\Rn}|(-\Delta)^{\frac{s}{2}}u|^2\,\mathrm{d}x,
	\qquad u\in\dot{H}^s(\Rn).
\end{equation} which provides the variational framework for our study in the homogeneous Sobolev space $\dot{H}^s(\Rn)$ of order $s$. Here $s\in(0,1)$, $t\in(0,2s)$, and $N>2s$ is a natural number, and the critical Hardy--Sobolev exponent is
\[\twostar\coloneqq\tfrac{2(N-t)}{N-2s},\qquad 2<\twostar<\tfrac{2N}{N-2s}.\] 

Throughout, we work in the same setting as in~\cite{CS25}. For the precise definitions of the fractional Sobolev spaces, the associated fractional Laplacian, and the properties used in this paper, we refer the reader to~\cite[\S~2.1]{CS25}.

The optimal constant $\mu_{s,t}>0$ is characterized variationally by
\begin{equation*}
	\mu_{s,t}(\Rn)
	=
	\inf_{u\in\Hs\setminus \{0\}}
	\frac{\|u\|_{\Hs}^{2}}
	{\big\||x|^{\frac{-t}{\twostar}}u\big\|_{\twostar}^{2}},
\end{equation*} and is attained by a radially symmetric, nonnegative, strictly decreasing minimizer $U_{s,t}$ (see e.g.~\cite{MRNA21}).  

All minimizers are obtained from $U_{s,t}$ through the natural scaling
\begin{equation}\label{dilation map}
	cU^\lambda_{s,t}(x)\coloneqq c\lambda^{\frac{N-2s}{2}}U_{s,t}(\lambda x), \qquad c\neq0,\, \lambda>0.
\end{equation}

Under the normalization $\||x|^{\frac{-t}{\twostar}}U_{s,t}\|_{\twostar}=1$, we define $V_{s,t}\coloneqq\mu_{s,t}^{\frac{1}{\twostar-2}}U_{s,t}$, so that the rescaled bubbles $V_{s,t}^\lambda$ satisfy \begin{equation}\label{bubble normalization}
	\|V^\lambda_{s,t}\|_{\hs}
	=
	\mu_{s,t}^{\frac{N-t}{2(2s-t)}},
	\quad
	\||x|^{\frac{-t}{\twostar}}V^\lambda_{s,t}\|_{\twostar}
	=
	\mu_{s,t}^{\frac{1}{\twostar-2}},\qquad \forall\lambda>0.
\end{equation}

The function $V_{s,t}$ is a weak solution of the normalized Euler--Lagrange
equation
\begin{equation}\label{maineq}
\begin{cases}
	(-\Delta)^s u = \dfrac{|u|^{\twostar-2}u}{|x|^t},
	\quad u\in\dot{H}^s(\mathbb{R}^N),\\[1mm]
	u>0\ \text{in }\Rn,
\end{cases}
\end{equation}
enjoying the uniqueness, growth and regularity properties described in~\cite[Theorem~1.1]{MN21}. Moreover, any other solution to \eqref{maineq} is given by $V_{s,t}^\lambda$ for some $\lambda>0$.

\begin{Remark}\label{remark on eigenfunctions}
We collect two structural properties of the bubbles $V=V^{\lambda}_{s,t}$ that will be used throughout the analysis.

\begin{enumerate}
\item[(1)] By \cite[Lemma~4.2]{MN21} and \cite[Theorem~6.2]{CS25}, the operator 
\[
\frac{\Flap}{|x|^{-t}V^{p-1}};\quad p+1=\twostar,
\]
has a discrete spectrum whose first two eigenpairs are $(1,V)$ and $(p,\dot V)$.  
The eigenvalues are non-decreasing, diverge to $\infty$, and do not depend on the scale $\lambda>0$.  
Moreover, \cite[Theorem~1.2]{MN21} shows that the third eigenvalue is strictly larger than $p$.

\item[(2)] The pointwise bound $|\dot V(x)|\lesssim_{N,s,t}\frac{V(x)}{\lambda}$ holds on $\Rn$ by \cite[Eq.~19, p.~301 \& Theorem~1.1]{MN21}.  
Using this together with \eqref{bubble normalization} and \eqref{fhs0}, we obtain
\[
\|\dot V\|_{L^{p+1}_{|x|^{-t}}}
   \lesssim_{N,s,t}\|\dot V\|_{\Hs}
   \lesssim_{N,s,t}\frac{\sqrt p}{\lambda}\|V\|_{L^{p+1}_{|x|^{-t}}}^{\frac{p+1}{2}}
   \lesssim_{N,s,t,}\frac{1}{\lambda}.
\]
\end{enumerate}
\end{Remark}

To measure the deviation of a function from satisfying~\eqref{maineq}, we introduce the
Euler--Lagrange deficit
\begin{equation}\label{def-ga}
	\Gamma(u)
	\coloneqq
	\Bigg\|(-\Delta)^s u -
	{|\cdot|^{-t}|u|^{\twostar-2}u}\Bigg\|_{\Fhs}.
\end{equation}

This deficit $\Gamma$ plays the same structural role as in the classical stability theory for
critical equations, and aligns with the concentration-compactness principle of Lions~\cite{PLL184, PLL284} and
the multi-bubble analysis initiated by Struwe and successors~\cite{S84, CGS89, R90, G98, S08, PP14}.

\medskip
Our earlier work~\cite[Theorem~1.4]{CS25} established a Struwe-type \emph{qualitative} stability result for~\eqref{maineq}. In particular, nonnegative sequences with nearly optimal Hardy--Sobolev energy and vanishing deficit admit a decomposition into finitely many well-separated bubbles, mirroring the classical bubbling patterns that arise in critical elliptic problems. For the reader's convenience, we recall the statement below.

\begin{theorem}[Qualitative stability]\label{SDHS}
Let $s\in(0,1)$, $t\in(0,2s)$, $N>2s$, and $N,\nu\in\N$. Let $\{u_k\}_{k\ge1}\subset\Hs$ be a sequence of nonnegative functions such that \[
	\Bigl(\nu-\tfrac12\Bigr)\mu_{s,t}^{\frac{N-t}{2s-t}}
	\le \|u_k\|_{\Hs}^2
	\le
	\Bigl(\nu+\tfrac12\Bigr)\mu_{s,t}^{\frac{N-t}{2s-t}},
	\qquad
	\Gamma(u_k)\to 0\ \text{as }k\to\infty.
\]

Then, up to a subsequence, there exist positive scales
$\lambda_k^1,\ldots,\lambda_k^\nu$ such that
\[
	\Big\|u_k-\sum_{i=1}^{\nu}V_{s,t}^{\lambda_k^i}\Big\|_{\Hs}\to0.
\]

If $u_k\rightharpoonup u_0\neq0$, then $\lambda_k^1=\lambda_0$ for all $k$; whereas if $u_k\rightharpoonup0$, no scale remains fixed. Moreover,
\[
	\lambda_k^i\to0\ \text{or}\ \lambda_k^i\to\infty,
	\qquad
	\big|\log(\lambda_k^i/\lambda_k^j)\big|\to\infty\quad (i\neq j),
\] except in the fixed-scale case.
\end{theorem}

The purpose of the present work is to develop the \emph{quantitative} analogue of Theorem~\ref{SDHS}. Restricting to the low-dimensional range $2s<N<6s-2t$, we prove that the deficit $\Gamma(u)$ \emph{linearly} controls the $\dot{H}^s$-distance from the multi-bubble set
\[
	\mathcal{M}
	=
	\Bigg\{\sum_{i=1}^{\nu}V_{s,t}^{\lambda_i}:\lambda_i>0\Bigg\}.
\]

This extends the one-bubble stability theory of~\cite{CS25} to arbitrary $\nu\ge1$ and
constitutes the first quantitative stability result (Theorem~\ref{energybound-implies-stability}) for the fractional Hardy--Sobolev
equation~\eqref{maineq} in the multi-bubble setting.

A fundamental difficulty in the fractional Hardy--Sobolev setting is that localization interacts poorly with both the nonlocal operator and the Hardy weight: any naive cutoff produces long-range errors, and classical commutator bounds do not capture the weighted structure of the problem. To overcome these obstacles, our analysis relies on two principal tools developed in Appendix~\ref{AA}. The first is a sharp cutoff estimate (Proposition~\ref{basic-cutoff-norm-estimate}) for logarithmic bump functions, which provides precise control of the fractional Laplacian on transition annuli. The second is a weighted Kato--Ponce type commutator estimate (Theorem~\ref{wkpv-estimate}), tailored to accommodate the Hardy singularity. Together, these results deliver robust bounds on localization errors and on mixed fractional terms arising from the interaction between cutoffs and the operator $(-\Delta)^s$. When coupled with our refined bump-function projection scheme (Lemma~\ref{bumps}), the weighted spectral gap inequality (Proposition~\ref{spectral-gap-inequality}), and finally the robustness of bubble coefficients (Proposition~\ref{bub-bub-man}), these analytic inputs yield the quantitative stability result (Theorem~\ref{quantitative-stability}) for families of weakly-interacting bubbles (cf.~\eqref{del-int-bub}). This stability principle then extends (\emph{via} Theorem~\ref{SDHS}) from the localized, projection-based framework to the full multi-bubble set, thereby completing the low-dimensional theory of quantitative stability for the fractional Hardy--Sobolev equation, as formulated in Theorem~\ref{energybound-implies-stability}. A further ingredient in our analysis is a global nearest-point projection onto the full multi-bubble set $\mathcal M_\nu$. This requires more than the local modulation arguments used in the stability literature: we establish a global parameter-matching lemma that confines arbitrary near-minimizers to a common local chart, where compactness and strict convexity yield attainment and uniqueness. In particular, Appendix~\ref{AC} provides the foundation underlying Theorem~\ref{quantitative-stability} by establishing rigorous projection framework, together with the associated tangent-space orthogonality and quantitative parameter stability, in the fractional Hardy--Sobolev setting.

\medskip
In a forthcoming companion work, we will address the multi-bubble scenario in the \emph{high-dimensional} regime $N \ge 6s - 2t$, where additional degeneracies appear and new techniques are required. Taken together, these results aim to place the fractional Hardy--Sobolev stability theory on the same level of precision as its classical counterpart while introducing analytic tools, in particular our weighted commutator estimate, that may prove useful across a wider class of nonlocal PDEs with singular potentials.

\subsection{Main results}
This section contains the principal quantitative stability results of the paper. We begin with a stability theorem for families of weakly (\emph{or} $\delta$)-interacting bubbles, giving us uniform control of perturbations orthogonal to the natural modulation directions under a smallness assumption in $\Hs$.

\begin{theorem}[Stability for $\delta$-interacting bubbles]\label{quantitative-stability}
For $s\in(0,1)$, $t\in(0,2s)$, and integers $N \in (2s,\,6s-2t)$ and $\nu\ge1$, there exist constants $\delta(N,\nu,s,t)>0$ and $C(N,\nu,s,t)>0$ such that the following holds. Suppose that $u\in\Hs$ satisfies
\begin{equation*}
    \Big\|u-\sum_{i=1}^{\nu} \widetilde{V}_i\Big\|_{\Hs} \le \delta,
\end{equation*}
where $(\widetilde{V}_i)_{i=1}^\nu$ is a $\delta$-interacting family of bubbles. Then there exists a $\delta'$-interacting family of bubbles $(V_i)_{i=1}^\nu$ such that
\begin{equation*}
    \bigg\|u-\sum_{i=1}^{\nu} V_i\bigg\|_{\Hs}
    \le C(N,\nu,s,t) \,\bigg\|\Flap u-\frac{|u|^{\twostar-2}u}{|x|^t}\bigg\|_{\Fhs},
\end{equation*}
where $\delta'\to 0$ as $\delta\to0$.
\end{theorem}

Combining Theorem~\ref{quantitative-stability} with the qualitative Struwe-decomposition of Theorem~\ref{SDHS} yields the global sharp statement, which forms the main theorem below.
    
\begin{theorem}[Sharp quantitative stability]\label{energybound-implies-stability}
Let $s\in(0,1)$, $t\in(0,2s)$, and choose integers $N\in (2s,\,6s-2t),\,\nu\ge1$. Let $u\in \Hs$ be a non-negative function satisfying
\begin{equation}\label{energy-bdd}
\left(\nu-\tfrac12\right)\mu_{s,t}^{\frac{N-t}{2s-t}}
\le\|u\|_{\Hs}^{2}\le\left(\nu+\tfrac12\right)\mu_{s,t}^{\frac{N-t}{2s-t}}.
\end{equation}

Then there exist a constant $C(N,\nu,s,t)>0$ and parameters
$\lambda_{1},\dots,\lambda_{\nu}>0$ such that
\[
    u=\sum_{i=1}^{\nu} V_{s,t}^{\lambda_{i}} + \rho,
\]
and the error term $\rho$ satisfies
\begin{equation}\label{def-bds-dis-lin}
    \|\rho\|_{\Hs} \le C(N,\nu,s,t)\, \Gamma(u).
\end{equation}

Moreover, this estimate is optimal in the sense that the conclusion fails if the term
$\|\rho\|_{\Hs}$ is replaced by $\|\rho\|_{\Hs}^{1-\varepsilon}$ for any $\varepsilon\in(0,1)$.

Further, for all $i\neq j$, the $\hs$-interaction between the bubbles can be estimated as \begin{equation}\label{Hs-interaction}\langle V_{s,t}^{\lambda_i},V_{s,t}^{\lambda_j}\rangle_{\Hs}=\int_{\Rn}\frac{V_{s,t}^{\lambda_i}(x)V_{s,t}^{\lambda_j}(x)^{\twostar-1}}{|x|^t}\dx\lesssim_{N,\nu,s,t}\Gamma(u).\end{equation}
\end{theorem}

    \subsection{Notations} The following notations will be used consistently throughout the paper, except where a different convention is explicitly specified.
	
	\begin{enumerate}
	\item The Schwartz space of rapidly decreasing functions is
	\[
	\mathcal{S}(\Rn, \C) \coloneqq \left\{ f \in C^\infty(\Rn, \C) \,:\, \forall \boldsymbol{\alpha}, \boldsymbol{\beta} \in \N_0^N,\, \|f\|_{\boldsymbol{\alpha},\boldsymbol{\beta}} < \infty \right\},
	\]
	where the semi-norm is given by
	$
	\|f\|_{\boldsymbol{\alpha},\boldsymbol{\beta}} \coloneqq \sup_{x \in \Rn} \big| x^{\boldsymbol{\alpha}} \partial^{\boldsymbol{\beta}} f(x) \big|
	$
	with $x^{\boldsymbol{\alpha}} = x_1^{\alpha_1} \cdots x_N^{\alpha_N}$ and $\partial^{\boldsymbol{\beta}} = \partial_{x_1}^{\beta_1} \cdots \partial_{x_N}^{\beta_N}$.\vspace{1mm}

	\item For $u \in \mathcal{S}(\Rn)$, the Fourier transform is defined by
	\[
	\mathcal{F}(u)(\xi) \coloneqq \Widehat{u}(\xi) \coloneqq (2\pi)^{-\frac{N}{2}} \int_{\Rn} e^{-i x \cdot \xi} u(x)\dx, \quad \xi \in \Rn.
	\]

	\item $\|\cdot\|_p$ denotes the standard $L^p(\Rn)$-norm. Given a positive weight $w \in L^1_{\rm loc}(\Rn)$, the weighted Lebesgue space is
	\[
	L^p_w(\Rn) \coloneqq \left\{ f : \Rn \to \C \,{\Big|}\, f \text{ is measurable and } \int_{\Rn} |f|^p w \dx < \infty \right\},
	\]
	with norm $\|f\|_{L^p_w}^p = \|w^{1/p} f\|_p^p =\int_{\Rn} |f|^p w \dx$.\vspace{1mm}

	\item Denote by $\Fhs$ the dual of the fractional homogeneous Sobolev space $\Hs$, endowed with the dual norm induced by $\|\cdot\|_{\hs}$. The duality pairing is denoted by ${}_{\fhs}\langle f, u\rangle_{\hs}\text{ for }f \in\Fhs, u \in \Hs.$ The Hilbert space inner product on $\Hs$ is denoted as $\langle \cdot, \cdot \rangle_{\hs}$ and is related to the usual $L^2$-inner product by $\langle u,v\rangle_{\hs} = \langle \flap u,\flap v\rangle=\int_{\Rn}\flap u(x)\overline{\flap v(x)}\dx$.\vspace{1mm}

	\item For nonnegative quantities $A,B$; we write \( A \lesssim B \) to mean \( A \le C B \) for some constant \( C > 0 \), and \( A \gtrsim B \) if \( B \lesssim A \). We also write \( A \sim B \) if both \( A \lesssim B \) and \( A \gtrsim B \). Dependence of the implicit constant on parameters \( a, b, \dots \) is indicated by \( \lesssim_{a,b,\dots} \), and similarly for \( \gtrsim \) and \( \sim \).\vspace{1mm}

    \item We write $C$ for a positive constant whose value may vary from line to line. When dependence on parameters $a,b,\dots$ is relevant, we use $C(a,b,\dots)$. In contrast, symbols of the form $C_{a,b,\dots}$ denote fixed, explicitly known constants from the literature.
\end{enumerate}

\subsection{Structure of the paper}
\S\ref{S2} develops the quantitative stability theory in the low-dimensional regime $2s < N < 6s-2t$. We first establish Theorem~\ref{quantitative-stability}, which provides linear control of the error for weakly-interacting bubbles. The localization scheme is introduced in \S\ref{S2.2}, and the spectral gap input is obtained in \S\ref{S2.3}. Interaction estimates and the robustness of bubble coefficients are derived in \S\ref{S2.4}, leading to the proof of Theorem~\ref{energybound-implies-stability} in \S\ref{S2.5}. Appendix~\ref{AA} collects the weighted fractional tools and commutator estimates used in the analysis, while Appendix~\ref{AB} records the two-bubble integrals required for the interaction bounds. Finally, Appendix~\ref{AC} establishes the global nearest-point projection onto $\mathcal M_\nu$, proving global parameter matching, attainment, uniqueness and the resulting tangent-space orthogonality.

\section{Stability in low dimensions \texorpdfstring{$2s<N<6s-2t$}{}}\label{S2}
\begin{Definition}[$\delta$-interacting bubbles]\label{del-int-bub}
	Let $V_i\coloneqq V_{s,t}^{\lambda_i}$ with scaling magnitude $\lambda_i>0$; $1\le i\le\nu$, be a family of bubbles. We say that the family is $\delta$-interacting for some $\delta>0$ if for all distinct $1\le i,\,j\le\nu$,
		\begin{equation}\label{delta interaction}
			\min\Bigg\{\frac{\lambda_i}{\lambda_j}, \frac{\lambda_j}{\lambda_i}\Bigg\} \le \delta.\end{equation}

        If we consider for some $\alpha_1,\dots,\alpha_\nu\in\R$, the family $\left(\alpha_iV_i\right)_{i=1}^{\nu}$, then the family is said to be $\delta$-interacting if in addition to \eqref{delta interaction} \begin{equation*}\max_{1\le i\le \nu}|\alpha_i-1|\le \delta\end{equation*} also holds. Moreover, for $\nu=1$, the $\delta$-interaction condition is understood vacuously, so that every single-bubble family is $\delta$-interacting for every $\delta>0$.
	\end{Definition}

    \begin{Remark}\label{bub-del}
        The notion of $\delta$-interaction between bubbles is naturally encoded by their interaction in the $\hs$ inner product.  For two bubbles $U_1=V_{s,t}^{\lambda_1}$ and $U_2=V_{s,t}^{\lambda_2}$, the corresponding $\hs$-interaction can be quantified explicitly.  Since both $U_1$ and $U_2$ solve \eqref{maineq}, one has \begin{equation*}
\langle U_1,U_2\rangle_{\hs}
 = \int_{\Rn}\frac{U_1^{\twostar-1}U_2}{|x|^t}\dx
 = \int_{\Rn}\frac{U_1 U_2^{\twostar-1}}{|x|^t}\dx
 \sim_{N,s,t}
 \min\!\left\{\frac{\lambda_1}{\lambda_2},\frac{\lambda_2}{\lambda_1}\right\}^{\frac{N-2s}{2}},
\end{equation*} as a consequence of Corollary~\ref{general-bubble-interaction-estimate}.  

In particular, membership in a
$\delta$-interacting family implies that the $\hs$ inner product between distinct bubbles is of order $\mathcal{O}(\delta^{\frac{N-2s}{2}})$.
\end{Remark}

    \subsection{Proof of Theorem~\ref{quantitative-stability}}
Fix a small parameter $0<\delta\ll1$, to be chosen later. We consider the set of superpositions of $\nu$ bubbles \[ \mathcal M_\nu \coloneqq \left\{ \sum_{i=1}^{\nu}\alpha_i V_{s,t}^{\lambda_i} \;\middle|\; \alpha_i\in\R,\ \lambda_i\in(0,\infty) \right\} \subset \Hs. \] 

Given $u\in \Hs$, we study the minimization problem \begin{equation}\label{eq:min-problem} \inf_{\sigma\in\mathcal M_\nu}\|u-\sigma\|_{\Hs}. \end{equation} 

For $\epsilon>0$, define the tubular neighborhood \begin{equation}\label{def:tubular-nhbd} U(\epsilon) \coloneqq \left\{ v\in\Hs \;\middle|\; \begin{array}{l} \text{there exists a $2\epsilon$-interacting family } (\alpha_i V_{s,t}^{\lambda_i})_{i=1}^\nu \\[2pt] \text{such that } \displaystyle \bigg\|v-\sum_{i=1}^\nu \alpha_i V_{s,t}^{\lambda_i}\bigg\|_{\Hs} <\epsilon \end{array} \right\}.\end{equation}

By Proposition~\ref{prop:AC-nearest-point}, there exists $\epsilon_0>0$, depending only on $N,\nu,s,t$, such that for every $0<\epsilon\le\epsilon_0$ and every $u\in U(\epsilon)$, the minimization problem~\eqref{eq:min-problem} admits a unique minimizing function $\sigma\in\mathcal M_\nu$, where uniqueness of the bubble parameters is understood up to permutation of the indices. Writing $u=\sigma+\rho$, the remainder satisfies~\eqref{AC-remainder-bound}. Moreover, for any reference family supplied by $u\in U(\epsilon)$, the minimizing parameters satisfy~\eqref{AC-parameter-stability}, and the minimizing family is $C\epsilon$-interacting, by \eqref{AC-interaction-stability}.

We denote this unique minimizer by
\[
\sigma=\sum_{i=1}^{\nu}\alpha_i V_i,
\qquad
V_i=V_{s,t}^{\lambda_i}.
\]

By the hypothesis of Theorem~\ref{quantitative-stability}, $\|u-\sum_{i=1}^{\nu}\widetilde V_i\|_{\Hs}\le\delta$, for some $\delta$-interacting family $(\widetilde V_i)_{i=1}^{\nu}$, where $\widetilde V_i=V_{s,t}^{\widetilde\lambda_i}$. Taking $\epsilon=2\delta$ in the definition~\ref{def:tubular-nhbd} of $U(\epsilon)$ implies $u\in U(2\delta)$. Proposition~\ref{prop:AC-nearest-point} then gives the unique minimizer $\sigma$. Since $\sum_i\widetilde V_i\in\mathcal M_\nu$, its minimizing property yields the estimate
\[
\|\rho\|_{\Hs}=\|u-\sigma\|_{\Hs}
=\inf_{\tau\in\mathcal M_\nu}\|u-\tau\|_{\Hs}
\le\big\|u-\sum_{i=1}^{\nu}\widetilde V_i\big\|_{\Hs}
\le\delta.
\]

Furthermore, applying~\eqref{AC-interaction-stability} with $\epsilon=2\delta$ and with $\widetilde\alpha_i=1$ gives, after a permutation of the indices,
\[
\max_{1\le i\le\nu}|\alpha_i-1|
\le 2C\delta,
\qquad
\min\left\{\frac{\lambda_i}{\lambda_j},
\frac{\lambda_j}{\lambda_i}\right\}
\le 2C\delta,
\qquad i\neq j.
\]

Consequently, the parameter $\delta'$ in
Theorem~\ref{quantitative-stability} may be taken to be $2C\delta$.

Since $V_i$ and $\dot{V}_i$ are the first two eigenfunctions of the operator $\frac{(-\Delta)^s}{|x|^{-t}V_i^{p-1}}$ corresponding to the eigenvalues $1$ and $p=\twostar-1$, respectively (cf.~Remark~\ref{remark on eigenfunctions}-$(1)$), it follows from~\eqref{AC-tangent-orthogonality} and integration by parts that the following orthogonality relations hold for all $1\le i\le\nu$.
\begin{align}
0=\langle\rho, V_i\rangle_{\dot{H}^s}
   &= \int_{\Rn}\frac{V_i^{p}\rho}{|x|^t}\dx, \label{Inn-bub}\\
0=\langle\rho, \dot{V}_i\rangle_{\dot{H}^s}
   &= p\int_{\Rn}\frac{V_i^{p-1}\dot{V}_i\rho}{|x|^t}\dx. \label{Inn-der-bub}
\end{align}

   Our objective is to bound $\|\rho\|_{\hs}$ in terms of the deficit $\Gamma(u)$. To this end, we exploit the orthogonality condition \eqref{Inn-bub}. Using
integration by parts, we compute
\begin{align}
\left\|\rho\right\|_{\hs}^2
&= \langle \rho,\rho\rangle_{\dot{H}^s}
 = \langle u-\sigma,\rho\rangle_{\dot{H}^s}
 = \left\langle (-\Delta)^{\frac{s}{2}}u, (-\Delta)^{\frac{s}{2}}\rho\right\rangle
 \nonumber\\
&={}_{\fhs}\bigg\langle (-\Delta)^{s} u-\frac{u|u|^{p-1}}{|x|^t},\rho\bigg\rangle_{\hs}
 +{}_{\fhs}\bigg\langle \frac{u|u|^{p-1}}{|x|^t}, \rho\bigg\rangle_{\hs}
 \nonumber\\
&\le \Gamma (u)\|\rho\|_{\hs}
 + \int_{\Rn}\frac{u|u|^{p-1}\rho}{|x|^t}\dx.
 \label{estimate-1}
\end{align}

   To estimate the final term in \eqref{estimate-1}, we invoke the following two
elementary inequalities, which will be used repeatedly in the sequel.
\begin{itemize}
\item For any $a,b\in\R$ and $p\ge 1$, there exists a constant
$C(p)>0$ such that
\begin{equation}\label{Ineq-1}
\left|(a+b)|a+b|^{p-1}-a|a|^{p-1}\right|
\le p|a|^{p-1}|b| + C(p)\Big(\chi_{p>2}|a|^{p-2}|b|^2 + |b|^p\Big),
\end{equation}
where $\chi_{p>2}=1$ if $p>2$ and $\chi_{p>2}=0$ for $1\le p\le 2$.\vspace{2mm}

\item For any $a_1,\dots,a_{\nu}\in\R$ and $p\ge 1$, one has
\begin{equation}\label{Ineq-2}
\Bigg|\Big(\sum_{i=1}^{\nu}a_i\Big)\Big|\sum_{i=1}^{\nu}a_i\Big|^{p-1}
-\sum_{i=1}^{\nu}a_i|a_i|^{p-1}\Bigg|
\lesssim_{p,\nu}
\sum_{1\le i\neq j\le \nu}|a_i|^{p-1}|a_j|.
\end{equation}
\end{itemize}
    
  Applying \eqref{Ineq-1} and \eqref{Ineq-2} with $a=\sigma$, $b=\rho$,
$a_i=\alpha_iV_i$, and $p=\twostar-1$, and noting that $N\in(2s,6s-2t)$ implies
$p>2$, we obtain the following pointwise estimate:
\begin{equation}\label{estimation-2}
\Bigg|\frac{u|u|^{p-1}}{|x|^t}
-\sum_{i=1}^{\nu}\alpha_i|\alpha_i|^{p-1}\frac{V_i^{p}}{|x|^t}\Bigg|
\le
p\frac{\sigma^{p-1}|\rho|}{|x|^t}
+ C(p,\nu)\bigg(
\frac{\sigma^{p-2}|\rho|^2}{|x|^t}
+ \frac{|\rho|^{p}}{|x|^t}
+ \sum_{1\le i\neq j\le \nu}\frac{V_i^{p-1}V_j}{|x|^t}
\bigg).
\end{equation}

    Therefore, using the orthogonality relations \eqref{Inn-bub}--\eqref{Inn-der-bub} in \eqref{estimation-2}, we estimate the last term of \eqref{estimate-1} as,
    \begin{align}
        \left|\int_{\Rn}\frac{u|u|^{p-1}\rho}{|x|^t}\dx\right| &\le p\int_{\Rn}\frac{\sigma^{p-1}|\rho|^2}{|x|^t}\dx + C(N,\nu,s,t)\Bigg(\int_{\Rn}\frac{\sigma^{p-2}|\rho|^3}{|x|^t}\dx\nonumber\\
        &\qquad+\int_{\Rn}\frac{|\rho|^{p+1}}{|x|^t}\dx +\sum_{1\le i\neq j\le \nu} \int_{\Rn}\frac{V_i^{p-1}V_j|\rho|}{|x|^t}\dx\Bigg).\label{estimation-3}
    \end{align}
    
    Now employing H\"older's inequality and fractional Hardy--Sobolev inequality~\eqref{fhs0}, we get,
    \begin{align}
    &\,\int_{\Rn}\frac{\sigma^{p-2}|\rho|^3}{|x|^t}\dx \le  \left(\int_{\Rn}\frac{\sigma^{p+1}}{|x|^t}\dx\right)^{\frac{p-2}{p+1}} \left(\int_{\Rn}\frac{|\rho|^{p+1}}{|x|^t}\dx\right)^{\frac{3}{p+1}}\lesssim_{N,\nu,s,t}\|\rho\|_{\hs}^3.\label{esti-4}\\
    &\,\int_{\Rn} \frac{V_i^{p-1}V_j|\rho|}{|x|^t}\dx \le \|\rho\|_{L^{p+1}_{|x|^{-t}}}\|V_i^{p-1}V_j\|_{L^{\frac{p+1}{p}}_{|x|^{-t}}} \lesssim_{N,s,t}\|\rho\|_{\hs} \|V_i^{p-1}V_j\|_{L^{\frac{p+1}{p}}_{|x|^{-t}}}\label{esti-5}.
    \end{align}
    
   Observe that the dimensional assumption $N<6s-2t$, equivalently $p>2$, guarantees
the validity of \eqref{esti-4}. Substituting \eqref{esti-4},\,\eqref{esti-5} into
\eqref{estimation-3}, we obtain
\begin{align}
\left|\int_{\Rn}\frac{u|u|^{p-1}\rho}{|x|^t}\dx\right|
&\le
p\int_{\Rn}\frac{\sigma^{p-1}\rho^2}{|x|^t}\dx
+C\Bigg(
\|\rho\|_{\hs}^3
+ \|\rho\|_{\hs}^{p+1}+\|\rho\|_{\hs}\sum_{i\neq j}\|V_i^{p-1}V_j\|_{L^{\frac{p+1}{p}}_{|x|^{-t}}}
\Bigg).
\label{estimation-4}
\end{align}

   For $i\neq j$, Corollary~\ref{general-bubble-interaction-estimate} with the notation $Q_{ij}=\min\left\{\tfrac{\lambda_i}{\lambda_j},\tfrac{\lambda_j}{\lambda_i}\right\}$, yields
\begin{equation}\label{bb}
\|V_i^{p-1}V_j\|_{L^{\frac{p+1}{p}}_{|x|^{-t}}}
=
\left(\int_{\Rn}\frac{V_i^{\frac{(p-1)(p+1)}{p}}V_j^{\frac{p+1}{p}}}{|x|^t}\dx\right)^{\frac{p}{p+1}}
\sim_{N,s,t}
Q_{ij}^{\frac{N-2s}{2}}
\sim_{N,s,t}
\int_{\Rn}\frac{V_i^p V_j}{|x|^t}\dx.
\end{equation}

Substituting \eqref{bb} into \eqref{estimation-4}, we obtain
\begin{align}
\left|\int_{\Rn}\frac{u|u|^{p-1}\rho}{|x|^t}\dx\right|
&\le
p\int_{\Rn}\frac{\sigma^{p-1}\rho^2}{|x|^t}\dx
+C\Bigg(
\|\rho\|_{\hs}^3+\|\rho\|_{\hs}^{p+1}+\|\rho\|_{\hs}\sum_{i\neq j}\int_{\Rn}\frac{V_i^pV_j}{|x|^t}\dx
\Bigg),
\label{estimation-5}
\end{align}

   To control the leading term in \eqref{estimation-5}, we invoke
Proposition~\ref{spectral-gap-inequality}, which ensures the existence of a positive
constant $C_0=C(N,\nu,s,t)<1$.
The remaining terms are estimated using \eqref{Est-5}.
Combining these bounds, we conclude that
\begin{align}\label{final0}
\left|\int_{\Rn}\frac{u|u|^{p-1}\rho}{|x|^t}\dx\right|
&\le
\left(C_0+\mathcal{o}(1)\right)\|\rho\|_{\hs}^2
+ C\left(
\|\rho\|_{\hs}^3
+\|\rho\|_{\hs}^{p+1}
+\|\rho\|_{\hs}\Gamma (u)
\right).
\end{align}

   Substituting \eqref{final0} into \eqref{estimate-1}, we arrive at
\begin{align}\label{final}
\left(1-C_0-\mathcal{o}(1)\right)\|\rho\|_{\hs}^2
\lesssim_{N,\nu,s,t}
\|\rho\|_{\hs}\Gamma (u)
+\|\rho\|_{\hs}^3
+\|\rho\|_{\hs}^{p+1}.
\end{align}

  Observe that the $\mathcal{o}(1)$ term in \eqref{final} originates from
Proposition~\ref{bub-bub-man}. Hence, by choosing $\varepsilon>0$ sufficiently
small, we can ensure that $C_0+\mathcal{o}(1)<1$, so that the left-hand side of
\eqref{final} is strictly positive. Since moreover $\|\rho\|_{\hs}<\delta\ll1$,
it follows that
\begin{equation}\label{rho stab}
\|\rho\|_{\hs} \lesssim_{N,\nu,s,t} \Gamma (u).
\end{equation}

 By a standard dichotomy argument, we may restrict attention to the regime $\Gamma(u)<1$. In this case, combining
\eqref{rho stab} with \eqref{Est-4}, we obtain
\begin{align*}
\Big\|u-\sum_{i}V_i\Big\|_{\hs}
&\le
\Big\|u-\sum_{i}\alpha_iV_i\Big\|_{\hs}
+\sum_{i}|\alpha_i-1|\Big\|V_i\Big\|_{\hs} \\
&\lesssim
(1+\mathcal{o}(1))\Gamma(u)
+\Gamma(u)^2
\lesssim
\Gamma(u),
\end{align*}
which concludes the proof.

The choice of the parameter $\delta$ for which
Theorem~\ref{quantitative-stability} holds is obtained by taking the minimum
among all admissible values of $\delta$ arising from
Proposition~\ref{spectral-gap-inequality}, Proposition~\ref{bub-bub-man} and the proof above.\qed

	\subsection{Localization via bump functions}\label{S2.2}
	\begin{Lemma}\label{bumps}
		Let $N>2s$ be the spatial dimension, $\nu\ge2$ be the number of bubbles, and $\varepsilon > 0$ be a fixed parameter. Assume also that the parameters $t$ and $s$ satisfy $0 < t < 2s < 2$. Then there exists a constant $\delta = \delta(N, \nu, \varepsilon, s, t) > 0$ such that the following holds.
		
		Given a $\delta$-interacting family $(V_i)_{i=1}^\nu$, where each $V_i = V_{s,t}^{\lambda_i}$ for some $\lambda_i > 0$, there exists an associated family of bump functions $(\Phi_i)_{i=1}^\nu$, with $\Phi_i : \Rn \rightarrow [0,1]$, satisfying the following properties for every $i$.\vspace{1mm}
		
		\begin{enumerate}
			\item $\int_{\left\{\Phi_i<1\right\}}\frac{V_i^{\twostar}}{|x|^t}\dx\le\varepsilon\mu_{s,t}^{\frac{N-t}{2s-t}}$.\vspace{1mm}
			\item $V_i>\frac{V_j}{\varepsilon}$ for all $j\neq i$ on the set $\{\Phi_i>0\}$.\vspace{1mm}
			\item $\forall j\neq i$ such that $\lambda_j\le\lambda_i$, one has $\sup_{\left\{\Phi_i>0\right\}}V_j\le(1+\varepsilon)\inf_{\left\{\Phi_i>0\right\}}V_j$.\vspace{1mm}
			\item $\left\||x|^\frac{t}{\twostar}\flap\Phi_i\right\|_{\frac{2\cdot\twostar}{\twostar-2}}<\varepsilon$.
		\end{enumerate}
	\end{Lemma}

\Proof Properties (1)--(3) follow directly from \cite[Lemma~3.9]{FG20} once specialized to the present (bubbles with fixed centers) setting. In the following, we restrict our attention to establishing property (4), for which Proposition~\ref{basic-cutoff-norm-estimate} serves as a substitute for \cite[Lemma~3.8]{FG20}.

           By the cutoff construction in \cite[Lemma~3.9]{FG20}, each bump function $\Phi_i$
           is given as a finite product of at most $\nu$ radial cutoffs. More precisely, for each $i$, we introduce the index set $I_i \coloneqq \{\,k\in\{1,\dots,\nu\}\,:\,\lambda_k>\lambda_i\,\}$. With this notation, the cutoff function $\Phi_i$ can be written as
\begin{equation}\label{def cutoff main}
\Phi_i = f_i \prod_{k\in I_i} f_k,
\end{equation}
where $f_i=\varphi_{r_i,R_i}$ and, for $k\in I_i$, $f_k=1-\varphi_{r_k,R_k}$, for
suitable radii $0<r_k<R_k$; see~\eqref{cut-off} for the definition of $\varphi_{r,R}$.
 
            Denote $q=\tfrac{2\cdot\twostar}{\twostar-2}$ and set
$w=|x|^{\frac{tq}{\twostar}}$. First observe that $w\in\mathcal{A}_q$, since
$N>2s>t$. Moreover, the functions $f_k$ satisfy $0\le f_k\le 1$.

We then invoke a standard weighted fractional Leibniz estimate, namely the endpoint $L^\infty$ Kato--Ponce inequality from~\cite[Theorem~1.1]{DG24} (cf.~\cite[Theorem~1.4]{GK96}), applied to the product of the functions $f_k$. This yields the bound
            \[
            \|\flap(f_1 f_2\cdots f_{m})\|_{L^q_w} \lesssim_{N,s,t} \sum_{k=1}^{m} \|\flap f_k\|_{L^q_w}.
            \]
            
            In particular, this implies
            \begin{equation}\label{eq:frac-product}
            	\left\|\,|x|^{\frac{t}{\twostar}} \flap \Phi_i \right\|_{q} \lesssim_{N,s,t} \sum_{k\in\{i\}\cup I_i} \left\|\,|x|^{\frac{t}{\twostar}} \flap f_k \right\|_{q}.
            \end{equation}

            By applying the logarithmic bound from Proposition~\ref{basic-cutoff-norm-estimate} to a single cutoff function \( f_k = \varphi_{r_k, R_k} \) (or \( 1 - \varphi_{r_k, R_k} \)) for sufficiently large ratios \( R_k/r_k \), and substituting this estimate into \eqref{eq:frac-product}, we obtain, using the fact that \({\texttt{\#}}(\{i\}\cup I_i) \le \nu \), the following bound:
            \begin{equation}\label{eq:max-bound}
            	\left\|\,|x|^{\frac{t}{\twostar}} \flap \Phi_i \right\|_{q} \lesssim_{N,\nu,s,t} \max_{k\in\{i\}\cup I_i} \left( \log \frac{R_k}{r_k} \right)^{-\frac12-\frac{1}{\twostar}}.
            \end{equation}
            
           We next relate the size of the transition annuli to the interaction parameter $\delta$. Since the family $(V_i)_{i=1}^\nu$ is $\delta$-interacting, any two scales $\lambda_i$ and $\lambda_j$ with $i\neq j$ satisfy $\min\!\left(\frac{\lambda_i}{\lambda_j},\frac{\lambda_j}{\lambda_i}\right)
           \le \delta$, which implies
           \[
           \min_{i\neq j}\max\!\left\{\frac{\lambda_i}{\lambda_j},\frac{\lambda_j}{\lambda_i}\right\}
           \ge \frac1\delta.
           \]
           
           After reordering the bubbles in a decreasing order of concentration scales, distinct bubbles correspond to scales separated by a factor of at least $1/\delta$.

           In the construction of \(\Phi_i\), each cutoff annulus \([r_k, R_k]\) is positioned between two such scales to isolate \(V_i\). Therefore, it follows that
           \begin{equation}\label{eq:logdelta}
           \min_{\substack{k\in\{i\}\cup I_i\\ 1 \le i \le \nu}} \frac{R_k}{r_k} \ge C \delta^{-1}
           \end{equation}
           for some constant \(C > 0\) depending only on \(\nu\).
           
            \noindent\textit{Choice of \(\delta(N,\nu,\varepsilon,s,t)\)}: Combining the estimates from \eqref{eq:max-bound} and \eqref{eq:logdelta}, we have
            \[
            \left\|\,|x|^{\frac{t}{\twostar}} \flap \Phi_i \right\|_{q} \lesssim_{N,\nu,s,t} \left(\log \frac{1}{\delta}\right)^{-\frac12-\frac{1}{\twostar}}.
            \]
            
            Given \(\varepsilon > 0\), it suffices to ensure that
            $
            \left(\log \frac{1}{\delta}\right)^{-\frac12-\frac{1}{\twostar}} \lesssim \varepsilon,
            $
            which is equivalent to
            \begin{equation}\label{epsilon dependence of delta}
            \delta \le \exp\Big(-C\varepsilon^{-\frac{2\cdot\twostar}{2+\twostar}}\Big).
            \end{equation}
            
            With this choice of \(\delta\) in \eqref{epsilon dependence of delta}, the estimate in property (4) holds uniformly for every \(i\).\qed\vspace{1.8mm}
	
In order to verify that the orthogonality conditions~\eqref{Inn-bub}-\eqref{Inn-der-bub} with respect to $V_i$ and $\dot{V}_i$ are almost preserved under localization, we establish the following estimate.	
	\begin{Lemma}\label{almost-eigenfunction}
		Let $\rho \in \Hs$ be such that $\rho \perp_{\hs} \operatorname{span}\{V_i,\dot{V}_i\}.$ Let $\left(\Phi_i\right)_{i=1}^\nu$ denote the family of bump functions associated with the parameter $\varepsilon>0$, as introduced in Lemma~\ref{bumps}. Then, one has 
\[
\langle \rho \Phi_i , V_i \rangle_{\hs} = \mathcal o(1)\,\|\rho\|_{\hs}, 
\qquad 
\langle \rho \Phi_i , \dot{V}_i \rangle_{\hs} = \mathcal o(1)\,\|\rho\|_{\hs}.
\]
Here, the notation $\mathcal o(1)$ stands for a quantity that vanishes in the limit $\varepsilon \to 0$.
\end{Lemma}

	\Proof Let $\Psi$ denote either $V_i$ or $\dot{V}_i$, and set $p+1=\twostar$ and $q=\frac{2\cdot\twostar}{\twostar-2}$. Since $\rho\perp_{\hs}\Psi$ if and only if $\langle\rho,\Psi\rangle_{\hs}=0$, we combine integration by parts with H\"older's inequality (using $\tfrac{1}{p+1}+\tfrac12+\tfrac{1}{q}=1$), Remark~\ref{remark on eigenfunctions}-$(2)$, Lemma~\ref{bumps}-(1), and the Hardy--Sobolev inequality~\eqref{fhs0} to obtain \begin{align*}
		\Bigg|\langle\rho\Phi_i,\Psi\rangle_{\hs}\Bigg|&=\Bigg|\langle\rho,\Psi\rangle_{\hs}-\langle\rho(1-\Phi_i),\Psi\rangle_{\hs}\Bigg|=\mu\left|\int_{\Rn}\rho(1-\Phi_i)\Psi\frac{V_i^{p-1}}{|x|^t}\dx\right|\\&\le\mu\|\rho\|_{\lst}\left(\int_{\Rn}\Psi^2\frac{V_i^{p-1}}{|x|^t}\dx\right)^{\frac12}\left(\int_{\left\{\Phi_i<1\right\}}\frac{V_i^{p+1}}{|x|^t}\dx \right)^{\frac1q}\\&\le\mathcal{o}(1)\|\rho\|_{\hs},
	\end{align*} where $\mu=1$ \emph{or} $p$ according to $\Psi=V_i$ \emph{or} $\dot{V}_i$ respectively.\qed
	
	\subsection{Spectral gap inequality}\label{S2.3}
	Using the localization from Lemma~\ref{bumps}, the proof of Proposition~\ref{spectral-gap-inequality} proceeds naturally: localize, apply the Rayleigh quotient inequality~\eqref{rayleigh-quotient-inequality} to each bubble, and sum the results for the full estimate.
	\begin{Proposition}\label{spectral-gap-inequality}
		Let $s\in(0,1),\, N>2s,\, \nu \ge 1$ and $0<t<2s$. There exists $\delta=\delta(N,\nu,s,t)>0$ such that, if $(\alpha_iV_i)_{i=1}^\nu$ is a $\delta$-interacting family of bubbles and $\rho \in \Hs$ satisfies the orthogonality conditions~\eqref{Inn-bub}, \eqref{Inn-der-bub}, then there exists a constant $C\coloneqq C(N,\nu,s,t) \in (0,1)$ such that
		
		$$
		\int_{\Rn} \frac{|\sigma|^{\twostar-2} \rho^2}{|x|^t} \dx \le \frac{C}{\twostar-1}\,\|\rho\|_{\Hs}^2,
		\hspace{2mm} \text{where } \sigma = \sum_{i=1}^\nu \alpha_i V_i.
		$$
	\end{Proposition}
	\Proof Let us denote $p+1=\twostar$ and $q= \frac{2\cdot\twostar}{\twostar - 2}$. By choosing $\varepsilon > 0$ sufficiently small in Lemma~\ref{bumps}, one guarantees the existence of a constant $\delta > 0$, where $\delta = \mathcal{o}(1)$ as established in~\eqref{epsilon dependence of delta}. Furthermore, under this choice of $\varepsilon$, there exists a family of bump functions $(\Phi_i)_{i=1}^\nu$ satisfying the properties specified in Lemma~\ref{bumps}. Here, the notation $\mathcal{o}(1)$ refers to a positive quantity that tends to zero as $\varepsilon \to 0$.

	Write \begin{align*}
		\int_{\Rn}\frac{|\sigma|^{p-1}\rho^2}{|x|^t}\dx&=\left(\int_{\left\{\sum_{i=1}^\nu\Phi_i\ge1\right\}}+\int_{\left\{\sum_{i=1}^\nu\Phi_i<1\right\}}\right)\frac{|\sigma|^{p-1}\rho^2}{|x|^t}\dx\\&\eqqcolon I + II.
	\end{align*}
	
	Notice by H\"older inequality (with $\tfrac{1}{p+1}+\tfrac1q=\tfrac12$), Hardy--Sobolev inequality~\eqref{fhs0} and property (1) in Lemma~\ref{bumps},\begin{align*}
		II\le\|\rho\|_{\lst}^{2}\left(\int_{\left\{\sum_{i=1}^\nu\hspace{-2mm}\Phi_i<1\right\}}\frac{|\sigma|^{p+1}}{|x|^t}\dx\right)^{\frac{2}{q}}\hspace{-1.7mm}\lesssim_{N,s,t}\left(\sum_{i=1}^\nu\int_{\left\{\Phi_i<1\right\}}\frac{V_i^{p+1}}{|x|^t}\dx\right)^{\frac2q}\|\rho\|_{\hs}^2=\mathcal{o}(1)\|\rho\|_{\hs}^2.
	\end{align*}
	
	Property (2) in Lemma~\ref{bumps}, yields $|\sigma|\le(1+\mathcal{o}(1))V_i$ on the region $\left\{\Phi_i>0\right\}$ and then \begin{align*}
		I&\le\sum_{i=1}^\nu\int_{\{\Phi_i>0\}}\frac{|\sigma|^{p-1}\rho^2}{|x|^t}\dx\le (1+\mathcal{o}(1))\sum_{i=1}^\nu\int_{\{\Phi_i>0\}}\frac{V_i^{p-1}\rho^2}{|x|^t}\dx\\&\le(1+\mathcal{o}(1))\left(\sum_{i=1}^\nu\int_{\Rn}\frac{V_i^{p-1}(\rho\Phi_i)^2}{|x|^t}\dx+\mathcal{o}(1)\|\rho\|_{\hs}^2\right).
	\end{align*}
	
	Since $\rho$ satisfies the orthogonality conditions~\eqref{Inn-bub}, \eqref{Inn-der-bub}, we have for each $i$, $\rho\perp_{\hs}\operatorname{span}\{V_i,\dot{V}_i\}$ and Lemma~\ref{almost-eigenfunction} applies. Now taking $f_\varepsilon=\frac{\rho\Phi_i}{\|\rho\|_{\hs}}$ in Proposition~\ref{spectral-gap-inequality-almost-eigenfunction}, we obtain \[\int_{\Rn}\frac{V_{i}^{p-1}(\rho\Phi_i)^2}{|x|^t}\dx\le\frac{1}{\eta_3}\int_{\Rn}\left|\flap (\rho\Phi_i)\right|^2\dx+\mathcal{o}(1)\|\rho\|_{\hs}^2.\]
	
	Next write $\flap(\rho\Phi_i)=\rho\flap\Phi_i+\Phi_i\flap\rho+\mathcal{C}[\rho,\Phi_i]$ and use Theorem~\ref{wkpv-estimate}, with the choices $$\alpha_1=0,\alpha_2=\frac s2,\hspace{2mm}p_1=p+1,p_2=q,\hspace{2mm} w=1\in\mathcal{A}_2,w_1=|x|^{-t}\in\mathcal{A}_{p+1},w_2=|x|^{qt/\twostar}\in\mathcal{A}_q.$$
	
	This entails via an application of H\"older inequality (with $\tfrac12=\tfrac1{p+1}+\tfrac1q$), the Hardy--Sobolev inequality~\eqref{fhs0} and property (4) in Lemma~\ref{bumps},\begin{align*}
		\left\|\flap(\rho\Phi_i)\right\|_2&\le\left\|\rho\flap\Phi_i\right\|_2+\left\|\Phi_i\flap\rho\right\|_2+\left\|\mathcal{C}[\rho,\Phi_i]\right\|_2\\&\le\left\|\rho\right\|_{\lst}\left\||x|^{\frac t\twostar}\flap\Phi_i\right\|_q\hspace{-1mm}+\left\|\Phi_i\flap\rho\right\|_2\hspace{-1mm}+C\left\|\rho\right\|_{L^{p_1}_{w_1}}\left\|\flap\Phi_i\right\|_{L^{p_2}_{w_2}}\\&\le\mathcal{o}(1)\left\|\rho\right\|_{\hs}+\left\|\Phi_i\flap\rho\right\|_2.
	\end{align*}
	
	Since the bump functions $(\Phi_i)_{i=1}^\nu$ have almost disjoint supports with controlled overlaps, \begin{align*}
		\sum_{i=1}^\nu\int_{\Rn}\frac{V_{i}^{p-1}(\rho\Phi_i)^2}{|x|^t}\dx&\le\frac{1}{\eta_3}\left(\sum_{i=1}^\nu\left\|\Phi_i\flap\rho\right\|_2^2+\mathcal{o}(1)\|\rho\|_{\hs}^2\right)+\mathcal{o}(1)\|\rho\|_{\hs}^2\\&\le\left(\frac{1}{\eta_3}+\mathcal{o}(1)\right)\|\rho\|_{\hs}^2.
	\end{align*}
	
	Thus we get \[I+II\le\left(\frac{1}{\eta_3}+\mathcal{o}(1)\right)\|\rho\|_{\hs}^2+\mathcal{o}(1)\|\rho\|_{\hs}^2=\left(\frac{1}{\eta_3}+\mathcal{o}(1)\right)\|\rho\|_{\hs}^2,\] which essentially proves the result since $\eta_3>p=\twostar-1$.\qed

\subsection{Taming interaction integrals and the robustness of bubble coefficients}\label{S2.4}
In this subsection we quantify how weakly-interacting bubbles influence one another and establish the stability of their coefficients under small perturbations. These bubble-bubble interaction bounds and coefficient robustness estimates furnish the final structural ingredient needed to complete the proof of Theorem~\ref{quantitative-stability}.

\begin{Proposition}\label{bub-bub-man}
Let $N>2s$ and $\nu\ge1$ be integers. For every $\varepsilon>0$ there exists
$\delta=\delta(N,\nu,\varepsilon,s,t)>0$ such that the following holds.
Assume that $(\alpha_iV_i)_{i=1}^\nu$ is a $\delta$-interacting family and that
\[
u=\sum_{i=1}^{\nu}\alpha_iV_i+\rho,
\]
where $\rho$ satisfies the orthogonality conditions \eqref{Inn-bub}-\eqref{Inn-der-bub} and
$\|\rho\|_{\Hs}<1$. Then the bubble coefficients are stable in the sense that, for each
$1\le i\le\nu$,
\begin{equation}\label{Est-4}
|\alpha_i-1|
\lesssim_{N,\nu,s,t}
\mathcal{o}(1)\,\|\rho\|_{\Hs}
+\Gamma(u)
+\|\rho\|_{\Hs}^{\min(2,\twostar-1)}.
\end{equation}

Moreover, bubble-bubble interactions are quantitatively controlled: for every distinct
pair $(i,j)$,
\begin{equation}\label{Est-5}
\int_{\Rn}\frac{V_i^{\twostar-1}V_j}{|x|^t}\dx
\lesssim_{N,\nu,s,t}
\mathcal{o}(1)\,\|\rho\|_{\Hs}
+\Gamma(u)
+\|\rho\|_{\Hs}^{\min(2,\twostar-1)}.
\end{equation}

Here $\mathcal{o}(1)$ denotes a quantity that tends to $0$ as
$\varepsilon\to0$.
\end{Proposition}

	\Proof We begin by ordering the bubbles according to their concentration scales,
\[
\lambda_1 > \lambda_2 > \cdots > \lambda_\nu,
\]
so that $V_1$ is the most concentrated profile around the origin. The proof proceeds by induction on the index $i$. The base case $i=1$ follows from the single-bubble analysis carried out in \cite[Theorem~1.5]{CS25}. Assume now that the conclusion of Proposition~\ref{bub-bub-man} holds for every index $j<i$, where $1<i\le\nu$. For the $i$-th bubble we introduce the auxiliary superposition
\[
\widetilde{\sigma}_i = \sum_{\substack{1\le j\le\nu \\ j\neq i}} \alpha_j V_j.
\]

Recall that $p+1=\twostar$, $q=\frac{2\cdot\twostar}{\twostar-2}$, and $\sigma=\sum_{i=1}^{\nu}\alpha_iV_i$.  

\medskip
Fix $\varepsilon>0$.  By Lemma~\ref{bumps}, there exists $\delta=\delta(N,\nu,\varepsilon,s,t)>0$ such that the $\delta$-interacting family $(V_i)_{i=1}^\nu$ \big(as $(\alpha_iV_i)_{i=1}^\nu$ is $\delta$-interacting\big) admits bump functions $(\Phi_i)_{i=1}^\nu$ satisfying the localization properties stated therein.  These functions verify $0\le \Phi_i\le 1$, and the interaction condition yields $|\alpha_i-1|\le \delta$ for every $1\le i\le\nu$. Since \eqref{epsilon dependence of delta} ensures $\delta=\mathcal{o}(1)$ as $\varepsilon\to 0$, by choosing $\varepsilon$ sufficiently small we may further assume $\tfrac{1}{2}<\min_{i}\alpha_i\le\max_i\alpha_i<\tfrac32$.

By invariance under the symmetries of the fractional Hardy--Sobolev functional, we may assume without loss of generality that $V_i=V_{s,t}^{1}$, i.e.\ $\lambda_i=1$; see \cite[\S\,2.2]{CS25}. Indeed, all quantities appearing in \eqref{Est-4} and \eqref{Est-5} are invariant under the natural scaling~\eqref{dilation map} introduced in the Introduction.

Set $r_i^-\coloneqq\max_{k\in I_i}r_k$, with the convention $r_i^-=0$ if $I_i=\varnothing$. By the explicit construction of the cutoffs in \eqref{def cutoff main}, we have $\{\Phi_i>0\}=\{x\in\Rn:r_i^-<|x|<R_i\}$, and hence $\supp\Phi_i=\{x\in\Rn:r_i^-\le |x|\le R_i\}$. Moreover, under the normalization $\lambda_i=1$, the choice of the cutoff radii gives $r_i^-=\mathcal{o}(1),R_i\ge1$. Consequently,
\begin{equation}\label{supp cond}
\{x\in\Rn:\mathcal o(1)\le |x|\le R_i\}=\supp\Phi_i.
\end{equation}

With these preparations, we consider the identity
\begin{align*}
  (\alpha_i-\alpha_i^p)\frac{V_i^p}{|x|^t}
  - p\frac{(\alpha_iV_i)^{p-1}\widetilde{\sigma}_i}{|x|^t}
  &= -(-\Delta)^s\rho
     + \Big( (-\Delta)^s u - \frac{u|u|^{p-1}}{|x|^t} \Big)
     - \sum_{\substack{1\le j\le\nu \\ j\neq i}}
       \frac{\alpha_jV_j^p}{|x|^t}   \\
  &\quad + p\frac{(\alpha_iV_i)^{p-1}\rho}{|x|^t}
     + \Big[
         \frac{(\sigma+\rho)|\sigma+\rho|^{p-1}}{|x|^t}
         - \frac{\sigma^p}{|x|^t}
         - p\frac{\sigma^{p-1}\rho}{|x|^t}
       \Big] \\
  &\quad + \Big[
         p\frac{\sigma^{p-1}\rho}{|x|^t}
         - p\frac{(\alpha_iV_i)^{p-1}\rho}{|x|^t}
       \Big]   \\
  &\quad + \Big[
        \frac{(\alpha_iV_i+\widetilde{\sigma}_i)^p}{|x|^t}
        - \frac{(\alpha_iV_i)^p}{|x|^t}
        - p\frac{(\alpha_iV_i)^{p-1}\widetilde{\sigma}_i}{|x|^t}
       \Big].
\end{align*}

The first term, $(\alpha_i-\alpha_i^p)\frac{V_i^p}{|x|^t}$, contains the coefficient error and is the starting point for establishing the estimate \eqref{Est-4}. The interaction term $\frac{V_i^{p-1}\widetilde{\sigma}_i}{|x|^t}$ will be responsible for deriving the mixed integral bound \eqref{Est-5}. Since our objective is to obtain a control in terms of $\Gamma(u)$, we introduce the fractional Laplacian explicitly on the right-hand side above.

Using \eqref{Ineq-1} together with the localization property in Lemma~\ref{bumps}-(2) and the bound $\max_i\alpha_i<\tfrac32$, we work on the region $\{\Phi_i>0\}$, where the interaction structure imposed by the lemma yields effective control of all nonlinear remainders. In this regime, we obtain the following pointwise estimates:
\begin{enumerate}
\item[\textup{(i)}] $\sum_{j\neq i}\frac{\alpha_j V_j^{p}}{|x|^t}
\le \varepsilon^{p-1}\frac{V_i^{p-1}\widetilde{\sigma}_i}{|x|^t}
= \mathcal{o}(1)\Big(\frac{V_i^{p-1}\widetilde{\sigma}_i}{|x|^t}\Big)$,

\item[\textup{(ii)}] $\left|
\frac{\sigma^{p-1}\rho-(\alpha_iV_i)^{p-1}\rho}{|x|^t}
\right|
\lesssim_{p-1}
\frac{(\alpha_iV_i)^{p-2}\widetilde{\sigma}_i|\rho|}{|x|^t}
+\chi_{p-1>2}\frac{(\alpha_iV_i)^{p-3}\widetilde{\sigma}_i^{2}|\rho|}{|x|^t}
+\frac{\widetilde{\sigma}_i^{p-1}|\rho|}{|x|^t}= \mathcal{o}(1)\Big(\frac{V_i^{p-1}|\rho|}{|x|^t}\Big)$,

\item[\textup{(iii)}] $\left|
\frac{(\alpha_iV_i+\widetilde{\sigma}_i)^p-(\alpha_iV_i)^p
-p(\alpha_iV_i)^{p-1}\widetilde{\sigma}_i}{|x|^t}
\right|
\lesssim_p
\chi_{p>2}\frac{(\alpha_iV_i)^{p-2}\widetilde{\sigma}_i^{2}}{|x|^t}
+\frac{\widetilde{\sigma}_i^{p}}{|x|^t}= \mathcal{o}(1)\Big(\frac{V_i^{p-1}\widetilde{\sigma}_i}{|x|^t}\Big)$,

\item[\textup{(iv)}] $\left|
\frac{(\sigma+\rho)|\sigma+\rho|^{p-1}-\sigma^p-p\sigma^{p-1}\rho}{|x|^t}
\right|
\lesssim_{p}
\frac{|\rho|^p}{|x|^t}
+\chi_{p>2}(1+\mathcal{o}(1))\frac{V_i^{p-2}\rho^2}{|x|^t}$.
\end{enumerate}

In the estimates above, we used the indicator $\chi_{\cdot>2}$ from \eqref{Ineq-1} to treat the two different dimensional regimes simultaneously. Since $p+1=\twostar=\frac{2(N-t)}{N-2s}$, we have $p>2 \iff 2s<N<6s-2t,\text{ and }p\le 2 \iff N\ge 6s-2t$. Thus the quadratic correction $\chi_{p>2}|a|^{p-2}|b|^2$ appears only in the low-dimensional regime $2s<N<6s-2t$ (equivalently $p>2$) and drops out when $N\ge 6s-2t$ (equivalently $p\le2$).

Combining the contributions in {(i)}-{(iv)} on the region where $\Phi_i>0$, we conclude that
\begin{align*}
\Bigg|
   (\alpha_i-\alpha_i^p)\frac{V_i^{\,p}}{|x|^t}
   -(p\alpha_i^{p-1}+\mathcal{o}(1))\frac{V_i^{\,p-1}\widetilde{\sigma}_i}{|x|^t}
   +(-\Delta)^s\rho
   -\Big((-\Delta)^su-\frac{|u|^{p-1}u}{|x|^t}\Big)
   -p\frac{(\alpha_iV_i)^{p-1}\rho}{|x|^t}
\Bigg|
\\[0.3em]
\lesssim_{N,s,t}
   \frac{|\rho|^{p}}{|x|^t}
   +\chi_{p>2}\frac{V_i^{\,p-2}\rho^2}{|x|^t}
   +\mathcal{o}(1)\bigg(\frac{V_i^{\,p-1}|\rho|}{|x|^t}\bigg).
\end{align*}

       We now test the last estimate against $\Psi\Phi_i$, where $\Psi\in\{V_i,\dot V_i\}$. Using the orthogonality relations \eqref{Inn-bub}-\eqref{Inn-der-bub}, and the triangle inequality, one observes
\begin{align}
&\left|\int_{\Rn}
   \Big[
      (\alpha_i-\alpha_i^{p})\frac{V_i^{p}}{|x|^t}
      -(p\alpha_i^{p-1}+\mathcal{o}(1))\frac{V_i^{p-1}\widetilde{\sigma}_i}{|x|^t}
   \Big]
   \Psi\Phi_i\dx\right|
   \nonumber\\
&\qquad\lesssim_{N,s,t}
   \left|\langle\rho,\Psi\Phi_i\rangle_{\dot{H}^s}\right|
   +\left|
      \prescript{}{(\dot{H}^s)'}{\langle}
      (-\Delta)^su-\frac{|u|^{p-1}u}{|x|^t},
      \Psi\Phi_i
      \rangle_{\dot{H}^s}
    \right|
   +\left|\int_{\Rn}\frac{\rho\,\Phi_i\,\Psi\,V_i^{p-1}}{|x|^t}\dx\right|
   \nonumber\\
&\qquad\quad
   +\int_{\Rn}\frac{|\rho|^p|\Psi|\Phi_i}{|x|^t}\dx
   +\chi_{p>2}\int_{\Rn}\frac{\rho^2|\Psi|\Phi_iV_i^{p-2}}{|x|^t}\dx
   +\mathcal{o}(1)\int_{\Rn}\frac{|\rho||\Psi|\Phi_iV_i^{p-1}}{|x|^t}\dx
   \nonumber\\
&\qquad\eqqcolon I_1+I_2+\cdots+\chi_{p>2}I_5+\mathcal{o}(1)I_6.
\label{Esti-7}
\end{align}

We begin by noting the following claim, which quantifies the contribution of each term appearing on the right-hand side of \eqref{Esti-7}.

\medskip
\noindent\textbf{Claim.} One has $I_1=\mathcal{o}(1)\|\rho\|_{\Hs}$, $I_2\lesssim\Gamma(u)$, $I_3=\mathcal{o}(1)\|\rho\|_{\Hs}$, $I_4\lesssim\|\rho\|^p_{\Hs}$, $I_5\lesssim\|\rho\|^2_{\Hs}$ and $I_6\lesssim\|\rho\|_{\Hs}$, with implicit constants depending only on $N,s,t$.

        The orthogonality relations \eqref{Inn-bub}-\eqref{Inn-der-bub} imply $\langle\rho,\Psi\rangle_{\hs}=0$. Using this and decomposing $I_1$ in terms of the commutator operator $\mathcal{C}[\cdot,\cdot]$, we have by the triangle inequality,\begin{align*}
        I_1=\left|\langle\rho,\Phi_i\Psi\rangle_{\hs}\right|&= \left|\langle\rho, (1-\Phi_i)\Psi\rangle_{\hs}\right|\\
           &=\left|\int_{\Rn}(-\Delta)^{\frac{s}{2}}\rho(-\Delta)^{\frac{s}{2}}\left((1-\Phi_i)\Psi\right)\dx\right|\\
           &\le\left|\int_{\Rn}\mathcal{C}[1-\Phi_i,\Psi](-\Delta)^{\frac{s}{2}}\rho\dx\right| + \left|\int_{\Rn}(1-\Phi_i)(-\Delta)^{\frac{s}{2}}\rho(-\Delta)^{\frac{s}{2}}\Psi\dx\right|\\
           &\qquad +\left|\int_{\Rn} \Psi(-\Delta)^{\frac{s}{2}}\rho(-\Delta)^{\frac{s}{2}}(1-\Phi_i)\dx\right|\eqqcolon J_1+J_2+J_3.
        \end{align*}

        \noindent\textit{Estimate for $J_1$.} Applying Theorem~\ref{wkpv-estimate} with $\alpha_1 = 0,\, \alpha_2 = \frac{s}{2},\,w_1 = |x|^{-t},\, w_2 = |x|^{\frac{qt}{p+1}},\,p_1 = p + 1,\,p_2 = q$, Lemma~\ref{bumps}-(4), H\"older's inequality (with exponents $\frac{1}{2} = \frac{1}{p+1} + \frac{1}{q}$), Remark~\ref{remark on eigenfunctions}-$(2)$, and \eqref{bubble normalization}, we find 
        \begin{align*}
            J_1=&\left|\int_{\Rn}\mathcal{C}[\Psi,1-\Phi_i] (-\Delta)^{\frac{s}{2}}\rho\dx\right| \le\left\|\mathcal{C}[\Psi,1-\Phi_i]\right\|_{2}\|(-\Delta)^{\frac{s}{2}}\rho\|_{2}\\
            &\hspace{5em}\lesssim\|\Psi\|_{L^{p+1}_{|x|^{-t}}}\left\||x|^\frac{t}{\twostar}\flap\Phi_i\right\|_{q} \|(-\Delta)^{\frac{s}{2}}\rho\|_{2}\\
            &\hspace{5em}\lesssim\varepsilon \|(-\Delta)^{\frac{s}{2}}\rho\|_{2}=\mathcal{o}(1)\|\rho\|_{\hs}.
        \end{align*}

        \noindent\textit{Estimate for $J_2$.} H\"older's inequality entails,
        \begin{equation*}
            J_2\le \|(-\Delta)^{\frac{s}{2}}\rho\|_2\left\|(1-\Phi_i)(-\Delta)^{\frac{s}{2}}\Psi\right\|_2\eqqcolon J_{2}'\|\rho\|_{\hs}.
        \end{equation*}
        
        To show $J_2'=\mathcal{o}(1)$, use $(1-\Phi_i)^2\le (1-\Phi_i)$ and integrate by parts in the following decomposition of $J_2'$ via the commutator operator.
        \begin{align*}
            |J_2'|^2&\le \int_{\Rn} (1-\Phi_i) (-\Delta)^{\frac{s}{2}}\Psi \cdot(-\Delta)^{\frac{s}{2}}\Psi\dx\\
            &=\int_{\Rn} \Psi(-\Delta)^{\frac{s}{2}}\left((1-\Phi_i)(-\Delta)^{\frac{s}{2}}\Psi\right)\dx\\
            &=\int_{\Rn}\Psi\cdot\Big(\mathcal{C}[1-\Phi_i,(-\Delta)^{\frac{s}{2}}\Psi]+(1-\Phi_i)(-\Delta)^{s}\Psi+ (-\Delta)^{\frac{s}{2}}\Psi (-\Delta)^{\frac{s}{2}}(1-\Phi_i)\Big)\dx\\
            &\eqqcolon J'_{21}+J'_{22}+J'_{23}.
        \end{align*}
        
        Now using H\"older's inequality, Theorem~\ref{wkpv-estimate} with the choices $w=|x|^{\frac{t}{p}}\in\mathcal A_{(p+1)/p},\,w_1=|x|^{\frac{2t}{p-1}}\in\mathcal A_q,w_2=1\in\mathcal A_2,\,\alpha_1=\frac{s}{2},\alpha_2=0,\,p_1=q,p_2=2$, Remark~\ref{remark on eigenfunctions}-$(2)$, \eqref{bubble normalization}, and finally Lemma~\ref{bumps}-(4), we deduce
        \begin{align*}
            |J'_{21}| &\le \int_{\Rn}|\Psi\cdot \mathcal{C}[1-\Phi_i,(-\Delta)^{\frac{s}{2}}\Psi]|\dx\\
            &\le\|\Psi\|_{L^{p+1}_{|x|^{-t}}}\left\|\mathcal{C}[1-\Phi_i,(-\Delta)^{\frac{s}{2}}\Psi]\right\|_{L^{\frac{p+1}{p}}_{|x|^{\frac{t}{p}}}}\\
            &\le\|\Psi\|_{L^{p+1}_{|x|^{-t}}}\left\||x|^{\frac{t}{p+1}}(-\Delta)^{\frac{s}{2}}\Phi_i\right\|_{q}\|(-\Delta)^{\frac{s}{2}}\Psi\|_{2}\lesssim \varepsilon\|V_i\|_{L^{p+1}_{|x|^{-t}}}^{\frac{p+3}{2}} =\mathcal{o}(1).
        \end{align*}

        In a similar fashion from the eigenvalue equation~\eqref{eigenvalue problem}, one gets
       \begin{equation*}
             J'_{22}=\int_{\Rn}(1-\Phi_i)\Psi(-\Delta)^{s}\Psi\dx= \mu \int_{\Rn}(1-\Phi_i)\frac{V_i^{p-1}|\Psi|^2}{|x|^t}\dx\lesssim \int_{\{\Phi_i<1\}}\frac{V_i^{p+1}}{|x|^t}=\mathcal{o}(1).    
       \end{equation*}

        Finally by H\"older's inequality with exponents $\frac{1}{2}+\frac{1}{q}+\frac{1}{p+1}=1$, Remark~\ref{remark on eigenfunctions}-$(2)$, \eqref{bubble normalization}, and Lemma~\ref{bumps}-(4), we obtain,
        \begin{align*}
            |J'_{23}|&\le\int_{\Rn}|\Psi\cdot(-\Delta)^{\frac{s}{2}}\Psi\cdot(-\Delta)^{\frac{s}{2}}(1-\Phi_i)|\dx\\&\le \|\Psi\|_{L^{p+1}_{|x|^{-t}}} \left\||x|^{\frac{t}{p+1}}(-\Delta)^{\frac{s}{2}}\Phi_i\right\|_{q}\|(-\Delta)^{\frac{s}{2}}\Psi\|_2=\mathcal{o}(1).
        \end{align*}

        Hence we have $J_2\le J_2'\|\rho\|_{\hs}=\mathcal{o}(1)\|\rho\|_{\hs}$. 

        \medskip
        \noindent\textit{Estimate for $J_3$.} H\"older's inequality together with Remark~\ref{remark on eigenfunctions}-$(2)$, \eqref{bubble normalization}, and Lemma~\ref{bumps}-(4) yield
        \begin{align*}
            J_3= \left|\int_{\Rn}  (-\Delta)^{\frac{s}{2}}\rho\cdot\Psi\cdot(-\Delta)^{\frac{s}{2}}(1-\Phi_i)\dx \right|&\le\|(-\Delta)^{\frac{s}{2}}\rho\|_{2}\left\||x|^{\frac{-t}{\twostar}}\Psi\right\|_{p+1} \left\||x|^\frac{t}{\twostar}\flap\Phi_i\right\|_{q}\\
            &\lesssim\varepsilon \|(-\Delta)^{\frac{s}{2}}\rho\|_{2}=\mathcal{o}(1)\|\rho\|_{\hs}.
        \end{align*}
        
        Combining the estimates on $J_1,\,J_2$ and $J_3$, we obtain the stated bound for $I_1$. For the second term in \eqref{Esti-7}, we have
\[
I_2
 = \left|\prescript{}{(\dot{H}^s)'}{\langle}
      (-\Delta)^s u - \frac{|u|^{p-1}u}{|x|^t},
      \Psi\Phi_i
    \rangle_{\dot{H}^s}\right|
 \le \Gamma(u)\,\left\|(-\Delta)^{\frac{s}{2}}(\Psi\Phi_i)\right\|_{2}
 \lesssim \Gamma(u).
\]

To justify the final bound in $I_2$, we again invoke Theorem~\ref{wkpv-estimate}, together with H\"older's inequality, Remark~\ref{remark on eigenfunctions}-$(2)$, \eqref{bubble normalization}, and Lemma~\ref{bumps}-(4).  These yield
\begin{align*}
\left\|(-\Delta)^{\frac{s}{2}}(\Psi\Phi_i)\right\|_{2}
&\le \|\mathcal{C}[\Psi,\Phi_i]\|_{2}
   + \|\Phi_i(-\Delta)^{\frac{s}{2}}\Psi\|_{2}
   + \|\Psi(-\Delta)^{\frac{s}{2}}\Phi_i\|_{2} \\
&\lesssim
     \|\Psi\|_{L^{p+1}_{|x|^{-t}}}
     \left\||x|^{\frac{t}{p+1}}(-\Delta)^{\frac{s}{2}}\Phi_i\right\|_{q}
     + \|(-\Delta)^{\frac{s}{2}}\Psi\|_{2}\lesssim 1 + \mathcal{o}(1).
\end{align*}

This proves the claimed bound for $I_2$. Now Lemma~\ref{almost-eigenfunction} directly gives \begin{equation*}
           I_3=\Big|\int_{\Rn}\frac{\rho\Phi_i\Psi V_i^{p-1}}{|x|^t}\dx\Big| =\tfrac{1}{\mu}\Big|\langle\rho\Phi_i,\Psi\rangle_{\hs}\Big|=\mathcal{o}(1)\|\rho\|_{\hs}.
       \end{equation*}
       
       For the other terms using H\"older's inequality, $\Phi_i\le1$, Remark~\ref{remark on eigenfunctions}-$(2)$, \eqref{bubble normalization}, and \eqref{fhs0}, we see
       \begin{align*}
           &I_4=\int_{\Rn}\frac{|\rho|^p|\Psi|\Phi_i}{|x|^t}\dx \lesssim \int_{\Rn}\frac{|\rho|^pV_i}{|x|^t}\dx\lesssim \|\rho\|_{L^{p+1}_{|x|^{-t}}}^p\|V_i\|_{L^{p+1}_{|x|^{-t}}} \lesssim\|\rho\|_{\hs}^p,\\
           &I_5=\int_{\Rn}\frac{\rho^2|\Psi|\Phi_iV_i^{p-2}}{|x|^t}\dx\lesssim\int_{\Rn}\frac{\rho^2V_i^{p-1}}{|x|^t}\dx \lesssim \|\rho\|_{L^{p+1}_{|x|^{-t}}}^2\|V_i\|_{L^{p+1}_{|x|^{-t}}}^{p-1}\lesssim  \|\rho\|_{\hs}^2,\\
           &I_6=\int_{\Rn}\frac{|\rho||\Psi|\Phi_iV_i^{p-1}}{|x|^t}\dx\lesssim\int_{\Rn}\frac{|\rho|V_i^p}{|x|^t}\dx \lesssim \|\rho\|_{L^{p+1}_{|x|^{-t}}}\|V_i\|_{L^{p+1}_{|x|^{-t}}}^p\lesssim\|\rho\|_{\hs}.
       \end{align*} 
       
       This concludes the claim. Combining the bounds for $I_1,\ldots,I_6$ established above and substituting them into
\eqref{Esti-7}, we obtain by using the hypothesis $\|\rho\|_{\hs}<1$,
\begin{align}
&\Bigg|
\int_{\Rn}
\frac{
(\alpha_i-\alpha_i^p)V_i^{p}
-\left(p\alpha_i^{p-1}+\mathcal{o}(1)\right)V_i^{\,p-1}\widetilde{\sigma}_i
}{|x|^{t}}
\,\Psi\Phi_i\dx
\Bigg|\nonumber\\&\quad
\lesssim_{N,s,t}
\mathcal{o}(1)\,\|\rho\|_{\hs}
+\Gamma(u)
+\|\rho\|_{\hs}^{\min(2,p)}.
\label{est-1}
\end{align}

       We now decompose
\[
\widetilde{\sigma}_i=\widetilde{\sigma}_{<i}+\widetilde{\sigma}_{>i},
\quad
\widetilde{\sigma}_{<i}=\sum_{j<i}\alpha_jV_j,\quad
\widetilde{\sigma}_{>i}=\sum_{j>i}\alpha_jV_j.
\]

By the induction hypothesis applied to \eqref{Est-5}, the interaction estimate is already available for all indices $j<i$.  Consequently, for $\Psi\in\{V_i,\dot V_i\}$ we get by using Remark~\ref{remark on eigenfunctions}, 
\begin{equation}\label{est-2}
\int_{\Rn}\frac{V_i^{p-1}\,\widetilde{\sigma}_{<i}\,|\Psi|\Phi_i}{|x|^t}\dx
\lesssim_{N,s,t}
\int_{\Rn}\frac{V_i^{p}\widetilde{\sigma}_{<i}}{|x|^t}\dx
\lesssim_{N,\nu,s,t}
\mathcal{o}(1)\,\|\rho\|_{\hs}
+\Gamma(u)
+\|\rho\|_{\hs}^{\min(2,p)}.
\end{equation}

Observe that, at the terminal index $i=\nu$, the contribution $\widetilde{\sigma}_{>\nu}$ vanishes identically, so that $\widetilde{\sigma}_\nu=\widetilde{\sigma}_{<\nu}$. Consequently, taking $\Psi=V_\nu$ in~\eqref{est-1} and using
\eqref{est-2}, we obtain
\begin{align*}
&|\alpha_\nu-\alpha_\nu^p|
\int_{\Rn}\frac{V_\nu^{p+1}\Phi_\nu}{|x|^t}\dx\\
&\le
\left|\int_{\Rn}
\frac{
(\alpha_\nu-\alpha_\nu^p)V_\nu^{p+1}
-\left(p\alpha_\nu^{p-1}+\mathcal{o}(1)\right)
V_\nu^p\widetilde{\sigma}_{<\nu}
}{|x|^t}\Phi_\nu\dx\right|+
\left(p\alpha_\nu^{p-1}+\mathcal{o}(1)\right)
\int_{\Rn}
\frac{V_\nu^p\widetilde{\sigma}_{<\nu}\Phi_\nu}{|x|^t}\dx\\
&\lesssim_{N,\nu,s,t}
\mathcal{o}(1)\,\|\rho\|_{\Hs}
+\Gamma(u)
+\|\rho\|_{\Hs}^{\min(2,p)}.
\end{align*}

By property~(1) of Lemma~\ref{bumps}, the integral on the left is bounded away from zero. Indeed, since $0\le\Phi_\nu\le1$, we have
\[
\int_{\mathbb R^N}\frac{V_\nu^{p+1}\Phi_\nu}{|x|^t}\dx
\ge
\int_{\{\Phi_\nu=1\}}
\frac{V_\nu^{p+1}}{|x|^t}\dx=
\int_{\mathbb R^N}\frac{V_\nu^{p+1}}{|x|^t}\dx
-
\int_{\{\Phi_\nu<1\}}
\frac{V_\nu^{p+1}}{|x|^t}\dx\ge
(1-\varepsilon)
\mu_{s,t}^{\frac{N-t}{2s-t}}.
\] 

Hence
\[
|\alpha_\nu-\alpha_\nu^p|
\lesssim_{N,\nu,s,t}
\mathcal{o}(1)\,\|\rho\|_{\Hs}
+\Gamma(u)
+\|\rho\|_{\Hs}^{\min(2,p)}.
\]

Since $\frac12<\alpha_\nu<\frac32$ and $p>1$, this yields the required estimate~\eqref{Est-4} for $i=\nu$.   

Next we control the contribution of $\widetilde{\sigma}_{>i}$. Note that,
\[
\sup_{\{\Phi_i>0\}}V_j
=\lambda_j^{\frac{N-2s}{2}}V_{s,t}(\lambda_jr_i^-)
=(1-\mathcal o(1))V_j(0),
\]
where $r_i^-=\inf\bigl\{|x|:\Phi_i(x)>0\bigr\}=\mathcal o(1)$ thanks to~\eqref{supp cond}. Here we used the continuity of the radially decreasing $V_{s,t}$ at the origin, together with the fact that $\lambda_j<\lambda_i=1$ for $j>i$. In particular, $0<\frac{V_{s,t}(\lambda_jr_i^-)}{V_{s,t}(0)}\le 1$, and the convergence to $1$ is uniform in $j>i$, since $\lambda_jr_i^-\le r_i^-=\mathcal o(1)$.

      Using Lemma~\ref{bumps}-(3), we obtain on the region where $\Phi_i>0$, the following estimate.
\begin{align*}
\widetilde{\sigma}_{>i}(x)
   = \sum_{j>i}\alpha_j V_j(x)
   &\ge \sum_{j>i}\alpha_j\inf_{\{\Phi_i>0\}}V_j
     \ge \frac{1}{1+\varepsilon}\sum_{j>i}\alpha_j\sup_{\{\Phi_i>0\}}V_j\ge\frac{1-\mathcal{o}(1)}{1+\varepsilon}\widetilde{\sigma}_{>i}(0).
\end{align*}

On the other hand, since $V_j$'s are radially symmetric and strictly decreasing, one has $\widetilde{\sigma}_{>i}(x)\le\widetilde{\sigma}_{>i}(0)$ for all $x\in\Rn$. Hence $(1-\mathcal{o}(1))\widetilde{\sigma}_{>i}(0)\Phi_i(x)\le\widetilde{\sigma}_{>i}(x)\Phi_i(x)
   \le\widetilde{\sigma}_{>i}(0)\Phi_i(x)$. Therefore,
\begin{equation}\label{Est-11}
\widetilde{\sigma}_{>i}(x)\Phi_i(x)=(1-\mathcal{o}(1))\widetilde{\sigma}_{>i}(0)\Phi_i(x),\quad\text{ on }\{x:\Phi_i(x)>0\}.
\end{equation}

       If $\alpha_i=1$, then \eqref{est-1} immediately closes the induction at level $j=i$, and the proof of Proposition~\ref{bub-bub-man} is complete. We therefore assume $\alpha_i\neq 1$. Since $\min_j\alpha_j>\tfrac12$, the coefficient deficit satisfies $\alpha_i-\alpha_i^{p}\neq 0$. In this case we introduce $\theta \coloneqq \frac{p\alpha_i^{p-1}\widetilde{\sigma}_{>i}(0)}{\alpha_i-\alpha_i^{p}}$, which is a well-defined fixed quantity and satisfies $\theta\neq0$.

       Using \eqref{Est-11} together with the decomposition
$\widetilde{\sigma}_i=\widetilde{\sigma}_{<i}+\widetilde{\sigma}_{>i}$ and the triangle
inequality, we rewrite the left-hand side of \eqref{est-1} in terms of the fixed parameter
$\theta$ as follows:
\begin{align}\label{est0}
&\Bigg|\int_{\Rn}
\frac{
(\alpha_i-\alpha_i^p)V_i^p
-\left(p\alpha_i^{p-1}+\mathcal{o}(1)\right)V_i^{p-1}\widetilde{\sigma}_i
}{|x|^t}
\,\Psi\Phi_i\dx
\Bigg|
\nonumber\\
&=|\alpha_i-\alpha_i^p|
\Bigg|
\int_{\Rn}
\frac{
V_i^p-(1\pm\mathcal{o}(1))\theta V_i^{p-1}
}{|x|^t}
\,\Psi\Phi_i\dx
-\frac{(1+\mathcal{o}(1))\theta}{\widetilde{\sigma}_{>i}(0)}
\int_{\Rn}
\frac{\widetilde{\sigma}_{<i}V_i^{p-1}}{|x|^t}
\,\Psi\Phi_i\dx
\Bigg|
\nonumber\\
&\ge
|\alpha_i-\alpha_i^p|
\Bigg|
\int_{\Rn}
\frac{
V_i^{p}-(1\pm\mathcal{o}(1))\theta V_i^{p-1}
}{|x|^t}
\,\Psi\Phi_i\dx
\Bigg|
-
(1+\mathcal{o}(1))p\alpha_i^{p-1}
\Bigg|
\int_{\Rn}
\frac{\widetilde{\sigma}_{<i}V_i^{p-1}}{|x|^t}
\,\Psi\Phi_i\dx
\Bigg|.
\end{align}

Combining \eqref{est-1} and \eqref{est-2} in \eqref{est0} and thanks to $\alpha_i<\tfrac{3}{2}$, we infer that
\begin{equation}\label{Est-10}
|\alpha_i-\alpha_i^p|
\left|
\int_{\Rn}
\frac{V_i^{p}-(1\pm\mathcal{o}(1))\theta V_i^{p-1}}{|x|^t}
\,\Psi\Phi_i\dx
\right|
\lesssim_{N,\nu,s,t}
\mathcal{o}(1)\,\|\rho\|_{\hs}
+\Gamma(u)
+\|\rho\|_{\hs}^{\min(2,p)}.
\end{equation}

       To estimate the left-hand side of \eqref{Est-10}, we apply the triangle inequality to obtain
\begin{align}\label{0est}
&\Bigg|\int_{\Rn}
\frac{V_i^{p}-(1\pm\mathcal{o}(1))\theta V_i^{p-1}}{|x|^t}
\,\Psi\Phi_i\dx
\Bigg|
\nonumber\\
&\ge
\Bigg|
\int_{\Rn}\frac{V_i^{p}\Psi}{|x|^t}\dx
-\theta(1\pm\mathcal{o}(1))
\int_{\Rn}\frac{V_i^{p-1}\Psi}{|x|^t}\dx
\Bigg|
\nonumber\\
&\quad
-\int_{\Rn}\frac{V_i^{p}|\Psi|(1-\Phi_i)}{|x|^t}\dx
-|\theta|(1\pm\mathcal{o}(1))
\int_{\Rn}\frac{V_i^{p-1}|\Psi|(1-\Phi_i)}{|x|^t}\dx
\nonumber\\
&=
\Bigg|
\int_{\Rn}\frac{V_i^{p}\Psi}{|x|^t}\dx
-\theta(1\pm\mathcal{o}(1))
\int_{\Rn}\frac{V_i^{p-1}\Psi}{|x|^t}\dx
\Bigg|
-\mathcal{o}(1)\eqqcolon|A[\Psi]-(1\pm\mathcal{o}(1))B[\Psi]\theta|-\mathcal{o}(1).
\end{align}

The last line in \eqref{0est} follows from Lemma~\ref{bumps}-(1) and Remark~\ref{remark on eigenfunctions}-$(2)$, which together imply that the contributions supported on the set $\{\Phi_i<1\}$ are negligible.  Indeed, one has
\begin{align*}
\int_{\Rn}\frac{V_i^{p}|\Psi|(1-\Phi_i)}{|x|^t}\dx
&\lesssim_{N,s,t}
\int_{\{\Phi_i<1\}}\frac{V_i^{p+1}}{|x|^t}\dx
= \mathcal{o}(1),
\\
\int_{\Rn}\frac{V_i^{p-1}|\Psi|(1-\Phi_i)}{|x|^t}\dx
&\lesssim_{N,s,t}
\int_{\{\Phi_i<1\}}\frac{V_i^{p}}{|x|^t}\dx
= \mathcal{o}(1).
\end{align*}

To justify the last estimate, we appeal directly to the definition~\eqref{def cutoff main} of the cutoff function
$\Phi_i$. For $\varepsilon<1$, this yields
\[
\int_{\{\Phi_i<1\}}\frac{V_i^{p}}{|x|^t}\dx
\le
\int_{|x|>\varepsilon^{-1}}\frac{V_i^{p}}{|x|^t}\dx
+\nu\int_{|x|<\varepsilon}\frac{V_i^{p}}{|x|^t}\dx
\lesssim_{N,s,t}
\varepsilon^{2s-t}+\nu\varepsilon^{N-t}
=\mathcal{o}(1),
\]
where the final bound follows from the computation in \eqref{lemma 2.3-(1) for p}.

Notice that \eqref{Est-10} closes the induction step at level $j=i$, and therefore yields
\eqref{Est-4} for all $1\le i\le\nu$, provided that the quantity
\begin{equation*}
A[\Psi]-(1\pm\mathcal{o}(1))\,B[\Psi]\theta
\end{equation*}
does not become too small simultaneously for the two choices $\Psi=V_i$ and
$\Psi=\dot V_i$.  

To verify this, we argue by contradiction. Suppose that
\[
A[\Psi]-(1\pm\mathcal{o}(1))\,B[\Psi]\theta=\pm\mathcal{o}(1)
\qquad\text{for all }\Psi\in\{V_i,\dot V_i\}.
\]

If $B[\Psi]\neq0$, this implies
\[
\theta=\frac{A[\Psi]}{B[\Psi]}\pm\mathcal{o}(1).
\]

Since $\theta\neq0$ is a fixed quantity defined above, we can then conclude
\[
\frac{A[V_i]}{B[V_i]}=\frac{A[\dot V_i]}{B[\dot V_i]}.
\]

We show below that this equality cannot occur, which yields the desired contradiction. Observe that, by \eqref{dilation map} and \eqref{bubble normalization}, the quantities
$A[\Psi]$ and $B[\Psi]$ can be computed explicitly for $\Psi\in\{V_i,\dot V_i\}$. 

Indeed,
\[
A[V_i]
= \int_{\Rn}\frac{V_i^{p+1}}{|x|^t}\dx
= \mu_{s,t}^{\frac{\twostar}{\twostar-2}},
\]
while
\[
A[\dot V_i]
= \int_{\Rn}\frac{V_i^{p}\dot V_i}{|x|^t}\dx
= \frac{1}{p+1}\,
\partial_{\lambda}
\int_{\Rn}\frac{(V_{s,t}^{\lambda})^{p+1}}{|x|^t}\dx\Big|_{\lambda=\lambda_i=1}
= 0.
\]

Moreover,
\[
B[V_i]
= \int_{\Rn}\frac{V_i^{p}}{|x|^t}\dx \neq 0,
\]
and
\[
B[\dot V_i]
= \int_{\Rn}\frac{V_i^{p-1}\dot V_i}{|x|^t}\dx
= \frac{1}{p}\,
\partial_{\lambda}
\int_{\Rn}\frac{V_{s,t}^{\lambda}(x)^{p}}{|x|^t}\dx\Big|_{\lambda=\lambda_i=1}
= -\frac{N-2s}{2p}
\int_{\Rn}\frac{V_{s,t}(x)^{p}}{|x|^t}\dx \neq 0.
\]

Indeed, by \cite[Theorem~1.1]{MN21}, the bubble satisfies the pointwise asymptotics $V_{s,t}(x)\sim (1+|x|^{N-2s})^{-1}$, so that $V_{s,t}(x)\sim 1$ for $|x|\le1$ and $V_{s,t}(x)\sim |x|^{-(N-2s)}$ for $|x|>1$, and hence
\begin{align}\label{lemma 2.3-(1) for p}
\int_{\Rn}\frac{V_{s,t}(x)^p}{|x|^t}\dx
&=\int_{|x|\le1}\frac{V_{s,t}(x)^p}{|x|^t}\dx
   +\int_{|x|>1}\frac{V_{s,t}(x)^p}{|x|^t}\dx\nonumber\\&
\sim_{N,s,t}\int_0^1 r^{N-1-t}\,{\rm d}r
   +\int_1^\infty r^{-1-(2s-t)}\,{\rm d}r=\tfrac{N+2s-2t}{(N-t)(2s-t)}.\end{align}

Thus, for each $\Psi\in\{V_i,\dot V_i\}$ one has $0<|B[\Psi]|<\infty$. Hence, the above computations yield
\[
\frac{A[V_i]}{B[V_i]}\neq 0=\frac{A[\dot V_i]}{B[\dot V_i]},
\] affirming the desired contradiction. Therefore, by choosing an appropriate $\Psi$ in \eqref{Est-10}, we obtain \eqref{Est-4}.

To conclude the proof of \eqref{Est-5}, we set $\Psi=V_i$ in \eqref{est-1}. By the triangle inequality this yields
\begin{align*}
&\Bigg|
\int_{\Rn}
\frac{\left(p\alpha_i^{p-1}+\mathcal{o}(1)\right)V_i^{\,p}\widetilde{\sigma}_i}{|x|^{t}}
\,\Phi_i\dx
\Bigg|
\\
&\qquad\lesssim_{N,s,t}
\mathcal{o}(1)\,\|\rho\|_{\hs}
+\Gamma(u)
+\|\rho\|_{\hs}^{\min(2,p)}
+|\alpha_i-\alpha_i^p|\,\|V_i\|_{L^{p+1}_{|x|^{-t}}}^{p+1}.
\end{align*}

Applying \eqref{Est-4}, we therefore obtain for every $j\neq i$,
\begin{equation}\label{Est100}
\int_{\Rn}\frac{V_i^{p}\widetilde{\sigma}_i\,\Phi_i}{|x|^t}\dx
\lesssim_{N,\nu,s,t}
\mathcal{o}(1)\,\|\rho\|_{\hs}
+\Gamma(u)
+\|\rho\|_{\hs}^{\min(2,p)}.
\end{equation}

We now remove the cutoff. By the induction hypothesis, the required estimate~\eqref{Est-5} is already available for $j<i$. It therefore remains to consider $j>i$. Under the normalization $\lambda_i=1$, the $\delta$-interaction condition implies $\lambda_j\le\delta=\mathcal{o}(1)$, and therefore the continuity of $V_{s,t}$ at the origin imply that
\[
V_j(x)\sim_{N,s,t}V_j(0),
\qquad |x|\le1,
\]
uniformly for $j>i$. Moreover, since $V_i\ge V_i(1)$ on $\overline{B_1}$, property~(1) of Lemma~\ref{bumps} yields
\[
\int_{|x|\le1}
\frac{V_i^pV_j(1-\Phi_i)}{|x|^t}\dx\lesssim_{N,s,t}
\frac{V_j(0)}{V_i(1)}\int_{\{\Phi_i<1\}}
\frac{V_i^{p+1}}{|x|^t}\dx
\lesssim_{N,s,t}\varepsilon V_j(0).
\]

On the other hand, the same uniform comparison gives
\[
\int_{|x|\le1}\frac{V_i^pV_j}{|x|^t}\dx\gtrsim_{N,s,t}
V_j(0)\int_{|x|\le1}\frac{V_i^p}{|x|^t}\dx
\gtrsim_{N,s,t}V_j(0).
\]

Consequently,
\[
\int_{|x|\le1}
\frac{V_i^pV_j(1-\Phi_i)}{|x|^t}\dx
\lesssim_{N,s,t}\varepsilon \int_{|x|\le1}\frac{V_i^pV_j}{|x|^t}\dx.
\]

Since $0\le\Phi_i\le1$, we have
\begin{align*}
\int_{|x|\le1}\frac{V_i^pV_j}{|x|^t}\dx
&=\int_{|x|\le1}
\frac{V_i^pV_j\Phi_i}{|x|^t}\dx+
\int_{|x|\le1}
\frac{V_i^pV_j(1-\Phi_i)}{|x|^t}\dx\\&\le
\int_{\Rn}
\frac{V_i^pV_j\Phi_i}{|x|^t}\dx
+C\varepsilon \int_{|x|\le1}\frac{V_i^pV_j}{|x|^t}\dx.
\end{align*}

Choosing $\varepsilon$ sufficiently small and absorbing the last term, we infer
\[
\int_{|x|\le1}\frac{V_i^pV_j}{|x|^t}\dx
\lesssim_{N,s,t}
\int_{\mathbb R^N}
\frac{V_i^pV_j\Phi_i}{|x|^t}\dx.
\]

Since $\alpha_j>\frac12$ and
$\widetilde{\sigma}_i=\sum_{k\neq i}\alpha_kV_k$,
estimate~\eqref{Est100} now gives, for every $j>i$,
\begin{equation}\label{Est101}
\int_{|x|\le1}\frac{V_i^pV_j}{|x|^t}\dx
\lesssim_{N,\nu,s,t}
\mathcal{o}(1)\|\rho\|_{\hs}
+\Gamma(u)
+\|\rho\|_{\hs}^{\min(2,p)}.
\end{equation}

Since $\lambda_i=1$, Corollary~\ref{loc-bub-bub-est} now entails for indices $j\neq i$,
\begin{equation}\label{Est102}
	\int_{\Rn}\frac{V_i^{p}V_j}{|x|^t}\dx
	\sim_{N,s,t}
	\int_{|x|\le1}\frac{V_i^{p}V_j}{|x|^t}\dx.
\end{equation}

In view of \eqref{Est102} and \eqref{Est101}, estimate \eqref{Est-5} follows. This completes the proof.\qed

       \subsection{Proof of Theorem~\ref{energybound-implies-stability}}
\label{S2.5}
We argue by contradiction. Let $\delta_0>0$ and $C_0>0$ be the constants in Theorem~\ref{quantitative-stability}. If the claimed estimate were false, then for every $k\ge1$ there would exist a nonnegative function $u_k\in\Hs$ satisfying~\eqref{energy-bdd} such that
\[
\Big\|u_k-\sum_{i=1}^{\nu}V_{s,t}^{\lambda_i}\Big\|_{\hs}
>k\,\Gamma(u_k)
\]
for every choice of $\lambda_1,\ldots,\lambda_\nu>0$.

Taking $\lambda_1=\cdots=\lambda_\nu=1$ and using the energy bound~\eqref{energy-bdd}, we infer that
\[
k\,\Gamma(u_k)
<
\left\|u_k-\nu V_{s,t}^{1}\right\|_{\hs}
\lesssim_{N,\nu,s,t}1.
\]

In particular, $\Gamma(u_k)\to0$ as $k\to\infty$. By Theorem~\ref{SDHS}, after passing to a subsequence, there exist positive scales
$\lambda_k^1,\ldots,\lambda_k^\nu$ such that
\[
\Big\|u_k-\sum_{i=1}^{\nu}
V_{s,t}^{\lambda_k^i}\Big\|_{\hs}
\longrightarrow0
\qquad\text{as }k\to\infty.
\]

Moreover, the scale-separation property in Theorem~\ref{SDHS} gives $\min\Big\{\frac{\lambda_k^i}{\lambda_k^j},\frac{\lambda_k^j}{\lambda_k^i}\Big\}\to0$ for $i\neq j$. Consequently, for all sufficiently large $k$, the family
$\bigl(V_{s,t}^{\lambda_k^i}\bigr)_{i=1}^{\nu}$
is $\delta_0$-interacting and
\[
\Big\|u_k-\sum_{i=1}^{\nu}
V_{s,t}^{\lambda_k^i}\Big\|_{\hs}
\le\delta_0.
\]

All the hypotheses of Theorem~\ref{quantitative-stability} are therefore satisfied. Applying that theorem, we obtain positive scales $\widetilde\lambda_k^1,\ldots,\widetilde\lambda_k^\nu$ such that
\[
\Big\|u_k-\sum_{i=1}^{\nu}
V_{s,t}^{\widetilde\lambda_k^i}\Big\|_{\hs}
\le C_0\,\Gamma(u_k).
\]

This contradicts the choice of $u_k$ whenever $k>C_0$. Hence there exist parameters $\lambda_1,\ldots,\lambda_\nu>0$ and a constant $C(N,\nu,s,t)>0$ such that
\[
u=\sum_{i=1}^{\nu}V_{s,t}^{\lambda_i}+\rho,
\qquad
\|\rho\|_{\hs}\le C(N,\nu,s,t)\Gamma(u),
\]
as stated in~\eqref{def-bds-dis-lin}.

Finally, the $\hs$-interaction estimate for any two distinct bubbles \eqref{Hs-interaction} follows directly from \eqref{Est-5} and \eqref{def-bds-dis-lin}.\qed

\subsubsection{Sharpness in low dimensions}
We work in the low-dimensional regime $2s<N<6s-2t$, equivalently $p=\twostar-1>2$. To establish the optimality of the linear rate in~\eqref{def-bds-dis-lin}, we argue by contradiction. Let $\mathcal M$ denote the multi-bubble set introduced above, and set
\[
d_{\mathcal M}(u)\coloneqq
\inf_{\tau\in\mathcal M}\|u-\tau\|_{\hs}.
\]

Suppose that there exists a fixed $0<\varepsilon<1$ such that
\begin{equation}\label{S-ass}
d_{\mathcal M}(u)^{1-\varepsilon}
\lesssim_{N,\nu,s,t}
\Gamma(u)
\end{equation}
for every nonnegative $u\in\Hs$ satisfying the energy bound~\eqref{energy-bdd}.

Fix a nonnegative, nonradial function $\phi\in\cc$, with $\|\phi\|_{\hs}>0$. For each $0<\kappa\ll1$, put
\[
\delta_\kappa\coloneqq
\kappa^{\frac{2}{N-2s}},
\qquad
\lambda_i^\kappa\coloneqq
\delta_\kappa^{i-1},
\qquad 1\le i\le\nu.
\]

Thus, for every $i\neq j$,
\[
Q_{ij}^\kappa
\coloneqq
\min\bigg\{
\frac{\lambda_i^\kappa}{\lambda_j^\kappa},
\frac{\lambda_j^\kappa}{\lambda_i^\kappa}
\bigg\}
=
\delta_\kappa^{|i-j|}
\le\delta_\kappa.
\]

In particular, the associated family of bubbles is $\delta_\kappa$-interacting.

Set
\[
V_i\coloneqq V_{s,t}^{\lambda_i^\kappa},
\qquad
\sigma\coloneqq\sum_{i=1}^{\nu}V_i,
\qquad
u_\kappa\coloneqq\sigma+\kappa\phi.
\]

Throughout the following estimates, the dependence of
$V_i$ and $\sigma$ on $\kappa$ is suppressed from the
notation.

Recall from Appendix~\ref{AC}, the closed linear subspace $\hs_{\mathrm{rad}}(\Rn)$ of $\Hs$, and then $\mathcal M\subset\hs_{\mathrm{rad}}(\Rn)$, as all the bubbles are radial. Since $\phi$ is nonradial,
\[
C_\phi\coloneqq
\inf_{v\in\hs_{\mathrm{rad}}}
\|\phi-v\|_{\hs}>0.
\]

Consequently, using $\sigma\in\mathcal M$, we obtain
\[
d_{\mathcal M}(u_\kappa)\le\|u_\kappa-\sigma\|_{\hs}=\kappa\|\phi\|_{\hs},\]
\[d_{\mathcal M}(u_\kappa)
\ge
\inf_{v\in\hs_{\mathrm{rad}}}
\|u_\kappa-v\|_{\hs}=\kappa
\inf_{v\in\hs_{\mathrm{rad}}}
\|\phi-v\|_{\hs}
=C_\phi\kappa.
\]

Thus
\begin{equation}\label{S-distance}
d_{\mathcal M}(u_\kappa)\sim\kappa.
\end{equation}

Using~\eqref{bubble normalization} and Corollary~\ref{general-bubble-interaction-estimate}, one has
\[
\|V_i\|_{\hs}^2
=\mu_{s,t}^{\frac{N-t}{2s-t}},
\qquad
\langle V_i,V_j\rangle_{\hs}
\lesssim_{N,s,t}
(Q_{ij}^\kappa)^{\frac{N-2s}{2}}
\]
whenever $i\neq j$. Thus,
\[
\|\sigma\|_{\hs}^2
=
\nu\,\mu_{s,t}^{\frac{N-t}{2s-t}}
+
\mathcal{O}\left(
\delta_\kappa^{\frac{N-2s}{2}}
\right),
\]
and hence
\[
\|u_\kappa\|_{\hs}^2
=
\|\sigma\|_{\hs}^2
+\kappa^2\|\phi\|^2_{\hs}
+2\kappa\langle\sigma,\phi\rangle_{\hs}
=
\nu\mu_{s,t}^{\frac{N-t}{2s-t}}
+\mathcal{O}\left(
\kappa+\kappa^2+
\delta_\kappa^{\frac{N-2s}{2}}
\right).
\]

Since $\delta_\kappa^{\frac{N-2s}{2}}=\kappa$ and $\|\sigma\|_{\hs}$ is uniformly bounded, $u_\kappa$ is nonnegative and satisfies~\eqref{energy-bdd} for $\kappa$ sufficiently small.

Recalling that
\[
\Gamma(u)=
\Bigl\|(-\Delta)^s u-\frac{u^p}{|x|^t}\Bigr\|_{\fhs},
\qquad p+1=\twostar,
\]
we decompose
\[
(-\Delta)^s u_\kappa-\frac{u_\kappa^p}{|x|^t}
=A_\kappa+B_\kappa+C_\kappa,
\]
where
\[
A_\kappa\coloneqq\kappa(-\Delta)^s\phi,\qquad
B_\kappa\coloneqq
\frac{\sigma^p-(\sigma+\kappa\phi)^p}{|x|^t},
\qquad
C_\kappa\coloneqq
\sum_{i=1}^{\nu}\frac{V_i^p}{|x|^t}
-\frac{\sigma^p}{|x|^t}.
\]

\medskip
\noindent\emph{Estimate for $A_\kappa$.}
By the definition of the operator norm on $\Fhs$,
\begin{equation}\label{S20}
\|A_\kappa\|_{\fhs}
=\sup_{\|\psi\|_{\hs}=1}
\kappa|\langle\phi,\psi\rangle_{\hs}|
\le\kappa\,\|\phi\|_{\hs}
\lesssim\kappa .
\end{equation}

\noindent\emph{Estimate for $B_\kappa$.}
Using the pointwise inequality~\eqref{Ineq-1},
\[
|(\sigma+\kappa\phi)^p-\sigma^p|
\le
p\sigma^{p-1}\kappa\phi
+C(p)(\sigma^{p-2}\kappa^2\phi^2+\kappa^p\phi^p),
\]
we obtain
\[
|B_\kappa|
\lesssim
\kappa\frac{\sigma^{p-1}\phi}{|x|^t}
+\kappa^2\frac{\sigma^{p-2}\phi^2}{|x|^t}
+\kappa^p\frac{\phi^p}{|x|^t}.
\]

By H\"older's inequality and the Hardy--Sobolev embedding~\eqref{fhs0}, we have for
$\|\psi\|_{\hs}=1$,
\begin{align*}
\Bigl|
\int_{\Rn}\frac{\sigma^{p-1}\phi\,\psi}{|x|^t}\dx
\Bigr|
&\le
\|\sigma\|_{L^{p+1}_{|x|^{-t}}}^{p-1}
\|\phi\|_{L^{p+1}_{|x|^{-t}}}
\|\psi\|_{L^{p+1}_{|x|^{-t}}}\\
&\lesssim
\|\sigma\|_{\hs}^{p-1}
\|\phi\|_{\hs}\|\psi\|_{\hs}
\lesssim 1 .
\end{align*}

Here we used the uniform bound $\|\sigma\|_{\hs}\lesssim\nu\|V_{s,t}\|_{\hs}$. The terms involving $\kappa^2$ and $\kappa^p$ are treated analogously, yielding for $\kappa\ll1$,
\begin{equation}\label{S3}
\|B_\kappa\|_{\fhs}\lesssim\kappa .
\end{equation}

\noindent\emph{Estimate for $C_\kappa$.}
Using~\eqref{Ineq-2},
\[
\Bigl|\sigma^p-\sum_{i=1}^{\nu}V_i^p\Bigr|
\lesssim
\sum_{i\ne j}V_i^{p-1}V_j ,
\]
we have, for $\psi\in\hs(\Rn)$ with $\|\psi\|_{\hs}=1$,
\[
\Bigl|{}_{\fhs}\langle C_\kappa,\psi\rangle_{\hs}\Bigr|
\lesssim
\sum_{i\ne j}
\int_{\Rn}\frac{V_i^{p-1}V_j|\psi|}{|x|^t}.
\]

Applying H\"older's inequality and the Hardy--Sobolev inequality~\eqref{fhs0} yields
\[
\int_{\Rn}\frac{V_i^{p-1}V_j|\psi|}{|x|^t}
\lesssim
\|V_i^{p-1}V_j\|_{L^{\frac{p+1}{p}}_{|x|^{-t}}}.
\]

Therefore, by Corollary~\ref{general-bubble-interaction-estimate} and the choice
$\delta_\kappa=\kappa^{\frac{2}{N-2s}}$,
\begin{equation}\label{S4}
\|C_\kappa\|_{\fhs}
\lesssim
\sum_{i\ne j}(Q_{ij}^\kappa)^{\frac{N-2s}{2}}
\lesssim
\delta_\kappa^{\frac{N-2s}{2}}
=\kappa.
\end{equation}

Combining~\eqref{S20}, \eqref{S3}, and~\eqref{S4}, we obtain
\begin{equation}\label{S5}
\Gamma(u_\kappa)\lesssim\kappa .
\end{equation}

Finally, invoking~\eqref{S-ass}, \eqref{S-distance}, and~\eqref{S5}, we obtain
\[
\kappa^{1-\varepsilon}
\lesssim
d_{\mathcal M}(u_\kappa)^{1-\varepsilon}
\lesssim
\Gamma(u_\kappa)
\lesssim\kappa,
\]
which is impossible as $\kappa\to0$. This contradiction shows that the linear rate in \eqref{def-bds-dis-lin} is sharp.

\appendix   
	\section{Few useful estimates}\label{AA}
	We seek to establish a Kato--Ponce type estimate for the fractional commutator $\mathcal{C}[f,g]$ (as defined in Theorem~\ref{wkpv-estimate}) in $L^p_w$ norms with weights $w$ belonging to the Muckenhoupt class $\mathcal{A}_p$, as outlined in the Introduction. The unweighted, one-dimensional case (also known as the fractional Leibniz rule) was first established by Kenig, Ponce, and Vega in the seminal work~\cite[Theorem~A.8]{KPV93}. Since our setting requires a weighted, higher-dimensional version, we outline a proof in Theorem~\ref{wkpv-estimate} below for completeness. While both approaches are grounded in Littlewood--Paley theory, our method circumvents the use of Coifman--Meyer paradifferential calculus by directly applying the LP theory to an alternative paraproduct decomposition~\eqref{paraproduct}.

\medskip
We begin by introducing some notations and a definition. \begin{Definition}[Muckenhoupt class]\label{muckenhoupt}
Let $w\in L^1_{\rm loc}(\Rn)$, $w\ge0$. For $p>1$, define the
$\mathcal{A}_p$-characteristic of $w$ by
\[
[w]_p\coloneqq
\sup_{B\text{ is a ball in }\Rn}
\left(\fint_B w(x)\dx\right)
\left(\fint_B w(x)^{\frac1{1-p}}\dx\right)^{p-1}.
\]
For $p=1$, define
\[
[w]_1\coloneqq
\Big\|\tfrac{Mw}{w}\Big\|_{L^\infty(\Rn)},
\]
where $M$ denotes the Hardy--Littlewood maximal function, i.e.,
for any $f\in L^1_{\rm loc}(\Rn)$,
\[
Mf(x)\coloneqq\sup_{B\ni x}\fint_B|f(y)|\dy.
\]
We say that $w\in\mathcal{A}_p$, for $p\ge1$, if and only if $[w]_p<\infty$.
\end{Definition}

    Define $\psi\in\cc$ to be a radial function satisfying $\psi(\xi)=1$ for $|\xi|\le1$, $\psi(\xi)=0$ for $|\xi|\ge 2$ and $0\le\psi\le1$.
	Let $\varphi(\xi) \coloneqq \psi(\xi) - \psi(2\xi)$, so that $\supp\varphi=\left\{ \xi \in \Rn : \frac{1}{2} \le |\xi| \le 2 \right\}$. For any $f \in \mathcal{S}(\Rn)$ and $j \in \mathbb{Z}$, define the dyadic rescaling $$
	\varphi_j(\xi) \coloneqq \varphi\left( \frac{\xi}{2^j} \right), \qquad \psi_j(\xi) \coloneqq \psi\left( \frac{\xi}{2^j} \right),
	$$ and introduce the Littlewood--Paley projection operators $P_j$ and $P_{\le j}$ via $$
	\Widehat{P_j f}(\xi) \coloneqq \varphi_j(\xi)\, \Widehat{f}(\xi), \qquad \Widehat{P_{\le j} f}(\xi) \coloneqq \psi_j(\xi)\, \Widehat{f}(\xi).
	$$
	
	We further define $P_{>j} \coloneqq \id - P_{\le j}$. When the context is clear, we will denote $f_j \coloneqq P_j f$, $f_{\le j} \coloneqq P_{\le j} f$, and so on for notational brevity.
	
	\begin{Lemma}\label{derivative localization}
	If $s>0,\,1<p<\infty$ and $w\in\mathcal{A}_p$, then for all $h\in\mathcal{S}(\Rn)$ the following holds. \[\Bigg\|\sqrt{\sum_{j\in\Z}\left|2^{-js}\flap h_{\le j}\right|^2}\Bigg\|_{L^p_w(\Rn)}\lesssim_{N,p,s,[w]_p}\left\|h\right\|_{L^p_w(\Rn)}.\]\end{Lemma}
    \Proof Lets consider the fattened LP projector $Q_k,\,k\in\Z$ such that $Q_kh=\Widehat{\chi_k}*h$ where $\chi_k(\xi)=\varphi(\frac\xi{2^{k+1}})+\varphi(\frac\xi{2^k})+\varphi(\frac\xi{2^{k-1}})$. Then it is easy to see that \begin{align*}
        \left|2^{-js}\flap h_{\le j}\right|\le\sum_{k\le j}2^{-js}\left|\flap h_k\right|=\sum_{k\le j}2^{-js}\left|\flap Q_kh_k\right|\lesssim_{N,s}\sum_{k\le j}2^{(k-j)s}Mh_k.
    \end{align*}

    Observe by Minkowski inequality, \begin{align*}\Bigg\|\sqrt{\sum_{j\in\Z}\Big|\sum_{k\le j}2^{(k-j)s}Mh_k\Big|^2}\Bigg\|_{L^p_w}=\Bigg\|\sqrt{\sum_{j\in\Z}\Big|\sum_{i\ge0}2^{-is}Mh_{j-i}\Big|^2}\Bigg\|_{L^p_w}\lesssim_{s}\Bigg\|\sqrt{\sum_{k\in\Z}|Mh_k|^2}\Bigg\|_{L^p_w}.
    \end{align*}

    Now an application of Fefferman--Stein inequality~\cite[Lemma~A.5/A.3(e)]{KPV93} followed by the square function estimate~\cite[Lemma~A.5/A.3(a)]{KPV93}, entails \begin{align*}\Bigg\|\sqrt{\sum_{k\in\Z}|Mh_k|^2}\Bigg\|_{L^p_w}\lesssim_{N,p,[w]_p}\Bigg\|\sqrt{\sum_{k\in\Z}|h_k|^2}\Bigg\|_{L^p_w}\lesssim_{N,p,[w]_p}\left\|h\right\|_{L^p_w}.
    \end{align*} 
    
    Combining, the three previous inequalities, the result follows.\qed
	\begin{theorem}[Weighted Kenig-Ponce-Vega Estimate]\label{wkpv-estimate}
		Let $\alpha \in (0,\frac12)$ and $\alpha_1, \alpha_2 \in [0,\alpha]$ with $\alpha_1 + \alpha_2 = \alpha$. Suppose that $p, p_1, p_2 \in (1, \infty)$ satisfy the H\"older relation $\tfrac{1}{p} = \tfrac{1}{p_1} + \tfrac{1}{p_2}$. Let $w_1, w_2 \in L^1_{\mathrm{loc}}(\Rn)$ be positive weight functions such that $w_1 \in \mathcal{A}_{p_1}$, $w_2 \in \mathcal{A}_{p_2}$, and $w = w_1^{p/p_1} w_2^{p/p_2}\in\mathcal{A}_p$.
		
		Define the commutator operator $$\mathcal{C}[f,g] \coloneqq\alap(fg) - g\alap f - f\alap g.$$
		
		Then, for all $f, g \in \mathcal{S}(\Rn)$, the following estimate holds. $$\left\lVert \mathcal{C}[f,g] \right\rVert_{L^p_w(\Rn)} \lesssim_{N,p,p_1,p_2,\alpha,\alpha_1,\alpha_2,[w]_p,[w_1]_{p_1},[w_2]_{p_2}} \left\lVert (-\Delta)^{\alpha_1}f \right\rVert_{L^{p_1}_{w_1}(\Rn)} \left\lVert (-\Delta)^{\alpha_2}g \right\rVert_{L^{p_2}_{w_2}(\Rn)}.
		$$\end{theorem}
	\Proof We begin with the Bony paraproduct decomposition~\cite{B81}; see also \cite[Lemma~2.7]{Tao07}. \begin{equation}\label{paraproduct}
		fg=\sum_{j\in\Z}f_j\widetilde{g}_j+\sum_{j\in\Z}f_{\le j-3}g_j+\sum_{j\in\Z}f_jg_{\le j-3}.
	\end{equation} where $\widetilde{g}_j=\sum_{|j-k|\le2}g_k$.
	
	\noindent\textbf{Step 1.} First we deal with the diagonal (high-high frequency) piece $\sum f_j\widetilde{g}_j$.
	
	Using Littlewood--Paley square function estimates~\cite[Lemma~A.5/A.3(a,b)]{KPV93} for the fractional operator $\alap$, one obtains \[\Bigg\|\alap\sum_{j\in\Z}f_j\widetilde{g}_j\Bigg\|_{L^p_w}\sim_{N,p,\alpha,[w]_p}\Bigg\|\sqrt{\sum_{k\in\Z}\Big|2^{2k\alpha}P_k\sum_{j\in\Z}f_j\widetilde{g}_j\Big|^2}\Bigg\|_{L^p_w}.\]
	
	Elementary support properties of the projectors $P_j$ then give $P_k(f_j\widetilde{g}_j)=0$ for all $j\le k-5$ and thus $P_k\sum_{j\in\Z}f_j\widetilde{g}_j=P_k\sum_{l>-5}f_{k+l}\widetilde{g}_{k+l}$. Denoting, $h_{k}=\sum_{l>-5}2^{2k\alpha}f_{k+l}\widetilde{g}_{k+l}$, and using the pointwise Bernstein's inequality $|P_jh_k(x)|\lesssim Mh_k(x);\,\forall j\in\Z$, one sees \[\Bigg\|\sqrt{\sum_{k\in\Z}\Big|2^{2k\alpha}P_k\sum_{j\in\Z}f_j\widetilde{g}_j\Big|^2}\Bigg\|_{L^p_w}=\Bigg\|\sqrt{\sum_{k\in\Z}\left|P_kh_k\right|^2}\Bigg\|_{L^p_w}\lesssim_{N}\Bigg\|\sqrt{\sum_{k\in\Z}\left|Mh_k\right|^2}\Bigg\|_{L^p_w}.\]
	
	Invoking the Fefferman--Stein inequality~\cite[Lemma~A.5/A.3(e)]{KPV93}, we achieve \begin{align*}
		\Bigg\|\sqrt{\sum_{k\in\Z}\left|Mh_{k}\right|^2}\Bigg\|_{L^p_w}&\lesssim_{N,p,[w]_p}\Bigg\|\sqrt{\sum_{k\in\Z}|h_{k}|^2}\Bigg\|_{L^p_w}.
	\end{align*}
	
	Now an application of Minkowski's inequality entails \begin{align*}\Bigg\|\sqrt{\sum_{k\in\Z}|h_{k}|^2}\Bigg\|_{L^p_w}=\Bigg\|\sqrt{\sum_{k\in\Z}\Big|\sum_{l>-5}2^{2k\alpha}f_{k+l}\widetilde{g}_{k+l}\Big|^2}\Bigg\|_{L^p_w}&\le\sum_{l>-5}\frac{1}{2^{2l\alpha}}\Bigg\|\sqrt{\sum_{k\in\Z}\left|2^{2(k+l)\alpha}f_{k+l}\widetilde{g}_{k+l}\right|^2}\Bigg\|_{L^p_w}\\&=\Bigg\|\sqrt{\sum_{k\in\Z}\left|2^{2k\alpha}f_{k}\widetilde{g}_{k}\right|^2}\Bigg\|_{L^p_w}\sum_{l>-5}\frac{1}{2^{2l\alpha}}\\&\sim_{\alpha}\Bigg\|\sqrt{\sum_{k\in\Z}\left|2^{2k\alpha_1}f_{k}\cdot2^{2k\alpha_2}\widetilde{g}_{k}\right|^2}\Bigg\|_{L^p_w}.\end{align*}
	
	Finally, weighted H\"older inequality \big(with exponent $p=\frac{p_1p_2}{p_1+p_2}$ and weight $w=w_1^{p/p_1}w_2^{p/p_2}$\big) and the square function estimates together imply \begin{align*}\Bigg\|\sqrt{\sum_{k\in\Z}\left|2^{2k\alpha_1}f_{k}\cdot2^{2k\alpha_2}\widetilde{g}_{k}\right|^2}\Bigg\|_{L^p_w}&\le\Bigg\|\sqrt{\sum_{k\in\Z}\left|2^{2k\alpha_1}f_{k}\right|^2}\sqrt{\sum_{k\in\Z}\left|2^{2k\alpha_2}\widetilde{g}_{k}\right|^2}\Bigg\|_{L^p_w}\\&\le\Bigg\|\sqrt{\sum_{k\in\Z}\left|2^{2k\alpha_1}f_{k}\right|^2}\Bigg\|_{L^{p_1}_{w_1}}\Bigg\|\sqrt{\sum_{k\in\Z}\left|2^{2k\alpha_2}\widetilde{g}_{k}\right|^2}\Bigg\|_{L^{p_2}_{w_2}}\\&\lesssim_{N,p_1,p_2,\alpha_1,\alpha_2,[w_1]_{p_1},[w_2]_{p_2}}\left\|(-\Delta)^{\alpha_1}f\right\|_{L^{p_1}_{w_1}}\left\|(-\Delta)^{\alpha_2}g\right\|_{L^{p_2}_{w_2}}.\end{align*}
	
	The last inequality for $\widetilde{g}_j$ follows by noticing that $\sum\left|2^{2j\alpha_2}\widetilde{g}_j\right|^2\lesssim_{\alpha_2}\sum|2^{2j\alpha_2}g_j|^2$. This settles the diagonal piece. Since the last two pieces in the decomposition~\eqref{paraproduct} are symmetric with respect to $f$ and $g$, it is enough to only deal with the low-high frequency piece $\sum f_{\le j-3}g_j$.
	
	\noindent\textbf{Step 2.} Start by writing \begin{align}\label{paraproduct1}\alap\sum_{j\in\Z}f_{\le j-3}g_j&=\sum_{j\in\Z}\left[
	\alap f_{\le j-3}g_j-f_{\le j-3}\alap g_j\right]\notag\\&\quad+f\alap g-\sum_{j\in\Z}f_{> j-3}\alap g_j.\end{align}
	
	First we take care of the last sum in \eqref{paraproduct1}.  Reversing the order of summation
	gives \begin{equation*}\sum_{j\in\mathbb Z}f_{>j-3}(-\Delta)^\alpha g_j
	=\sum_{j\in\mathbb Z}f_j(-\Delta)^\alpha g_{\le j+2}.\end{equation*}

	\noindent\textit{Case 1.} $\alpha_1>0$.
	Since $\alpha=\alpha_1+\alpha_2$, we have
	\begin{align*}
		\sum_{j\in\Z} f_j(-\Delta)^\alpha g_{\le j+2}
		=2^{4\alpha_1}\sum_{j\in\Z}2^{2j\alpha_1}f_j\cdot2^{-2(j+2)\alpha_1}
		(-\Delta)^{\alpha_1}
		((-\Delta)^{\alpha_2}g)_{\le j+2}.
	\end{align*}

	Hence, by the Cauchy--Schwarz inequality, weighted H\"older's
	inequality, and the Littlewood--Paley inequality, \begin{align*}
		\Bigg\|\sum_{j\in\Z} f_j(-\Delta)^\alpha g_{\le j+2}
		\Bigg\|_{L^p_w}\hspace{-3.5mm}&\lesssim
		\Bigg\|\sqrt{\sum_{j\in\Z}
		|2^{2j\alpha_1}f_j|^2}
		\Bigg\|_{L^{p_1}_{w_1}}
		\Bigg\|\sqrt{\sum_{j\in\Z}
		\left|2^{-2(j+2)\alpha_1}
		(-\Delta)^{\alpha_1}
		((-\Delta)^{\alpha_2}g)_{\le j+2}\right|^2}\Bigg\|_{L^{p_2}_{w_2}}\\&\lesssim
		\|(-\Delta)^{\alpha_1}f\|_{L^{p_1}_{w_1}}
		\Bigg\|\sqrt{\sum_{j\in\Z}
		\left|2^{-2j\alpha_1}
		(-\Delta)^{\alpha_1}
		((-\Delta)^{\alpha_2}g)_{\le j}\right|^2}\Bigg\|_{L^{p_2}_{w_2}}.
	\end{align*}
	
	Applying Lemma~\ref{derivative localization} with
	$s=2\alpha_1$ and $h=(-\Delta)^{\alpha_2}g$, yields
	\[
	\Bigg\|\sqrt{\sum_{j\in\Z}
	\left|2^{-2j\alpha_1}
	(-\Delta)^{\alpha_1}
	((-\Delta)^{\alpha_2}g)_{\le j}\right|^2}\Bigg\|_{L^{p_2}_{w_2}}
	\lesssim
	\|(-\Delta)^{\alpha_2}g\|_{L^{p_2}_{w_2}}.
	\]
	
	Consequently,
	\[\Bigg\|\sum_{j\in\Z}f_{>j-3}(-\Delta)^\alpha g_j\Bigg\|_{L^p_w}
		\lesssim
		\|(-\Delta)^{\alpha_1}f\|_{L^{p_1}_{w_1}}
		\|(-\Delta)^{\alpha_2}g\|_{L^{p_2}_{w_2}}.\]
	
	\noindent\textit{Case 2.} $\alpha_1=0$ (hence $\alpha_2=\alpha$).
	We write \begin{align*}
		\sum_{j\in\Z}f_j(-\Delta)^\alpha g_{\le j+2}=\sum_{j\in\Z}f_j(-\Delta)^\alpha g_{\le j-3}+\sum_{j\in\Z}f_j\sum_{j-3<k<j+3}(-\Delta)^\alpha g_k=T_1+T_2,
	\end{align*}

	For $T_1$, the definitions of $\varphi_j$ and $\psi_j$ give \[\supp\Widehat{f_j}
	\subseteq\{\xi:2^{j-1}\le|\xi|\le2^{j+1}\}, \text{ and } \supp
	\Widehat{(-\Delta)^\alpha g_{\le j-3}}
	\subseteq B(0,2^{j-2}).\]
	
	Therefore one obtains $\supp
	\Widehat{f_j(-\Delta)^\alpha g_{\le j-3}}
	\subseteq
	\left\{\xi:
	2^{j-2}\le|\xi|\le9\cdot2^{j-2}\right\}$ and hence these products have Fourier supports in annuli comparable
	to $2^j$, with uniformly finite overlap. The weighted square-function
	estimate, followed by the pointwise bound
	\[|(-\Delta)^\alpha g_{\le j-3}|=|P_{\le j-3}(-\Delta)^\alpha g|
	\lesssim M((-\Delta)^\alpha g),\quad\forall j\in\Z,\] yields \begin{align*}
		\|T_1\|_{L^p_w}
		&\lesssim
		\Bigg\|\sqrt{\sum_{j\in\Z}
		|f_j(-\Delta)^\alpha g_{\le j-3}|^2}\Bigg\|_{L^p_w}\\
		&\lesssim
		\Bigg\|\sqrt{\sum_{j\in\Z}|f_j|^2}\cdot
		M((-\Delta)^\alpha g)\Bigg\|_{L^p_w}\\
		&\le
		\Bigg\|\sqrt{\sum_{j\in\Z}|f_j|^2}
		\Bigg\|_{L^{p_1}_{w_1}}
		\Big\|M((-\Delta)^\alpha g)\Big\|_{L^{p_2}_{w_2}}\lesssim
		\|f\|_{L^{p_1}_{w_1}}
		\|(-\Delta)^\alpha g\|_{L^{p_2}_{w_2}}.
	\end{align*}
	
	For $T_2$, since $j-3<k<j+3$ is equivalent to
	$k=j+d$ with $d\in\{-2,-1,0,1,2\}$, pointwise
	Cauchy--Schwarz, followed by a weighted H\"older inequality and the Littlewood--Paley inequality gives \begin{align*}
		\|T_2\|_{L^p_w}
		&\le\sum_{d=-2}^{2}
		\Bigg\|\sum_{j\in\Z} f_j(-\Delta)^\alpha g_{j+d}
		\Bigg\|_{L^p_w}\\
		&\le\sum_{d=-2}^{2}\Bigg\|\sqrt{\sum_{j\in\Z}|f_j|^2}
		\sqrt{\sum_{j\in\Z}|(-\Delta)^\alpha g_{j+d}|^2}
		\Bigg\|_{L^p_w}\\
		&\le\sum_{d=-2}^{2}
		\Bigg\|\sqrt{\sum_{j\in\Z}|f_j|^2}
		\Bigg\|_{L^{p_1}_{w_1}}
		\Bigg\|\sqrt{\sum_{j\in\Z}|P_{j+d}(-\Delta)^\alpha g|^2}\Bigg\|_{L^{p_2}_{w_2}}\lesssim
		\|f\|_{L^{p_1}_{w_1}}
		\|(-\Delta)^\alpha g\|_{L^{p_2}_{w_2}},
	\end{align*}
	
	Cases 1 and 2 together prove the required estimate for $\sum f_{>j-3}\alap g_j$.

	To deal with the first sum $\sum\left[
	\alap f_{\le j-3}g_j-f_{\le j-3}\alap g_j\right]$, we first need to define some fattened Littlewood--Paley projectors. Denote $\widetilde{\varphi}(\xi)\coloneqq\sum_{-2\le j\le2}\varphi(2^j\xi)$ and the dyadic rescaling $\widetilde{\varphi}_j=\widetilde{\varphi}(\frac{\cdot}{2^j})$. Then it is apparent that $\supp\widetilde{\varphi}=\{\xi\in\Rn\colon1/8\le|\xi|\le8\}$ and $\widetilde{\varphi}\equiv1$ on the annulus $\{\xi\in\Rn\colon{1}/{4}\le|\xi|\le4\}$. Now we define the fattened projector $\widetilde{P}_j$ as, $\Widehat{\widetilde{P}_jf}=\widetilde{\varphi}_j\Widehat{f}$ for $j\in\Z$. Also denote $\Widehat{\widetilde{\psi}}(\xi)\coloneqq|\xi|^{2\alpha}\widetilde{\varphi}(\xi)\in\cc$.
	
	A little bit of Fourier analysis combined with the support properties of the projector $\widetilde{P}_j$ and the fundamental theorem of calculus now give \begin{align*}
		&\sum_{j\in\Z}\left[
		\alap f_{\le j-3}g_j-f_{\le j-3}\alap g_j\right]\\&=\sum_{j\in\Z}\left[
		\alap\widetilde{P}_j(f_{\le j-3}g_j)-f_{\le j-3}\alap\widetilde{P}_jg_j\right]\\&=\sum_{j\in\Z}2^{j(N+2\alpha)}\int_{\Rn}\widetilde{\psi}(2^jy)g_j(x-y)\left[f_{\le j-3}(x-y)-f_{\le j-3}(x)\right]\dy\\&=\sum_{j\in\Z}2^{j(N+2\alpha)}\int_{\Rn}\widetilde{\psi}(2^jy)g_j(x-y)\int_{0}^{1}\nabla f_{\le j-3}(x-ty)\cdot (-y){\,\rm d}t\dy.
	\end{align*}
	
	Using the local consistency of the band-limited functions $\nabla f_{\le j}$ and $g_j$ (cf. \cite[Lemma~3.1]{Tao07}), we have \begin{align*}\left|\nabla f_{\le j-3}(x-ty)\right|&\lesssim_{N}\langle2^{j-3}(-ty)\rangle^NM\nabla f_{\le j-3}(x),\\ |g_j(x-y)|&\lesssim_{N}\langle2^j(-y)\rangle^NMg_j(x).\end{align*}
	
	Here, $\langle\cdot\rangle$ is the \emph{Japanese bracket}, defined as $\langle x\rangle\coloneqq\sqrt{1+|x|^2}$. Using these pointwise estimates in the last integral and finally applying Cauchy--Schwarz, one gets \begin{align*}
		&\Bigg|\sum_{j\in\Z}\left[
		\alap f_{\le j-3}g_j-f_{\le j-3}\alap g_j\right]\Bigg|\\&\lesssim_{N}\sum_{j\in\Z}2^{j(N+2\alpha)}\int_{\Rn}|\widetilde{\psi}(2^jy)|\,\langle2^j|y|\rangle^NMg_j(x)\,\langle2^{j-3}|ty|\rangle^N|y|M\nabla f_{\le j-3}(x)\dy\\&\lesssim_{N,\alpha_1} \sum_{j\in\Z}2^{2j\alpha_2}Mg_j(x)\cdot2^{(j-3)(2\alpha_1-1)}M\nabla f_{\le j-3}(x)\\&\lesssim_{N,\alpha_1}\sqrt{\sum_{j\in\Z}|M(2^{2j\alpha_2}g_j)|^2}\sqrt{\sum_{j\in\Z}|M(2^{j(2\alpha_1-1)}\nabla f_{\le j})|^2}.
	\end{align*}
	
	Note the second estimate above relies on the fact that $|y|(1+|y|^2)^N\widetilde{\psi}(y)\in L^1(\Rn)$ as $\widetilde{\psi}\in\mathcal{S}(\Rn)$. Now applying the weighted H\"older inequality followed by the Fefferman--Stein inequality, Littlewood--Paley inequality and finally Lemma~\ref{derivative localization} (with $s=1-2\alpha_1>0,\,h=\mathcal{R}(-\Delta)^{\alpha_1}f$) yields \begin{align*}
		&\Bigg\|\sum_{j\in\Z}\left[
		\alap f_{\le j-3}g_j-f_{\le j-3}\alap g_j\right]\Bigg\|_{L^p_w}\\&\lesssim_{N,\alpha_1}\Bigg\|\sqrt{\sum_{j\in\Z}|M(2^{j(2\alpha_1-1)}\nabla f_{\le j})|^2}\Bigg\|_{L^{p_1}_{w_1}}\Bigg\|\sqrt{\sum_{j\in\Z}|M(2^{2j\alpha_2}g_j)|^2}\Bigg\|_{L^{p_2}_{w_2}}\\&\lesssim_{N,p_1,p_2,\alpha_1,[w_1]_{p_1},[w_2]_{p_2}}\Bigg\|\sqrt{\sum_{j\in\Z}|2^{j(2\alpha_1-1)}\nabla f_{\le j}|^2}\Bigg\|_{L^{p_1}_{w_1}}\Bigg\|\sqrt{\sum_{j\in\Z}|2^{2j\alpha_2}g_j|^2}\Bigg\|_{L^{p_2}_{w_2}}
        \\&\lesssim_{N,p_1,p_2,\alpha_1,\alpha_2,[w_1]_{p_1},[w_2]_{p_2}}\left\|(-\Delta)^{\alpha_1}f\right\|_{L^{p_1}_{w_1}}\left\|(-\Delta)^{\alpha_2}g\right\|_{L^{p_2}_{w_2}}.
	\end{align*}

Observe that in the last inequality above we rewrote $
\nabla = (-\Delta)^{\frac{1-2\alpha_1}{2}} \, (-\Delta)^{-\frac{1}{2}}\nabla \, (-\Delta)^{\alpha_1},
$ so that the operator $\mathcal{R}\coloneqq (-\Delta)^{-\frac{1}{2}}\nabla$ appears. This is precisely the classical vector-valued Riesz transform, which is bounded on $L^p_w(\Rn)$ for all $1<p<\infty$ and $w\in\mathcal{A}_p$, i.e., $
\|\mathcal{R}\|_{L^p_w\to L^p_w} \lesssim_{N,p,[w]_p}1$. Finally combining all the preceding estimates yields the desired conclusion.\qed

	\begin{Proposition}\label{spectral-gap-inequality-almost-eigenfunction} Let $\varepsilon>0$ and $f_\varepsilon\in\Hs$ be such that $\langle f_\varepsilon,V_i\rangle_{\hs}=\mathcal{o}(1),\,\langle f_\varepsilon,\dot{V}_i\rangle_{\hs}=\mathcal{o}(1)$. Then the following inequality holds. \[\int_{\Rn}f_\varepsilon^2\frac{V_i^{\twostar-2}}{|x|^t}\dx\le\frac{1}{\eta_3}\int_{\Rn}\left|\flap f_\varepsilon\right|^2\dx+\mathcal{o}(1),\] where $\eta_3$ denotes the third eigenvalue of the linearized operator $\frac{\Flap}{{|x|^{-t}}{V_i^{\twostar-2}}}$; $V_i=V_{s,t}^{\lambda_i}$ for some $\lambda_i>0$, and $\mathcal{o}(1)$ is any quantity that tends to zero as $\varepsilon\to0$.
	\end{Proposition}
	
	\Proof Let $p=\twostar-1$ and consider the eigenvalue problem \begin{equation}\Flap\psi=\mu\frac{V_i^{p-1}}{|x|^t}\psi,\hspace{2mm}\psi\in\Hs.\label{eigenvalue problem}\end{equation}
	
	By the Remark~\ref{remark on eigenfunctions}-$(1)$, the operator $\frac{\Flap}{|x|^{-t}V_i^{p-1}}$ has a discrete spectrum with first two eigen-pairs being $(1, V_i)$ and $(p, \dot{V}_i)$ respectively. Moreover, the third eigenvalue $\eta_3$ satisfies $\eta_3>p$, and admits the Rayleigh quotient characterization \begin{equation}\label{rayleigh-quotient-inequality}\eta_3 \le \frac{\langle \psi, \psi \rangle_{\hs}}{\int_{\Rn} \psi^2 \frac{V_{i}^{p-1}}{|x|^t}\dx}, \quad \text{for all } \psi \perp_{\hs} \operatorname{span}\left\{V_{i}, \dot{V}_{i}\right\},\end{equation} with equality if and only if $\psi$ is the third eigenfunction.
	
	Take $$\psi=f_\varepsilon-\langle f_\varepsilon,V_{i}\rangle_{\hs}\frac{V_{i}}{\|V_{i}\|_{\hs}^2}-\langle f_\varepsilon,\dot{V}_{i}\rangle_{\hs}\frac{\dot{V}_{i}}{\|\dot{V}_{i}\|^2_{\hs}}.$$ 
	
	It is now clear that $\psi\perp_{\hs}\operatorname{span}\{V_{i},\dot{V}_{i}\}$ and therefore using this $\psi$ in the Rayleigh quotient inequality~\eqref{rayleigh-quotient-inequality} gives the required estimate.\qed
	
	\begin{Proposition}\label{basic-cutoff-norm-estimate}
		Let $s\in(0,1), N>2s, t\in(0,2s)$ and $0<r<R$. Define the cutoff function $\varphi_{r,R}\colon\Rn\rightarrow[0,1]$ by \begin{equation}\label{cut-off}\varphi_{r,R}(x)=\begin{cases}
				1 & \text{if } |x|<r,\\
				\displaystyle\frac{\log \frac R{|x|}}{\log \frac Rr} & \text{if } r\le|x|\le R,\\
				0 &\text{if } |x|>R.
		\end{cases}\end{equation}
		
		Then for all sufficiently large $\frac Rr$, we have the estimate \begin{equation*}
			\bigg\||x|^{\frac t\twostar}\flap\varphi_{r,R}\bigg\|_{\frac{2\cdot\twostar}{\twostar-2}}\lesssim_{N,s}\left(\log\frac Rr\right)^{-\frac12-\frac1\twostar}.
		\end{equation*}
	\end{Proposition}
	
	\Proof To simplify notations, we denote $q=\frac{2\cdot\twostar}{\twostar-2}\text{ and }\varphi_{r,R}$ simply by $\varphi$ in what follows. We start by decomposing the integral into several regions and then estimate over each region individually.
	
	Write \begin{align*}
		&\int_{\Rn}|x|^{\frac {qt}\twostar}\big|\flap\varphi(x)\big|^{q}\dx\\&=\bigg(\int_{|x|<r}+\int_{r\le|x|\le R}+\int_{|x|>R}\bigg)|x|^{\frac {qt}\twostar}\big|\flap\varphi(x)\big|^{q}\dx\eqqcolon I_1+I_2+I_3.
 	\end{align*}
	
	From the definition of the fractional Laplacian, one has \begin{align*}\flap\varphi(x)&=C_{N,s}\,\pv\bigg(\int_{|y|<r}+\int_{r\le|y|\le R}+\int_{|y|>R}\bigg)\frac{\varphi(x)-\varphi(y)}{|x-y|^{N+s}}\dy\\&\eqqcolon J_1(x)+J_2(x)+J_3(x).\end{align*}
	
	Using this we have by the triangle inequality, \begin{align*}
    I_1^{\frac{1}{q}}=\bigg(\int_{|x|<r}|x|^{\frac {qt}\twostar}|J_1(x)+J_2(x)+J_3(x)|^{q}\dx\bigg)^{\frac1q}\le\sum_{k=1}^{3}\Big\||x|^{\frac{t}{\twostar}}J_k\Big\|_{L^{q}\left(B_r\right)}\eqqcolon \sum_{1\le k\le3}I_{1k}.
	\end{align*} and similarly decompose $I_2$ and $I_3$. It is now immediate from the definition of $\varphi$ that $I_{11}=0$. 

    \noindent\emph{Estimate for $I_{13}$.}\begin{align*}
		I_{13}^{q}&=\int_{|x|<r}|x|^{\frac {qt}\twostar}|J_3(x)|^{q}\dx\\&=C_{N,s}\int_{|x|<r}\bigg(|x|^{\frac t\twostar}\int_{|y|>R}\frac{1}{|x-y|^{N+s}}\dy\bigg)^{q}\dx\\
		&\lesssim_{N,s}\int_{|x|<r}\bigg(|x|^{\frac t\twostar}\int_{|y-x|>R-r}\frac{1}{|y-x|^{N+s}}\dy\bigg)^{q}\dx\sim_{N,s} \bigg(\frac{r}{R-r}\bigg)^{qs}=\bigg(\frac R r-1\bigg)^{-qs}.
	\end{align*}

    \noindent\emph{Estimate for $I_{12}$.}\begin{align*}
		I_{12}^{q}&=\int_{|x|<r}\left(|x|^{\frac t\twostar}|J_2(x)|\right)^{q}\dx\\
		&=C_{N,s}\int_{|x|<r}\Bigg(|x|^{\frac t\twostar}\int_{r<|y|<R}\Bigg(1-\frac{\log\frac R{|y|}}{\log\frac Rr}\Bigg)\frac{{\rm d}y}{|x-y|^{N+s}}\Bigg)^{q}\dx\\
		&=C_{N,s}\int_{|x|<r}\Bigg(|x|^{\frac t\twostar}\frac{1}{\log\frac Rr}\int_{r<|y|<R}\frac{\log |y|-\log r}{|x-y|^{N+s}}\dy\Bigg)^{q}\dx\\
		&=C_{N,s}\int_{|x|<r}\Bigg(|x|^{\frac t\twostar}\frac{1}{\log\frac Rr}\Bigg(\int_{r<|y|<2r}+\int_{2r<|y|<R}\Bigg)\frac{\log |y|-\log r}{|x-y|^{N+s}}\dy\Bigg)^{q}\dx.	
	\end{align*}
	
	Now we estimate the above two integrals inside $I_{12}$ as follows. Notice for $|x|<r$, \begin{align*}
		\int_{r<|y|<2r}\frac{\log |y|-\log r}{|y-x|^{N+s}}\dy&\le\int_{r<|y|<2r}\frac{|y-x|}{r}\frac{{\rm d}y}{|y-x|^{N+s}}\le\frac1r\int_{|y-x|<3r}\frac{{\rm d}y}{|y-x|^{N+s-1}}\sim_{N}\frac{1}{r^s}.
	\end{align*}

	\begin{align*}
		\int_{2r<|y|<R}\frac{\log |y|-\log r}{|y-x|^{N+s}}\dy\le\frac{1}{r^s}\int_{|z|\ge2}\frac{\log|z|}{(|z|-1)^{N+s}}\dz\le\frac{2^{N+s}\omega_{N-1}}{r^s}\int_2^\infty\frac{\log\alpha}{\alpha^{1+s}}\,{\rm d}\alpha\sim_{N,s}\frac{1}{r^s}.\end{align*}
        
    Here $\omega_{N-1}$ denotes the surface measure of the sphere $\mathbb S^{N-1}$ in $\Rn$. An integration by parts shows $\int_2^\infty\tfrac{\log\alpha}{\alpha^{1+s}}\,{\rm d}\alpha=2^{-s}\left(\tfrac{\log2} {s}+s^{-2}\right)$.

	Therefore \begin{align*}
		I_{12}\lesssim_{N,s}\left(\int_{|x|<r}\left(|x|^{\frac t\twostar}\frac{1}{\log\frac Rr}\left(\frac{1}{r^s}+\frac{1}{r^s}\right)\right)^q\dx\right)^{\frac1q}
		\sim_{N,s}\left(\log\frac Rr\right)^{-1}.
	\end{align*}
	
	Finally, combining these estimates on $I_{11},\,I_{12}$ and $I_{13}$, we obtain $$I_1\lesssim_{N,s}\left[{\left(\frac R r-1\right)^{-s}}+\left(\log\frac Rr\right)^{-1}\right]^q.$$
	
	Thus for all sufficiently large $\frac{R}{r}$, we have $
	I_1 \lesssim_{N,s} (\log\frac{R}{r})^{-q}$. We now turn to the estimates for $I_2$ and $I_3$ and a similar argument yields for all sufficiently large $\frac Rr$, $$
	I_2\lesssim_{N,s} \left( \log \frac{R}{r} \right)^{1-q},\hspace{2mm}I_3 \lesssim_{N,s} \left( \log \frac{R}{r} \right)^{-q}.
	$$
	
	The desired result is now immediate from $(I_1+I_2+I_3)^{1/q}\lesssim_{N,s}(\log\frac{R}{r})^{-1+1/q}$. For brevity, we briefly discuss the estimate for the dominant component, $I_{22}$ of $I_2$. Write $A=\frac{R}{r},\,q=\frac{2\cdot\twostar}{\twostar-2}$ and then $\frac{qt}{\twostar}=qs-N$.

	Notice \begin{align*}
		I_{22}&\le C_{N,s}\frac{1}{\log A}\bigg(\int_{1\le|x|\le A}|x|^{qs-N}\bigg(\int_{1\le|y|\le A}\frac{|\log|x|-\log|y||}{|x-y|^{N+s}}\dy\bigg)^q\dx\bigg)^{\frac 1q}.
	\end{align*}
	
	Decompose for $1\le|x|\le A$, \begin{align*}&\int_{1\le|y|\le A}\frac{|\log|x|-\log|y||}{|x-y|^{N+s}}\dy\\&\le\bigg(\int_{\frac12\le|y|\le\frac{|x|}2}+\int_{\frac{|x|}2<|y|<2|x|}+\int_{2|x|\le|y|\le2A}\bigg)\frac{|\log|x|-\log|y||}{|x-y|^{N+s}}\dy\\&\eqqcolon K_1(x)+K_2(x)+K_3(x).\end{align*}

   \noindent\emph{Estimates for $K_1(x),\,K_3(x)$.} Plainly by the triangle inequality, a change of variable, and an integration by parts, one sees\begin{align*}
       K_1(x)&=\int_{\frac12\le|y|\le\frac{|x|}2}\tfrac{\log|x|-\log|y|}{|x-y|^{N+s}}\dy\le\frac{2^{N+s}}{|x|^s}\int_{|z|\le\frac12}\log\frac{1}{|z|}\dz=\tfrac{2^s\omega_{N-1}\left(\frac{\log2}{N}+N^{-2}\right)}{|x|^s}\lesssim_{N,s}\frac1{|x|^s},\\ K_3(x)&=\int_{2|x|\le|y|\le2A}\tfrac{\log|y|-\log|x|}{|x-y|^{N+s}}\dy\le\frac{1}{|x|^s}\int_{|z|\ge2}\tfrac{\log|z|}{(|z|-1)^{N+s}}\dz=\tfrac{2^N\omega_{N-1}\left(\frac{\log2}{s}+s^{-2}\right)}{|x|^s}\lesssim_{N,s}\frac{1}{|x|^s}.
   \end{align*}

   \noindent\emph{Estimate for $K_2(x)$.} For this we need an elementary inequality. It is easy to check that $\forall a>0,\,\text{and } b\ge\tfrac{a}{2},\,|\log b-\log a|\le\tfrac{\log 4}{a}|b-a|$. Invoking this with $a=|x|,\,b=|y|$, one sees\begin{align*}
       K_2(x)&=\int_{\frac{|x|}2<|y|<2|x|}\frac{|\log|x|-\log|y||}{|x-y|^{N+s}}\dy\\&\le\frac{\log4}{|x|}\int_{\frac{|x|}2<|y|<2|x|}\frac{||x|-|y||}{|x-y|^{N+s}}\dy\le\frac{\log4}{|x|}\int_{|x-y|<3|x|}\frac{\dy}{|x-y|^{N+s-1}}\lesssim_N|x|^{-s}.
   \end{align*}

Now plugging these estimates in $I_{22}$ gives\[I_{22}\lesssim_{N,s}\frac{1}{\log A}\left(\int_{1\le|x|\le A}|x|^{qs-N}\cdot|x|^{-qs}\dx\right)^{\frac1q}\lesssim_{N,s}(\log A)^{\frac1q-1}.\]

This gives the stated bound of $I_2$.\qed

\section{Integrals involving two bubbles}\label{AB}
The estimates here are direct adaptations of \cite[Appendix~B]{FG20}: after rescaling to a single-bubble scale, we split the integrals into regions where the profiles show power-type behavior and compute each piece explicitly. The only difference is that, in the absence of an explicit formula for $V_{s,t}$, we use the sharp pointwise bound \big($\sim\tfrac{1}{1+|x|^{N-2s}}$\big) from \cite[Theorem~1.1]{MN21}, which allows all steps to be repeated with only cosmetic changes and with the same scaling laws and logarithmic corrections. We therefore simply state the resulting bounds used in the paper.

	\begin{Lemma}
		Let $s\in(0,1),\,t\in(0,2s),\,N > 2s,\,\lambda \in (0,1]$, and $\alpha, \beta \ge 0$ be such that $\alpha + \beta = \twostar$. Then for every $\varepsilon > 0$ the following estimate holds:
		
		$$
		\int_{\mathbb{R}^N} \frac{V_{s,t}(x)^\alpha V_{s,t}^\lambda(x)^\beta}{|x|^t} \dx \sim 
		\begin{cases}
			C(N,\varepsilon)\, \displaystyle\lambda^{\frac{(N - 2s)\min(\alpha, \beta)}{2}} & \text{if } |\alpha - \beta| \ge \varepsilon, \\
			C(N)\, \lambda^{\frac{N - t}{2}}\left(1 + \log\frac{1}{\lambda}\right) & \text{if } \alpha = \beta.
		\end{cases}
		$$
	\end{Lemma}

	\begin{Corollary}\label{general-bubble-interaction-estimate} For $i\neq j;\,\lambda_i,\lambda_j>0$, define $Q_{ij}\coloneqq\min\left\{\frac{\lambda_i}{\lambda_j},\frac{\lambda_j}{\lambda_i}\right\}$. Then $$
		\int_{\mathbb{R}^N} \frac{V^{\lambda_i}_{s,t}(x)^\alpha V_{s,t}^{\lambda_j}(x)^\beta}{|x|^t} \dx \sim_{N,\varepsilon} 
		\begin{cases}
			Q_{ij}^{\frac{(N - 2s)\min(\alpha, \beta)}{2}} & \text{if } |\alpha - \beta| \ge \varepsilon, \\
			Q_{ij}^{\frac{N - t}{2}}\left(1 + \log\frac{1}{Q_{ij}}\right) & \text{if } \alpha = \beta.
		\end{cases}
		$$
	\end{Corollary}
 
The next corollary captures this coupling quantitatively, reflecting a key feature of two-bubble interactions at different scales. More precisely, the full interaction integral is asymptotically equivalent to its restriction to the natural core region associated with either bubble.

 \begin{Corollary}\label{loc-bub-bub-est}
Let $V_i = V_{s,t}^{\lambda_i}$ and $V_j = V_{s,t}^{\lambda_j}$, and assume
$p+1=\twostar$. Then the mixed interaction between $V_i$ and $V_j$ is asymptotically
localized, in the sense that
\[
\int_{|x|\le\frac1{\lambda_j}}\frac{V_iV_j^{p}}{|x|^t}\dx
\sim_{N,s,t}
\int_{\Rn}\frac{V_iV_j^{p}}{|x|^t}\dx
=
\int_{\Rn}\frac{V_i^{p}V_j}{|x|^t}\dx
\sim_{N,s,t}
\int_{|x|\le\frac1{\lambda_i}}\frac{V_i^{p}V_j}{|x|^t}\dx .
\]
\end{Corollary}

\Proof
By the sharp growth estimates for $V_{s,t}$, one readily verifies that
\[
\int_{|x|\le\frac1{\lambda_i}}\frac{V_i^{p}V_j}{|x|^t}\dx
\sim_{N,s,t}
Q_{ij}^{\frac{N-2s}{2}}.
\]

The conclusion then follows directly from
Corollary~\ref{general-bubble-interaction-estimate} together with
Remark~\ref{bub-del}.
\qed

\section{Nearest-point projection onto the multi-bubble set}\label{AC}

In the classical local critical Sobolev setting, Bahri--Coron established a nearest-point construction for well-separated bubbles. Their results provide parameter matching and compactness within a prescribed admissible region, yielding a unique interior minimizer, tangent-space orthogonality and continuous dependence of the parameters; see~\cite[Proposition~7, Lemmas~A.1--A.2]{Bahri-Coron-88}. In fractional settings, multi-bubble nearest-point projections have been invoked in~\cite{A23,CKW25}, where the minimization over finite sums of fractional Talenti bubbles and the associated tangent-space orthogonality are used. The fractional Bahri--Coron parameter-selection result is also stated in~\cite{ACH18}, with its proof referred to the classical construction, and employed in the fractional Brezis--Nirenberg analysis of~\cite{DNK23CalVar}. Rigorous attainment results for nonlocal nearest-point projections are available, however, only for single-bubble configurations; see~\cite{DNK23JFA,CS25}. To the best of our knowledge, a complete proof of global attainment and uniqueness for the nearest-point projection onto the full multi-bubble set $\mathcal M_\nu$ of the fractional Hardy--Sobolev family has remained unavailable. The present appendix fills this gap by establishing global parameter matching, attainment and uniqueness for this projection.

\medskip
Throughout this appendix, $N>2s$, $s\in(0,1)$, $t\in(0,2s)$, and $N,\nu\in\N$ are fixed. All constants may depend on $N,\nu,s,t$, but are independent of the reference parameters. We use the notation $V\coloneqq V_{s,t}$ and the Hilbert space $\Hs$ with inner product $\langle\cdot,\cdot\rangle_{\hs}$.

\subsection{The dilation group and its logarithmic representation}

For $a\in\R$, define
\begin{equation*}
	\mathcal U(a)f(x)
	\coloneqq e^{\frac{N-2s}{2}a}f(e^a x),
	\qquad f\in\Hs.
\end{equation*}

By the scaling relation~\eqref{dilation map},
\begin{equation}\label{AC-dilation-identities}
	V_{s,t}^{\lambda}=\mathcal U(\log\lambda)V,
	\qquad
	\mathcal U(a)\mathcal U(b)=\mathcal U(a+b),
	\qquad
	\|\mathcal U(a)f\|_{\Hs}=\|f\|_{\Hs}.
\end{equation}

Thus $(\mathcal U(a))_{a\in\R}$ is a unitary group on $\Hs$. The Fourier scaling formula is
\begin{equation}\label{AC-fourier-dilation}
	\widehat{\mathcal U(a)f}(\xi)
	=e^{-\frac{N+2s}{2}a}\widehat f(e^{-a}\xi),
	\qquad f\in\Hs.
\end{equation}

\begin{Lemma}\label{lem:AC-dilation-weak}
	The group $a\mapsto\mathcal U(a)$ is strongly continuous on $\Hs$. Moreover,
	for every $f,g\in\Hs$,
	\begin{equation}\label{AC-weak-dilations}
		\lim_{|a|\to\infty}
		\langle g,\mathcal U(a)f\rangle_{\Hs}=0.
	\end{equation}
\end{Lemma}
\Proof
For $\lambda>0$, let $D_\lambda f(x)=\lambda^{\frac{N-2s}{2}}f(\lambda x)$ be the dilation used in~\cite[\S~2.3]{CS25}. The change of variable $\lambda=e^a$ gives $\mathcal U(a)=D_{e^a}$. Strong continuity of this action is proven in~\cite[\S~2.3]{CS25}. If $|a_k|\to\infty$, set $\lambda_k=e^{a_k}$ and $\delta_k=1$ in~\cite[Lemma~2.1]{CS25}. Its separation condition becomes $|\log(\lambda_k/\delta_k)|=|a_k|\to\infty$, so that $\langle g,\mathcal U(a_k)f\rangle_{\hs}\to0$ for every $f,g\in\Hs$. Since this holds for every sequence with $|a_k|\to\infty$, it proves~\eqref{AC-weak-dilations}.\qed

In~\cite[Theorem~1.1(iii)]{MN21}, the dilation orbit is parameterized by $z_\tau(x)=\tau^{\frac{2s-N}{2}}V(x/\tau)$, $\tau>0$, where $V=V_{s,t}$. Thus $z_\tau=\mathcal U(-\log\tau)V$ and $V_{s,t}^{\lambda}=\mathcal U(\log\lambda)V=z_{\lambda^{-1}}$. Theorem~1.1(iii) states that $\tau\mapsto z_\tau$ is a regular $C^2$ curve in $\Hs$. Since $\lambda\mapsto\lambda^{-1}$ and $a\mapsto e^a$ are smooth diffeomorphisms, $a\mapsto\mathcal U(a)V$ is a regular $C^2$ curve in $\Hs$. Put
\begin{equation*}
	Z\coloneqq\frac{d}{da}\bigg|_{a=0}\mathcal U(a)V,
	\qquad
	\widetilde Z\coloneqq\frac{d^2}{da^2}\bigg|_{a=0}\mathcal U(a)V.
\end{equation*}

Then $Z,\widetilde Z\in\Hs$, $Z\neq0$, and the group property gives
\begin{equation}\label{AC-orbit-derivatives}
	\frac{d}{da}\mathcal U(a)V=\mathcal U(a)Z,
	\qquad
	\frac{d^2}{da^2}\mathcal U(a)V=\mathcal U(a)\widetilde Z.
\end{equation}

Since the norm of $\mathcal U(a)V$ is constant, differentiating its square
at $a=0$ yields
\begin{equation}\label{AC-VZ-orthogonal}
	\langle V,Z\rangle_{\Hs}=0.
\end{equation}

In the notation of Remark~\ref{remark on eigenfunctions}, the dot denotes differentiation with respect to the scale parameter $\lambda$. Since $\mathcal U(a)V=V_{s,t}^{e^a}$, the chain rule gives
\[
Z=\frac{d}{da}\bigg|_{a=0}V_{s,t}^{e^a}
=\partial_\lambda\bigg|_{\lambda=1} V_{s,t}^{\lambda}
=\dot V_{s,t}^{1}.
\]

More generally, using the group property in~\eqref{AC-dilation-identities},
\begin{align}
	\mathcal U(\log\lambda)Z
	=\left.\frac{d}{da}\right|_{a=0}
	\mathcal U(\log\lambda)\mathcal U(a)V
	=\left.\frac{d}{da}\right|_{a=0}V_{s,t}^{\lambda e^a}=\lambda\,\partial_\lambda V_{s,t}^{\lambda}
	=\lambda\dot V_{s,t}^{\lambda}.
	\label{AC-log-derivative}
\end{align}

By~\cite{MRNA21}, $V$ is radial, nonnegative and strictly decreasing; in addition, the normalized equation~\eqref{maineq} gives $V>0$ in $\Rn$. Let us also denote the closed radial subspace of $\Hs$ by $\hs_{\mathrm{rad}}(\Rn)$.

\begin{Lemma}\label{lem:AC-radial-log-isometry}
	There exists a unitary map
	$\mathcal J:\hs_{\mathrm{rad}}(\Rn)\to L^2(\R)$ such that,
	\begin{equation}\label{AC-radial-translation}
		\mathcal J(\mathcal U(a)f)(r)=\mathcal Jf(r-a)
		\qquad\text{for every }a\in\R.
	\end{equation}
\end{Lemma}

\Proof
For $f\in\hs_{\mathrm{rad}}(\Rn)$, write $\widehat f(\xi)=F_f(|\xi|)$. Define $(\mathcal Jf)(r)\coloneqq |\mathbb S^{N-1}|^{1/2}e^{\frac{N+2s}{2}r}F_f(e^r)$. Then $\mathcal J$ is an isometry. Indeed, by polar coordinates and the Fourier representation of the homogeneous Sobolev norm,
\begin{align*}
	\|f\|_{\hs}^2
	=|\mathbb S^{N-1}|
	\int_0^\infty \rho^{N+2s-1}|F_f(\rho)|^2\,\mathrm{d}\rho=|\mathbb S^{N-1}|
	\int_{\R}e^{(N+2s)r}|F_f(e^r)|^2\,\mathrm{d}r
	=\|\mathcal Jf\|_{L^2(\R)}^2.
\end{align*}

We verify that $\mathcal J$ is onto. Let $g\in L^2(\R)$ and set $F(\rho)\coloneqq |\mathbb S^{N-1}|^{-1/2}\rho^{-\frac{N+2s}{2}}g(\log\rho)$ for $\rho>0$. The change of variables $r=\log\rho$ gives
\begin{equation*}
|\mathbb S^{N-1}|\int_0^\infty
\rho^{N+2s-1}|F(\rho)|^2\,\mathrm d\rho
=\int_\R|g(r)|^2\,\mathrm dr.
\end{equation*}

Thus $F$ belongs to the weighted space $L^2\big((0,\infty),\rho^{N+2s-1}\,\mathrm d\rho\big)$. Since the weight is positive and smooth on every compact subinterval of $(0,\infty)$, truncation to annuli followed by mollification shows that there are $F_n\in C_c^\infty((0,\infty))$ such that
\begin{equation}\label{AC-profile-approximation}
\int_0^\infty \rho^{N+2s-1}|F_n(\rho)-F(\rho)|^2\,\mathrm d\rho
\longrightarrow0.
\end{equation}

Define $f_n$ by $\widehat f_n(\xi)\coloneqq F_n(|\xi|)$. Because $F_n$ is smooth and compactly supported in $(0,\infty)$, the function $\xi\mapsto F_n(|\xi|)$ belongs to $C_c^\infty(\Rn\setminus\{0\})$. Its inverse Fourier transform $f_n$ is therefore a radial Schwartz function. In particular, $f_n$ belongs to $\hs_{\mathrm{rad}}(\Rn)$. By the Fourier representation of the $\Hs$ norm and~\eqref{AC-profile-approximation},
\[
\|f_n-f_m\|_{\hs}^2
=|\mathbb S^{N-1}|\int_0^\infty
\rho^{N+2s-1}|F_n(\rho)-F_m(\rho)|^2\,\mathrm d\rho.
\]

Thus $(f_n)$ is Cauchy in $\hs_{\mathrm{rad}}(\Rn)$. Let $f$ denote its limit. By the definition of $\mathcal J$ and~\eqref{AC-profile-approximation},
\[
\|\mathcal Jf_n-g\|_{L^2(\R)}^2
=|\mathbb S^{N-1}|\int_0^\infty
\rho^{N+2s-1}|F_n(\rho)-F(\rho)|^2\,\mathrm d\rho\longrightarrow0.
\]

Since $\mathcal J$ is an isometry, it is continuous. Now continuity yields $\mathcal Jf=g$, and therefore it is unitary.

Using~\eqref{AC-fourier-dilation}, we obtain
\begin{align*}
	\mathcal J(\mathcal U(a)f)(r)
	=|\mathbb S^{N-1}|^{1/2}
	e^{\frac{N+2s}{2}r}
	e^{-\frac{N+2s}{2}a}F_f(e^{r-a})=(\mathcal Jf)(r-a).
\end{align*}

This proves~\eqref{AC-radial-translation}.\qed

\subsection{A compactness lemma for finite exponential sums}
Put $q\coloneqq\mathcal JV$ and let $\widehat q$ denote its one-dimensional unitary Fourier transform. Since $V\neq0$ and $\mathcal J$ is unitary, $q\in L^2(\R)$ is nonzero. Consequently, $w\coloneqq|\widehat q|^2$ belongs to $L^1(\R)$, is nonnegative, and is not identically zero. We can therefore choose a bounded measurable set $E\subset\R$ and a constant $c_E>0$ such that
\begin{equation}\label{AC-E-choice}
	|E|>0,
	\qquad w(\xi)\ge c_E
	\quad\text{for a.e. }\xi\in E.
\end{equation}

Indeed, the sets $E_{m,R}\coloneqq\left\{\xi\in[-R,R]\mid w(\xi)\ge\frac1m\right\}$, for $m,R\in\N$, cover $\{w>0\}$, up to a null set, so at least one has positive measure.

By Lemma~\ref{lem:AC-radial-log-isometry}, the one-dimensional Fourier
transform, and Plancherel's theorem, for arbitrary real coefficients
$\beta_j,\gamma_i$ and logarithmic scales $b_j,a_i$ we have
\begin{align}
	\Bigg\|
	\sum_{j=1}^{m}\beta_j\mathcal U(b_j)V
	-\sum_{i=1}^{\nu}\gamma_i\mathcal U(a_i)V
	\Bigg\|_{\Hs}^2=
	\int_{\R}w(\xi)
	\Bigg|
	\sum_{j=1}^{m}\beta_je^{-ib_j\xi}
	-\sum_{i=1}^{\nu}\gamma_ie^{-ia_i\xi}
	\Bigg|^2\,\mathrm{d}\xi.
	\label{AC-exponential-identity}
\end{align}

In particular, by~\eqref{AC-E-choice},
\begin{align}
	\Bigg\|
	\sum_{j=1}^{m}\beta_j\mathcal U(b_j)V
	-\sum_{i=1}^{\nu}\gamma_i\mathcal U(a_i)V
	\Bigg\|_{\Hs}\ge c_E^{1/2}
	\Bigg\|
	\sum_{j=1}^{m}\beta_je^{-ib_j\xi}
	-\sum_{i=1}^{\nu}\gamma_ie^{-ia_i\xi}
	\Bigg\|_{L^2(E)}.
	\label{AC-exponential-lower-bound}
\end{align}

\begin{Lemma}\label{lem:AC-exponential-compactness}
	Let $E\subset\R$ be bounded and measurable, with $|E|>0$. Fix
	$m\in\N$ and $R>0$. The family
	\begin{equation}\label{AC-normalized-exponential-family}
		\Bigg\{\sum_{j=1}^{m}c_je^{-id_j\xi}\,\bigg|\,d_j\in[-R,R],\ c_j\in\C,\Big\|\sum_{j=1}^{m}c_je^{-id_j\xi}\Big\|_{L^2(E)}=1\Bigg\}
	\end{equation}
	is relatively compact in $L^2(E)$.
\end{Lemma}

\Proof After reordering, assume $-R\le d_1\le d_2\le\cdots\le d_m\le R$. For $1\le j\le m$, let
\begin{equation*}
	\Psi_j(d_1,\ldots,d_j;\xi)
	\coloneqq [d_1,\ldots,d_j]
	\big(d\mapsto e^{-id\xi}\big),
\end{equation*}
where the divided difference is understood in its confluent sense when some nodes coincide. The Hermite--Genocchi formula gives
\begin{equation}\label{AC-HG-formula}
	\Psi_j(d_1,\ldots,d_j;\xi)
	=\int_{\Delta_{j-1}}
	(-i\xi)^{j-1}
	e^{-i\xi\sum_{\ell=1}^{j}\theta_\ell d_\ell}
	\,\mathrm{d}\theta,
\end{equation}
where $\Delta_{j-1}$ is the standard $(j-1)$-simplex, equipped with its usual measure. Since $E$ is bounded, the integrand in~\eqref{AC-HG-formula} is uniformly bounded on $E$ for all nodes in $[-R,R]$. Dominated convergence therefore shows that $(d_1,\ldots,d_j)\mapsto\Psi_j(d_1,\ldots,d_j;\cdot)$ is continuous from the compact node set into $L^2(E)$.

We next show that, for every ordered node tuple, the functions
\begin{equation}\label{AC-Newton-basis}
	\Psi_1(d_1;\cdot),\ldots,
	\Psi_m(d_1,\ldots,d_m;\cdot)
\end{equation}
are linearly independent in $L^2(E)$. Let $\delta_1<\cdots<\delta_k$ be the distinct nodes, with multiplicities $n_1,\ldots,n_k$. The confluent Newton divided-difference functions in~\eqref{AC-Newton-basis} are related by an invertible triangular change of basis to
\[
\left\{
\xi^\ell e^{-i\delta_j\xi}:
1\le j\le k,\ 0\le\ell<n_j
\right\}.
\]

To verify the linear independence of the latter family, suppose
\[
\sum_{j=1}^{k}P_j(\xi)e^{-i\delta_j\xi}=0
\quad\text{for a.e. }\xi\in E,
\]
where $P_j$ is a polynomial of degree less than $n_j$. The left-hand side is an entire function of $\xi$. Since it vanishes on a set of positive measure, its zeros have an accumulation point, and the identity theorem implies that it vanishes identically. Fixing $j$, apply
\[
\prod_{\ell\neq j}
\left(\frac{d}{d\xi}+i\delta_\ell\right)^{n_\ell}
\]
to this identity. Every term except the $j$-th is annihilated. On the $j$-th term, after factoring out $e^{-i\delta_j\xi}$, the resulting operator on $P_j$ is a product of operators of the form $\frac{d}{d\xi}+i(\delta_\ell-\delta_j)$. Each is invertible on the finite-dimensional polynomial space of degree less than $n_j$, since $\delta_\ell-\delta_j\neq0$. Hence $P_j=0$. This holds for every $j$, proving the claimed linear independence.

For an ordered node tuple, let
\begin{equation*}
	G(d_1,\ldots,d_m)
	\coloneqq
	\left(\int_E
	\Psi_j(d_1,\ldots,d_j;\xi)
	\overline{\Psi_\ell(d_1,\ldots,d_\ell;\xi)}\,\mathrm{d}\xi
	\right)_{j,\ell=1}^{m}.
\end{equation*}

The entries of $G$ depend continuously on the nodes. The linear independence just proved shows that $G$ is positive definite at every ordered node tuple. Since the set of ordered node tuples in $[-R,R]^m$ is compact, there exists $\kappa=\kappa(E,m,R)>0$ such that
\begin{equation}\label{AC-Gram-uniform}
	z^*G(d_1,\ldots,d_m)z\ge\kappa|z|^2
\end{equation}
for every $z\in\C^m$ and every ordered node tuple.

For every ordered node tuple, the Newton--Hermite interpolation identity (with confluent divided differences at repeated nodes) gives, for each
$1\le j\le m$,
\[
e^{-id_j\xi}
=\sum_{\ell=1}^{j}
\Psi_\ell(d_1,\ldots,d_\ell;\xi)
\prod_{r=1}^{\ell-1}(d_j-d_r),
\]
where the empty product for $\ell=1$ is $1$. Hence each of the exponentials $e^{-id_j\xi}$ belongs to the span of the functions in~\eqref{AC-Newton-basis}, including when nodes coincide. Thus every exponential sum in~\eqref{AC-normalized-exponential-family} can be written as
\[
g(\xi)=\sum_{j=1}^{m}z_j
\Psi_j(d_1,\ldots,d_j;\xi).
\]

If $\|g\|_{L^2(E)}=1$, then~\eqref{AC-Gram-uniform} gives $|z|\le\kappa^{-1/2}$. Consequently, the normalized family is contained in the image of the compact set
\[
\left\{(d_1,\ldots,d_m):-R\le d_1\le\cdots\le d_m\le R\right\}
\times\overline B_{\C^m}(0,\kappa^{-1/2})
\]
under the continuous map
\[
(d,z)\longmapsto
\sum_{j=1}^{m}z_j\Psi_j(d_1,\ldots,d_j;\cdot)
\quad\text{into }L^2(E).
\]

The image is compact. This proves the lemma.\qed

We shall use the following uniform form of the Riemann--Lebesgue lemma. If $K$ is a compact subset of $L^1(E)$, then
\begin{equation}\label{AC-uniform-RL}
	\lim_{|a|\to\infty}
	\sup_{g\in K}
	\left|\int_E e^{-ia\xi}g(\xi)\,\mathrm{d}\xi\right|=0.
\end{equation}

To see this, extend every $g\in K$ by zero outside $E$. Given $\eta>0$, choose a finite $\eta$-net $g_1,\ldots,g_M$ for $K$ in $L^1(E)$. The Riemann--Lebesgue lemma applies to each $g_j$, while
\[
\left|\int_E e^{-ia\xi}(g-g_j)(\xi)\,\mathrm{d}\xi\right|
\le\|g-g_j\|_{L^1(E)}<\eta.
\]

Taking the maximum over the finite net and then letting $\eta\downarrow0$ proves~\eqref{AC-uniform-RL}.

\subsection{Global matching of nearby multi-bubble configurations}
\begin{Lemma}\label{lem:AC-global-matching}
	Let $\epsilon_k\downarrow0$ and suppose that
	\[
	\widetilde\sigma_k
	=\sum_{i=1}^{\nu}
	\widetilde\alpha_{i,k}\mathcal U(a_{i,k})V
	\]
	is $2\epsilon_k$-interacting, in the sense of
	Definition~\ref{del-int-bub}. Let
	\[
	\sigma_k
	=\sum_{j=1}^{\nu}\beta_{j,k}\mathcal U(b_{j,k})V
	\in\mathcal M_\nu,
	\qquad \beta_{j,k}\in\R,
	\quad b_{j,k}\in\R,
	\]
	be arbitrary. If
	\begin{equation}\label{AC-global-matching-hypothesis}
		\|\sigma_k-\widetilde\sigma_k\|_{\Hs}\longrightarrow0,
	\end{equation}
	then, after passing to a subsequence, there exist permutations $\pi_k$
	of $\{1,\ldots,\nu\}$ such that
	\begin{equation}\label{AC-global-matching-conclusion}
		\beta_{\pi_k(i),k}-\widetilde\alpha_{i,k}\longrightarrow0,
		\qquad
		b_{\pi_k(i),k}-a_{i,k}\longrightarrow0
		\quad(1\le i\le\nu).
	\end{equation}
\end{Lemma}

\Proof
By Definition~\ref{del-int-bub} and the choice of the logarithmic scales,
\begin{equation}\label{AC-reference-separation}
	|a_{i,k}-a_{j,k}|\ge L_k,
	\qquad
	L_k\coloneqq\log\frac1{2\epsilon_k}\longrightarrow\infty
	\quad(i\neq j),
\end{equation}
and
\begin{equation}\label{AC-reference-amplitudes}
	|\widetilde\alpha_{i,k}-1|\le2\epsilon_k,
	\qquad
	\widetilde\alpha_{i,k}\longrightarrow1.
\end{equation}

Set
\[
P_k(\xi)\coloneqq\sum_{j=1}^{\nu}\beta_{j,k}e^{-ib_{j,k}\xi},
\qquad
Q_k(\xi)\coloneqq
\sum_{i=1}^{\nu}\widetilde\alpha_{i,k}e^{-ia_{i,k}\xi}.
\]

By~\eqref{AC-exponential-lower-bound} and~\eqref{AC-global-matching-hypothesis},
\begin{equation}\label{AC-PQ-close}
	\|P_k-Q_k\|_{L^2(E)}\longrightarrow0.
\end{equation}

Since the coefficients of $Q_k$ are uniformly bounded,
$\sup_k\|Q_k\|_{L^2(E)}<\infty$. Hence
\begin{equation}\label{AC-P-bounded}
	\sup_k\|P_k\|_{L^2(E)}<\infty.
\end{equation}

\medskip
\noindent\textbf{Step 1.}
Since there are only finitely many pairs of competing indices, after passing to a subsequence we may assume that, for every pair $j,\ell$, the sequence $|b_{j,k}-b_{\ell,k}|$ is either bounded or tends to infinity. Define
\[
j\sim\ell
\quad\Longleftrightarrow\quad
(b_{j,k}-b_{\ell,k})_k\text{ is bounded}.
\]

This is an equivalence relation. Let its equivalence classes be $C_1,\ldots,C_m$, where $1\le m\le\nu$. For each class $C_r$, choose an index $j_r\in C_r$ and put
\[
c_{r,k}\coloneqq b_{j_r,k},
\qquad
d_{j,k}\coloneqq b_{j,k}-c_{r,k}\quad(j\in C_r).
\]

For every $r$, there is $R_r<\infty$ such that
\begin{equation}\label{AC-cluster-bounded}
	|d_{j,k}|\le R_r
	\qquad(j\in C_r)
\end{equation}
for all sufficiently large $k$. Distinct cluster centres satisfy
\begin{equation}\label{AC-centres-separated}
	|c_{r,k}-c_{\ell,k}|\longrightarrow\infty
	\qquad(r\neq\ell).
\end{equation}

By a further subsequence, we may also assume that, for every pair $(r,i)$, the sequence $|c_{r,k}-a_{i,k}|$ is either bounded or tends to infinity, retaining~\eqref{AC-cluster-bounded} and~\eqref{AC-centres-separated}.

Define the contribution of the $r$-th cluster by
\begin{equation*}
	S_{r,k}(\xi)
	\coloneqq e^{-ic_{r,k}\xi}G_{r,k}(\xi),
	\qquad
	G_{r,k}(\xi)
	\coloneqq\sum_{j\in C_r}\beta_{j,k}e^{-id_{j,k}\xi}.
\end{equation*}

Then
\begin{equation*}
	P_k=\sum_{r=1}^{m}S_{r,k}.
\end{equation*}

After passing to another subsequence, each cluster is either zero in $L^2(E)$ for every $k$, or nonzero in $L^2(E)$ for every $k$. We call the latter clusters active. For an inactive cluster, $G_{r,k}$ is an entire exponential sum which vanishes almost everywhere on the positive-measure set $E$; hence the identity theorem gives $G_{r,k}\equiv0$ on $\C$. Thus $S_{r,k}(\xi)=e^{-ic_{r,k}\xi}G_{r,k}(\xi)$ vanishes identically. In particular, applying~\eqref{AC-exponential-identity} to the cluster alone, with $m=|C_r|$ and all reference coefficients equal to zero, gives
\[
\Big\|\sum_{j\in C_r}\beta_{j,k}\mathcal U(b_{j,k})V\Big\|_{\hs}^2
=\int_{\R}w(\xi)|S_{r,k}(\xi)|^2\,\mathrm{d}\xi=0.
\]

Hence an inactive cluster represents the zero element of $\Hs$ and may be discarded from the sum of bubbles defining $\sigma_k$. Therefore, for an inactive cluster in $P_k$, one has $S_{r,k}\equiv0$. For an active cluster, put
\[
A_{r,k}\coloneqq\|G_{r,k}\|_{L^2(E)}>0,
\qquad
g_{r,k}\coloneqq A_{r,k}^{-1}G_{r,k}.
\]

By~\eqref{AC-cluster-bounded} and
Lemma~\ref{lem:AC-exponential-compactness}, applied with $m=|C_r|$ and
$R=R_r$, the family $(g_{r,k})_k$ is relatively compact in $L^2(E)$,
and hence also in $L^1(E)$, since $E$ has finite measure.

\medskip
\noindent\textbf{Step 2.}
We claim that for distinct active clusters $r,\ell$,
\begin{equation}\label{AC-cluster-orthogonality}
	\langle S_{r,k},S_{\ell,k}\rangle_{L^2(E)}
	=o(1)A_{r,k}A_{\ell,k}.
\end{equation}

Indeed, after dividing by $A_{r,k}A_{\ell,k}$, the inner product is
\[
\int_Ee^{-i(c_{r,k}-c_{\ell,k})\xi}
g_{r,k}(\xi)\overline{g_{\ell,k}(\xi)}\,\mathrm{d}\xi.
\]

The products $g_{r,k}\overline{g_{\ell,k}}$ form a relatively compact family in $L^1(E)$: this follows from the relative compactness in $L^2(E)$ and the continuity of the map $(f,g)\mapsto f\overline g$ from $L^2(E)\times L^2(E)$ to $L^1(E)$. Thus~\eqref{AC-uniform-RL} and~\eqref{AC-centres-separated} prove~\eqref{AC-cluster-orthogonality}.

Since the inner product is Hermitian, the sum of the off-diagonal terms in the expansion of $\|P_k\|_{L^2(E)}^2$ is the sum over ordered pairs $r\neq\ell$. Since there are finitely many clusters, there is a sequence $\eta_k\to0$ such that, for every distinct active pair,
\[
|\langle S_{r,k},S_{\ell,k}\rangle_{L^2(E)}|
\le\eta_k A_{r,k}A_{\ell,k}.
\]

Consequently, if $\mathcal A$ denotes the set of active clusters, then
\begin{align*}
	\|P_k\|_{L^2(E)}^2
	\ge \sum_{r\in\mathcal A}A_{r,k}^2
	-\eta_k\sum_{\substack{r,\ell\in\mathcal A\\r\neq\ell}}
	A_{r,k}A_{\ell,k}\ge \bigl(1-(m-1)\eta_k\bigr)
	\sum_{r\in\mathcal A}A_{r,k}^2.
\end{align*}

Here we used
\[
\sum_{\substack{r,\ell\in\mathcal A\\r\neq\ell}}
A_{r,k}A_{\ell,k}
\le (m-1)\sum_{r\in\mathcal A}A_{r,k}^2,
\]
which follows from Cauchy--Schwarz. For sufficiently large $k$, $1-(m-1)\eta_k\ge1/2$. Thus, for all sufficiently large $k$,
\[
\sum_{r\in\mathcal A}A_{r,k}^2
\le 2\|P_k\|_{L^2(E)}^2.
\]

In view of~\eqref{AC-P-bounded}, this gives
\begin{equation}\label{AC-cluster-norm-bounded}
	\sup_k A_{r,k}<\infty
	\qquad\text{for every active cluster }r.
\end{equation}

\medskip
\noindent\textbf{Step 3.}
Fix $i\in\{1,\ldots,\nu\}$. By~\eqref{AC-PQ-close},
\begin{equation}\label{AC-test-PQ}
	\langle P_k,e^{-ia_{i,k}\cdot}\rangle_{L^2(E)}
	-\langle Q_k,e^{-ia_{i,k}\cdot}\rangle_{L^2(E)}
	\longrightarrow0.
\end{equation}

On the other hand,
\begin{align*}
	\langle Q_k,e^{-ia_{i,k}\cdot}\rangle_{L^2(E)}
	=\widetilde\alpha_{i,k}|E|+\sum_{\ell\neq i}
	\widetilde\alpha_{\ell,k}
	\int_Ee^{-i(a_{\ell,k}-a_{i,k})\xi}\,\mathrm{d}\xi.
\end{align*}

By~\eqref{AC-reference-separation}, the frequencies in the sum tend to infinity in absolute value. The Riemann--Lebesgue lemma and~\eqref{AC-reference-amplitudes} therefore imply
\begin{equation}\label{AC-reference-test-limit}
	\langle Q_k,e^{-ia_{i,k}\cdot}\rangle_{L^2(E)}
	\longrightarrow |E|.
\end{equation}

Combining~\eqref{AC-test-PQ} and~\eqref{AC-reference-test-limit},
\begin{equation}\label{AC-P-test-limit}
	\langle P_k,e^{-ia_{i,k}\cdot}\rangle_{L^2(E)}
	\longrightarrow |E|>0.
\end{equation}

We claim that there is an active cluster $r$ for which $(c_{r,k}-a_{i,k})_k$ is bounded. If not, by the subsequence choice in Step~1, $|c_{r,k}-a_{i,k}|\to\infty$ for every active cluster. For such a cluster,
\[
\langle S_{r,k},e^{-ia_{i,k}\cdot}\rangle_{L^2(E)}
=A_{r,k}
\int_Ee^{-i(c_{r,k}-a_{i,k})\xi}g_{r,k}(\xi)\,\mathrm{d}\xi.
\]

The family $(g_{r,k})_k$ is relatively compact in $L^1(E)$, so
\eqref{AC-uniform-RL}, together with~\eqref{AC-cluster-norm-bounded}, shows that this expression tends to zero. Summing over the finitely many active clusters would give
\[
\langle P_k,e^{-ia_{i,k}\cdot}\rangle_{L^2(E)}\longrightarrow0,
\]
contradicting~\eqref{AC-P-test-limit}. This proves the claim.

\medskip
\noindent\textbf{Step 4.}
An active cluster cannot have bounded distance from two distinct reference scales. Indeed, if $|c_{r,k}-a_{i,k}|$ and $|c_{r,k}-a_{j,k}|$ were both bounded for $i\neq j$, then $|a_{i,k}-a_{j,k}|$ would be bounded, contrary to~\eqref{AC-reference-separation}.

It follows that the $\nu$ reference scales require at least $\nu$ distinct active clusters. Since there are at most $\nu$ candidate clusters (and there are exactly $\nu$ candidate indices in total), we necessarily have $m=\nu$, every cluster is active, and every cluster contains exactly one candidate index. After permuting the candidate indices, we may therefore assume
\begin{equation}\label{AC-matched-bounded-scales}
	|b_{i,k}-a_{i,k}|\le R_i
	\qquad(1\le i\le\nu)
\end{equation}
for fixed constants $R_i<\infty$. Moreover, for a singleton cluster, $A_{i,k}=|\beta_{i,k}|\,|E|^{1/2}$, and~\eqref{AC-cluster-norm-bounded} gives
\begin{equation}\label{AC-matched-bounded-amplitudes}
	\sup_k|\beta_{i,k}|<\infty
	\qquad(1\le i\le\nu).
\end{equation}

\medskip
\noindent\textbf{Step 5.}
By~\eqref{AC-matched-bounded-scales} and~\eqref{AC-matched-bounded-amplitudes}, after passing to a subsequence we may assume that, for every $i$,
\[
b_{i,k}-a_{i,k}\longrightarrow h_i\in\R,
\qquad
\beta_{i,k}\longrightarrow\beta_i\in\R.
\]

Set
\[
r_{i,k}
\coloneqq
\beta_{i,k}\mathcal U(b_{i,k}-a_{i,k})V
-\widetilde\alpha_{i,k}V.
\]

Strong continuity of the dilation group and~\eqref{AC-reference-amplitudes} give
\begin{equation}\label{AC-ri-strong-limit}
	r_{i,k}\longrightarrow
	r_i\coloneqq\beta_i\mathcal U(h_i)V-V
	\quad\text{strongly in }\Hs.
\end{equation}

We can write
\[
\sigma_k-\widetilde\sigma_k
=\sum_{i=1}^{\nu}\mathcal U(a_{i,k})r_{i,k}.
\]

Expanding its squared norm,
\begin{align}
	\|\sigma_k-\widetilde\sigma_k\|_{\hs}^2
	=\sum_{i=1}^{\nu}\|r_{i,k}\|_{\hs}^2+2\sum_{1\le i<j\le\nu}
	\left\langle
	\mathcal U(a_{i,k})r_{i,k},
	\mathcal U(a_{j,k})r_{j,k}
	\right\rangle_{\hs}.
	\label{AC-matched-norm-expansion}
\end{align}

We show that every cross term in~\eqref{AC-matched-norm-expansion} tends to zero. Fix $i\neq j$ and put $\tau_k\coloneqq a_{j,k}-a_{i,k}$. By unitarity and the group property,
\[
\left\langle
\mathcal U(a_{i,k})r_{i,k},
\mathcal U(a_{j,k})r_{j,k}
\right\rangle_{\hs}=\left\langle
r_{i,k},\mathcal U(\tau_k)r_{j,k}
\right\rangle_{\hs}.
\]

Adding and subtracting $r_i$ in the first entry and $r_j$ in the second entry, we obtain
\begin{align*}
&\left\langle r_{i,k},\mathcal U(\tau_k)r_{j,k}\right\rangle_{\hs}
-\left\langle r_i,\mathcal U(\tau_k)r_j\right\rangle_{\hs}\\
&=\left\langle r_{i,k}-r_i,\mathcal U(\tau_k)r_{j,k}\right\rangle_{\hs}
+\left\langle r_i,\mathcal U(\tau_k)(r_{j,k}-r_j)\right\rangle_{\hs}.
\end{align*}

The Cauchy--Schwarz inequality and the unitarity of $\mathcal U(\tau_k)$ therefore give
\begin{align*}
&\left|
\left\langle r_{i,k},\mathcal U(\tau_k)r_{j,k}\right\rangle_{\hs}
-\left\langle r_i,\mathcal U(\tau_k)r_j\right\rangle_{\hs}
\right|\\
&\quad\le
\|r_{i,k}-r_i\|_{\hs}\,\|r_{j,k}\|_{\hs}
+\|r_i\|_{\hs}\,\|r_{j,k}-r_j\|_{\hs}.
\end{align*}

By~\eqref{AC-ri-strong-limit}, the two differences on the right tend to zero. The sequence $(r_{j,k})_k$ is bounded because it converges strongly to $r_j$. Hence the entire right-hand side tends to zero.

It remains to consider the fixed-vector term. By~\eqref{AC-reference-separation}, $|\tau_k|=|a_{j,k}-a_{i,k}|\ge L_k\to\infty$. Applying Lemma~\ref{lem:AC-dilation-weak} to the fixed vectors $r_i,r_j\in\Hs$ yields $\left\langle r_i,\mathcal U(\tau_k)r_j\right\rangle_{\hs}\to0$.

Combining the last two limits proves
\[
\left\langle
\mathcal U(a_{i,k})r_{i,k},
\mathcal U(a_{j,k})r_{j,k}
\right\rangle_{\hs}\longrightarrow0.
\]

There are only finitely many pairs $1\le i<j\le\nu$, so the entire off-diagonal sum in~\eqref{AC-matched-norm-expansion} tends to zero.

It follows from~\eqref{AC-global-matching-hypothesis} and
\eqref{AC-matched-norm-expansion} that
\[
0=\lim_{k\to\infty}
\|\sigma_k-\widetilde\sigma_k\|_{\hs}^2
=\sum_{i=1}^{\nu}\|r_i\|_{\hs}^2.
\]

Hence $r_i=0$ for every $i$, that is,
\begin{equation}\label{AC-single-bubble-identification}
	\beta_i\mathcal U(h_i)V=V.
\end{equation}

Taking norms in~\eqref{AC-single-bubble-identification} gives $|\beta_i|=1$. Since $V$ and all its dilates are strictly positive, the case $\beta_i=-1$ is impossible. Hence $\beta_i=1$.

It remains to show that $h_i=0$. If $\mathcal U(h_i)V=V$, Lemma~\ref{lem:AC-radial-log-isometry} with $q=\mathcal JV$, gives $q(\cdot-h_i)=q$ in $L^2(\R)$. If $h_i\neq0$, then $q$ is periodic with a nonzero period. A nonzero periodic function cannot belong to $L^2(\R)$: the integrals of $|q|^2$ over the disjoint intervals of length $|h_i|$ would all be equal, and their sum could be finite only if each were zero. Since $q\neq0$, this is impossible. Therefore $h_i=0$.

We have thus proved $\beta_{i,k}\to1$, and $b_{i,k}-a_{i,k}\to0$. Together with $\widetilde\alpha_{i,k}\to1$, this is~\eqref{AC-global-matching-conclusion}.\qed

\begin{Corollary}\label{cor:AC-uniform-matching}
	For every $\eta>0$, there exists $\epsilon_\eta>0$, depending only on $N,\nu,s,t,\eta$, such that the following holds. Let $0<\epsilon\le\epsilon_\eta$, and let $\widetilde\sigma=\sum_{i=1}^{\nu}\widetilde\alpha_i\mathcal U(a_i)V$ be a $2\epsilon$-interacting family. If $\sigma=\sum_{j=1}^{\nu}\beta_j\mathcal U(b_j)V\in\mathcal M_\nu$ satisfies $\|\sigma-\widetilde\sigma\|_{\Hs}\le2\epsilon$, then, after a permutation of the candidate indices,
	\begin{equation}\label{AC-uniform-matching-conclusion}
		\max_{1\le i\le\nu}
		\left(
		|\beta_i-\widetilde\alpha_i|+|b_i-a_i|
		\right)<\eta.
	\end{equation}
\end{Corollary}

\Proof
Suppose that the assertion is false for some $\eta>0$. Then there exist sequences $\epsilon_k\downarrow0$, $2\epsilon_k$-interacting reference families $\widetilde\sigma_k$, and competitors $\sigma_k\in\mathcal M_\nu$ such that $\|\sigma_k-\widetilde\sigma_k\|_{\hs}\le2\epsilon_k$, but~\eqref{AC-uniform-matching-conclusion} fails for every permutation. Since $2\epsilon_k\to0$, Lemma~\ref{lem:AC-global-matching} applies and produces permutations for which all the amplitude and logarithmic-scale differences tend to zero. This contradicts the failure of~\eqref{AC-uniform-matching-conclusion}.\qed

\subsection{Local nondegeneracy and strict convexity}
We now use the local tangent-space argument analogous to the finite-dimensional projection step in~\cite[Proposition~7]{Bahri-Coron-88}. The global matching result above is what permits us to apply this local argument to the infimum over the multi-bubble set $\mathcal M_\nu$.

\begin{Lemma}\label{lem:AC-uniform-differential}
	There exist constants $L_0>0$ and $c_0>0$, depending only on
	$N,\nu,s,t$, such that the following holds. Suppose that
	$|a_i-a_j|\ge L_0$ for $i\neq j$ and
	$\frac12\le\widetilde\alpha_i\le\frac32$. Define
	\begin{equation}\label{AC-local-chart}
		\Phi(h,\tau)
		\coloneqq
		\sum_{i=1}^{\nu}
		(\widetilde\alpha_i+h_i)
		\mathcal U(a_i+\tau_i)V,
		\qquad (h,\tau)\in\R^{2\nu}.
	\end{equation}
	Then
	\begin{equation}\label{AC-differential-lower-bound}
		\|D\Phi(0,0)[h,\tau]\|_{\Hs}
		\ge c_0\sqrt{|h|^2+|\tau|^2}.
	\end{equation}
\end{Lemma}

\Proof
Differentiating~\eqref{AC-local-chart} at the origin gives
\begin{equation}\label{AC-differential-formula}
	D\Phi(0,0)[h,\tau]
	=\sum_{i=1}^{\nu}
	\left(h_i\mathcal U(a_i)V
	+\widetilde\alpha_i\tau_i\mathcal U(a_i)Z\right).
\end{equation}

By~\eqref{AC-VZ-orthogonal} and the unitarity of $\mathcal U(a)$, $\langle\mathcal U(a_i)V,\mathcal U(a_i)Z\rangle_{\hs}=0$. Thus the diagonal terms in the square of the norm in~\eqref{AC-differential-formula} equal $\sum_{i=1}^{\nu}\Big(h_i^2\|V\|_{\hs}^2+\widetilde\alpha_i^2\tau_i^2\|Z\|_{\hs}^2\Big)$.

Set
\[
m_*^2\coloneqq
\min\left\{\|V\|_{\hs}^2,\,\tfrac14\|Z\|_{\hs}^2\right\}>0.
\]

The diagonal contribution is therefore at least
$m_*^2(|h|^2+|\tau|^2)$. Every off-diagonal term is an inner product of the form
\[
\langle\mathcal U(a_i)X,\mathcal U(a_j)Y\rangle_{\hs}
=\langle X,\mathcal U(a_j-a_i)Y\rangle_{\hs},
\qquad X,Y\in\{V,Z\}.
\]

By Lemma~\ref{lem:AC-dilation-weak},
\[
\omega(L)\coloneqq
\max_{X,Y\in\{V,Z\}}
\sup_{|a|\ge L}
|\langle X,\mathcal U(a)Y\rangle_{\hs}|
\longrightarrow0
\quad\text{as }L\to\infty.
\]

Using $2|xy|\le x^2+y^2$ and the bounds on $\widetilde\alpha_i$, the sum of all off-diagonal terms is bounded in absolute value by $C_\nu\omega(L_0)(|h|^2+|\tau|^2)$, where $C_\nu$ depends only on $\nu$. Choose $L_0$ sufficiently large that $C_\nu\omega(L_0)\le m_*^2/2$. We then obtain
\[
\|D\Phi(0,0)[h,\tau]\|_{\hs}^2
\ge\frac{m_*^2}{2}(|h|^2+|\tau|^2).
\]

The desired estimate~\eqref{AC-differential-lower-bound} follows with $c_0=m_*/\sqrt2$.\qed

\begin{Lemma}\label{lem:AC-chart-estimates}
	There exist $r>0$, $M>0$, and $K>0$, depending only on $N,\nu,s,t$, with the following properties whenever the hypotheses of Lemma~\ref{lem:AC-uniform-differential} hold.
	\begin{enumerate}
		\item[(1)] For every $\theta=(h,\tau)\in\R^{2\nu}$ with $|\theta|\le r$ and every
		$\xi\in\R^{2\nu}$,
		\begin{equation}\label{AC-chart-differential-bound}
			\|D\Phi(\theta)[\xi]\|_{\Hs}
			\ge\frac{c_0}{2}|\xi|.
		\end{equation}
		\item[(2)] For every $\theta,\eta\in\overline B_r(0)$,
		\begin{equation}\label{AC-chart-inverse-bound}
			\|\Phi(\theta)-\Phi(\eta)\|_{\Hs}
			\ge\frac{c_0}{2}|\theta-\eta|.
		\end{equation}
		\item[(3)] For every $|\theta|\le r$ and $\xi,\zeta\in\R^{2\nu}$,
		\begin{equation}\label{AC-chart-second-derivative}
			\|D^2\Phi(\theta)[\xi,\zeta]\|_{\Hs}
			\le M|\xi||\zeta|,
		\end{equation}
		while
		\begin{equation}\label{AC-chart-first-derivative}
			\|D\Phi(\theta)\|_{\mathrm{op}}\le K.
		\end{equation}
	\end{enumerate}
	The constants are independent of the reference scales $a_i$.
\end{Lemma}

\Proof
Write $\theta=(h,\tau)$, $\xi=(v,w)$, and $\zeta=(\widehat v,\widehat w)$. Since the orbit $a\mapsto\mathcal U(a)V$ is $C^2$ and~\eqref{AC-orbit-derivatives} holds, direct differentiation gives
\begin{align}
	D\Phi(\theta)[v,w]
	&=\sum_{i=1}^{\nu}
	\left(v_i\mathcal U(a_i+\tau_i)V
	+(\widetilde\alpha_i+h_i)w_i
	\mathcal U(a_i+\tau_i)Z\right),
	\label{AC-chart-first-formula}\\
	D^2\Phi(\theta)[\xi,\zeta]
	&=\sum_{i=1}^{\nu}\Big[
	(v_i\widehat w_i+\widehat v_iw_i)
	\mathcal U(a_i+\tau_i)Z
	+(\widetilde\alpha_i+h_i)w_i\widehat w_i
	\mathcal U(a_i+\tau_i)\widetilde Z
	\Big].
	\label{AC-chart-second-formula}
\end{align}

Since $|\widetilde\alpha_i-1|\le1/2$, we have $|\widetilde\alpha_i|\le3/2$. For $\theta=(h,\tau)\in\overline B_{1/2}(0)$, one has $|h_i|\le|h|\le|\theta|\le1/2$, and hence $|\widetilde\alpha_i+h_i|\le\frac32+\frac12=2$. Thus the coefficients $|\widetilde\alpha_i+h_i|$ are uniformly bounded on the fixed ball $\overline B_{1/2}(0)$. Since $\mathcal U(a)$ is an isometry, the norms of $\mathcal U(a)V$, $\mathcal U(a)Z$, and $\mathcal U(a)\widetilde Z$ are independent of $a$. Applying Cauchy--Schwarz to the finite sums in~\eqref{AC-chart-first-formula}--\eqref{AC-chart-second-formula} gives uniform constants $K,M$ satisfying~\eqref{AC-chart-first-derivative} and~\eqref{AC-chart-second-derivative} on $\overline B_{1/2}(0)$. 

By the fundamental theorem of calculus and~\eqref{AC-chart-second-derivative}, $\|D\Phi(\theta)-D\Phi(0)\|_{\mathrm{op}}\le M|\theta|$. Combining this with Lemma~\ref{lem:AC-uniform-differential},
\[
\|D\Phi(\theta)[\xi]\|_{\hs}
\ge(c_0-Mr)|\xi|.
\]

Choose $r\le1/2$ small enough that $Mr\le c_0/2$. This proves~\eqref{AC-chart-differential-bound}. For $\theta,\eta\in\overline B_r(0)$, the line segment joining them is contained in $\overline B_r(0)$. Hence
\begin{align*}
	\Phi(\theta)-\Phi(\eta)
	=D\Phi(0)[\theta-\eta]+\int_0^1
	\big(D\Phi(\eta+q(\theta-\eta))-D\Phi(0)\big)
	[\theta-\eta] \,\mathrm{d}q.
\end{align*}

The norm of the Taylor remainder in the last expression is bounded by $Mr|\theta-\eta|$. Applying Lemma~\ref{lem:AC-uniform-differential} and reducing $r$ if necessary gives~\eqref{AC-chart-inverse-bound}.\qed

\begin{Proposition}\label{prop:AC-strict-convexity}
	There exist $r>0$ and $\epsilon_1>0$, depending only on
	$N,\nu,s,t$, such that the following holds. Let
	$0<\epsilon\le\epsilon_1$, let $u\in U(\epsilon)$, and choose a
	reference family
	\[
	\widetilde\sigma
	=\sum_{i=1}^{\nu}\widetilde\alpha_iV_{s,t}^{\widetilde\lambda_i}
	=\sum_{i=1}^{\nu}\widetilde\alpha_i\mathcal U(a_i)V,
	\qquad a_i=\log\widetilde\lambda_i,
	\]
	provided by the definition~\ref{def:tubular-nhbd} of $U(\epsilon)$. With $\Phi$ as in
	\eqref{AC-local-chart}, the function
	\begin{equation}\label{AC-distance-functional}
		F_u(\theta)\coloneqq
		\frac12\|u-\Phi(\theta)\|_{\Hs}^2
	\end{equation}
	is uniformly strictly convex on $\overline B_r(0)$. More precisely,
	\begin{equation}\label{AC-strict-convexity-estimate}
		D^2F_u(\theta)[\xi,\xi]
		\ge\frac{c_0^2}{8}|\xi|^2
		\qquad(|\theta|\le r,\ \xi\in\R^{2\nu}).
	\end{equation}
\end{Proposition}

\Proof
We first decrease the radius $r$ in Lemma~\ref{lem:AC-chart-estimates}, if necessary, so that $MKr\le\frac{c_0^2}{16}$. We then choose
$\epsilon_1$ sufficiently small that
\[
\log\frac1{2\epsilon_1}\ge L_0,
\qquad 2\epsilon_1\le\frac12,
\qquad M\epsilon_1\le\frac{c_0^2}{16}.
\]

Then the hypotheses of Lemmas~\ref{lem:AC-uniform-differential} and~\ref{lem:AC-chart-estimates} hold for the reference family. In particular, \eqref{AC-chart-inverse-bound} and~\eqref{AC-chart-first-derivative} hold with uniform constants.

For $\xi\in\R^{2\nu}$, differentiation of~\eqref{AC-distance-functional}
gives
\begin{equation*}
	D^2F_u(\theta)[\xi,\xi]
	=\|D\Phi(\theta)[\xi]\|_{\hs}^2
	-\langle u-\Phi(\theta),
	D^2\Phi(\theta)[\xi,\xi]\rangle_{\hs}.
\end{equation*}

By~\eqref{AC-chart-differential-bound},
\[
\|D\Phi(\theta)[\xi]\|_{\hs}^2
\ge\frac{c_0^2}{4}|\xi|^2.
\]

Furthermore, by~\eqref{AC-chart-second-derivative},
\[
\left|
\langle u-\Phi(\theta),
D^2\Phi(\theta)[\xi,\xi]\rangle_{\hs}
\right|
\le M\|u-\Phi(\theta)\|_{\hs}|\xi|^2.
\]

Since $u\in U(\epsilon)$ and $\|u-\Phi(0)\|_{\Hs}<\epsilon$, we also have
\[
\|u-\Phi(\theta)\|_{\hs}
\le\|u-\Phi(0)\|_{\hs}
+\|\Phi(0)-\Phi(\theta)\|_{\hs}
<\epsilon+K r.
\]

Consequently,
\[
D^2F_u(\theta)[\xi,\xi]
\ge\left(\frac{c_0^2}{4}-M\epsilon-MKr\right)|\xi|^2.
\]

According to the choices of $r,\epsilon_1$ fixed above, this yields~\eqref{AC-strict-convexity-estimate}.\qed

\subsection{Existence and uniqueness of the nearest point}
We now combine the global matching result of Lemma~\ref{lem:AC-global-matching} and its uniform consequence, Corollary~\ref{cor:AC-uniform-matching}, with the local strict-convexity estimate of Proposition~\ref{prop:AC-strict-convexity}. The matching result places every sufficiently good competitor in a common local parameter chart, while strict convexity yields uniqueness in that chart. This is the step that identifies the minimizer obtained in a local parameter chart with the infimum over the full set $\mathcal M_\nu$ in~\eqref{eq:min-problem}.

\begin{Proposition}\label{prop:AC-nearest-point}
	There exists $\epsilon_0>0$, depending only on $N,\nu,s,t$, such that for
	every $0<\epsilon\le\epsilon_0$ and every $u\in U(\epsilon)$, the
	minimization problem~\eqref{eq:min-problem} admits a unique minimizing
	function $\sigma\in\mathcal M_\nu$, where uniqueness of the bubble
	parameters is understood up to permutation of the indices. Moreover,
	\begin{equation}\label{AC-remainder-bound}
		\|u-\sigma\|_{\Hs}<\epsilon.
	\end{equation}
	If
	\[
	u=\sigma+\rho,
	\qquad
	\sigma=\sum_{i=1}^{\nu}\alpha_iV_{s,t}^{\lambda_i},
	\]
	and
	\[
	\widetilde\sigma
	=\sum_{i=1}^{\nu}
	\widetilde\alpha_iV_{s,t}^{\widetilde\lambda_i}
	\]
	is any reference family supplied by $u\in U(\epsilon)$, then, after a
	permutation,
	\begin{equation}\label{AC-parameter-stability}
		|\alpha_i-\widetilde\alpha_i|
		+\left|\log\frac{\lambda_i}{\widetilde\lambda_i}\right|
		\le C\epsilon,
		\qquad 1\le i\le\nu,
	\end{equation}
	where $C$ depends only on $N,\nu,s,t$. In particular,
	\begin{equation}\label{AC-interaction-stability}
		\max_i|\alpha_i-1|\le C\epsilon,
		\qquad
		\min\left\{\frac{\lambda_i}{\lambda_j},
		\frac{\lambda_j}{\lambda_i}\right\}
		\le C\epsilon
		\quad(i\neq j).
	\end{equation}
	Finally, the remainder is orthogonal to the tangent space at $\sigma$:
	\begin{equation}\label{AC-tangent-orthogonality}
		\langle\rho,V_{s,t}^{\lambda_i}\rangle_{\Hs}=0,
		\qquad
		\langle\rho,\dot V_{s,t}^{\lambda_i}\rangle_{\Hs}=0,
		\qquad 1\le i\le\nu.
	\end{equation}
\end{Proposition}

\Proof Fix $u\in U(\epsilon)$ and choose a reference family
\[
\widetilde\sigma
=\sum_{i=1}^{\nu}
\widetilde\alpha_i\mathcal U(a_i)V,
\qquad a_i=\log\widetilde\lambda_i,
\]
such that
\begin{equation}\label{AC-reference-distance}
	\|u-\widetilde\sigma\|_{\hs}<\epsilon.
\end{equation}

In particular, since $\widetilde\sigma\in\mathcal M_\nu$,
\begin{equation}\label{AC-inf-less-epsilon}
	d\coloneqq\inf_{\sigma\in\mathcal M_\nu}
	\|u-\sigma\|_{\hs}<\epsilon.
\end{equation}

Let $r$ and $\epsilon_1$ be as in
Proposition~\ref{prop:AC-strict-convexity}. Choose $\eta>0$ so small that
\begin{equation}\label{AC-eta-choice}
	\sqrt{\nu}\,\eta<\frac r2.
\end{equation}

Let $\epsilon_2>0$ be the corresponding constant in Corollary~\ref{cor:AC-uniform-matching} for this choice of $\eta$. Set $\epsilon_0\coloneqq\min\{\epsilon_1,\epsilon_2\}$.

\medskip
\noindent\textit{Existence.}
Choose a minimizing sequence
\[
\sigma_k=\sum_{j=1}^{\nu}
\beta_{j,k}\mathcal U(b_{j,k})V\in\mathcal M_\nu
\]
such that
\[
\|u-\sigma_k\|_{\hs}\longrightarrow d.
\]

By~\eqref{AC-reference-distance} and~\eqref{AC-inf-less-epsilon},
\[
\limsup_{k\to\infty}
\|\sigma_k-\widetilde\sigma\|_{\hs}
\le d+\|u-\widetilde\sigma\|_{\hs}<2\epsilon.
\]

In particular, for all sufficiently large $k$, $\|\sigma_k-\widetilde\sigma\|_{\hs}\le2\epsilon$ and therefore, Corollary~\ref{cor:AC-uniform-matching} applies. After permuting the indices of each $\sigma_k$, its parameters can be written as
\[
\beta_{i,k}=\widetilde\alpha_i+h_{i,k},
\qquad
b_{i,k}=a_i+\tau_{i,k},
\]
where by~\eqref{AC-uniform-matching-conclusion}, $|h_{i,k}|+|\tau_{i,k}|<\eta$ for every $i$. Therefore, we have \begin{equation*}
	|(h_k,\tau_k)|^2
	=\sum_{i=1}^{\nu}(h_{i,k}^2+\tau_{i,k}^2)
	<\nu\eta^2<\frac{r^2}{4},
\end{equation*}
where the last inequality follows from~\eqref{AC-eta-choice}. Thus the parameter vectors lie in the compact set $\overline B_{\sqrt{\nu}\eta}(0)\subset B_{r/2}(0)$. After passing to a subsequence,
\[
(h_k,\tau_k)\longrightarrow(h_*,\tau_*)
\in\overline B_{\sqrt{\nu}\eta}(0)
\subset B_{r/2}(0).
\]

By continuity of $\Phi$ as an $\Hs$-valued map,
\[
\sigma_k=\Phi(h_k,\tau_k)
\longrightarrow\Phi(h_*,\tau_*)\eqqcolon\sigma_*
\quad\text{strongly in }\Hs.
\]

Since
\[
\sigma_*
=\sum_{i=1}^{\nu}
(\widetilde\alpha_i+h_{*,i})
\mathcal U(a_i+\tau_{*,i})V,
\]
all its coefficients are real and all its scales are positive. Hence
$\sigma_*\in\mathcal M_\nu$. By continuity of the norm,
\[
\|u-\sigma_*\|_{\hs}
=\lim_{k\to\infty}\|u-\sigma_k\|_{\hs}=d.
\]

Thus the infimum in~\eqref{eq:min-problem} is attained.

\medskip
\noindent\textit{Uniqueness.}
Suppose that $\sigma_1,\sigma_2\in\mathcal M_\nu$ are two minimizers. By~\eqref{AC-inf-less-epsilon} and~\eqref{AC-reference-distance},
\[
\|\sigma_j-\widetilde\sigma\|_{\hs}
\le\|\sigma_j-u\|_{\hs}
+\|u-\widetilde\sigma\|_{\hs}<2\epsilon,
\qquad j=1,2.
\]

Corollary~\ref{cor:AC-uniform-matching} therefore gives, after
permutations,
\[
\sigma_1=\Phi(\theta_1),
\qquad
\sigma_2=\Phi(\theta_2),
\qquad
|\theta_1|,|\theta_2|<\frac r2.
\]

If $\theta_1\neq\theta_2$, strict convexity from Proposition~\ref{prop:AC-strict-convexity} gives
\[
F_u\left(\frac{\theta_1+\theta_2}{2}\right)
<\frac12F_u(\theta_1)+\frac12F_u(\theta_2)=\frac{d^2}{2}.
\]

But $\Phi((\theta_1+\theta_2)/2)\in\mathcal M_\nu$, contradicting the definition of $d$. Thus $\theta_1=\theta_2$ after permutation, and $\sigma_1=\sigma_2$. If the same minimizing function admitted another representation by $\nu$ bubbles, that representation would satisfy the same estimate $\|\sigma_*-\widetilde\sigma\|_{\hs}<2\epsilon$. Hence Corollary~\ref{cor:AC-uniform-matching} places that representation in the same chart after a permutation. The local inverse estimate~\eqref{AC-chart-inverse-bound} then shows that the two parameter vectors coincide. Thus the minimizing parameters are unique up to permutation.

\medskip
\noindent\textit{Distance and parameter estimates.}
Since $\sigma_*$ is a minimizer and $\widetilde\sigma$ is an admissible competitor,
\[
\|u-\sigma_*\|_{\hs}
\le\|u-\widetilde\sigma\|_{\hs}<\epsilon,
\]
which proves~\eqref{AC-remainder-bound}.

Write $\sigma_*=\Phi(\theta_*)$ with $\theta_*=(h_*,\tau_*)$. By minimality and the triangle inequality,
\begin{align*}
	\|\Phi(\theta_*)-\Phi(0)\|_{\hs}
	\le\|\Phi(\theta_*)-u\|_{\hs}
	+\|u-\Phi(0)\|_{\hs}\le2\|u-\Phi(0)\|_{\hs}<2\epsilon.
\end{align*}

Applying~\eqref{AC-chart-inverse-bound} with $\eta=0$ gives
\[
\frac{c_0}{2}|\theta_*|
\le\|\Phi(\theta_*)-\Phi(0)\|_{\hs}<2\epsilon.
\]

Therefore
\begin{equation}\label{AC-theta-bound}
	|\theta_*|<\frac4{c_0}\epsilon.
\end{equation}

In particular,
\[
|\alpha_i-\widetilde\alpha_i|
+\left|\log\frac{\lambda_i}{\widetilde\lambda_i}\right|
\le C\epsilon,
\]
which proves~\eqref{AC-parameter-stability}.

Since the reference family is $2\epsilon$-interacting,
\[
|\widetilde\alpha_i-1|\le2\epsilon,
\qquad
\min\left\{
\frac{\widetilde\lambda_i}{\widetilde\lambda_j},
\frac{\widetilde\lambda_j}{\widetilde\lambda_i}
\right\}\le2\epsilon.
\]

The amplitude estimate follows from
\[
|\alpha_i-1|
\le|\alpha_i-\widetilde\alpha_i|
+|\widetilde\alpha_i-1|
\le C\epsilon.
\]

For the scales, write
$\tau_i^*=\log(\lambda_i/\widetilde\lambda_i)$. Then
\begin{align*}
	\Big|\log\frac{\lambda_i}{\lambda_j}\Big|
	\ge\Big|\log\frac{\widetilde\lambda_i}{\widetilde\lambda_j}\Big|-|\tau_i^*|-|\tau_j^*|\ge\log\frac1{2\epsilon}-C\epsilon.
\end{align*}

Exponentiating, and reducing $\epsilon_0$ if necessary, gives
\[
\min\left\{
\frac{\lambda_i}{\lambda_j},
\frac{\lambda_j}{\lambda_i}
\right\}
\le2\epsilon e^{C\epsilon}
\le C'\epsilon.
\]

This proves~\eqref{AC-interaction-stability}.

\medskip
\noindent\textit{Orthogonality of the remainder.}
We reduce $\epsilon_0$ once more, if necessary, so that $2\epsilon_0+4\epsilon_0/c_0<1/2$. Since $|\widetilde\alpha_i-1|\le2\epsilon$ and $|h_{*,i}|\le|\theta_*|<4\epsilon/c_0$ by~\eqref{AC-theta-bound}, this guarantees $\alpha_i>1/2$ for every $i$. The existence argument gives $\theta_*=(h_*,\tau_*)\in\overline B_{r/2}(0)\subset B_r(0)$. Hence $\theta_*$ is an interior point of the local chart, and the first variations of $F_u$ vanish at $\theta_*$. Writing
\[
\rho=u-\sigma_*,
\qquad
\sigma_*=
\sum_{i=1}^{\nu}\alpha_i
\mathcal U(a_i+\tau_i^*)V,
\]
we differentiate $F_u$ with respect to $h_i$ to obtain $0=\partial_{h_i}F_u(\theta_*)=-\langle\rho,\mathcal U(a_i+\tau_i^*)V\rangle_{\hs}$. Consequently, $\langle\rho,V_{s,t}^{\lambda_i}\rangle_{\hs}=0$. 

Differentiating $F_u$ with respect to $\tau_i$ gives $0=\partial_{\tau_i}F_u(\theta_*)=-\alpha_i\langle\rho,\mathcal U(a_i+\tau_i^*)Z\rangle_{\hs}$. Now by~\eqref{AC-log-derivative} we have $\mathcal U(a_i+\tau_i^*)Z=\lambda_i\dot V_{s,t}^{\lambda_i}$. Since $\alpha_i>1/2,\lambda_i>0$, we obtain $\langle\rho,\dot V_{s,t}^{\lambda_i}\rangle_{\hs}=0$. 

This completes the proof of~\eqref{AC-tangent-orthogonality}, which is precisely the orthogonality conditions~\eqref{Inn-bub}--\eqref{Inn-der-bub} used in the proof of Theorem~\ref{quantitative-stability}.\qed

\begin{Lemma}
	The minimizing parameters depend locally Lipschitz-continuously on $u$ as
	long as the reference family is fixed. More precisely, fix
	$0<\epsilon\le\epsilon_0$, where $\epsilon_0$ is as in
	Proposition~\ref{prop:AC-nearest-point}, and fix a $2\epsilon$-interacting
	reference configuration $\widetilde\sigma$. Let $\Phi$ be the chart
	associated with this same reference configuration. Suppose that
	$\|u-\widetilde\sigma\|_{\Hs}<\epsilon$ and
	$\|v-\widetilde\sigma\|_{\Hs}<\epsilon$. Then
	$u,v\in U(\epsilon)$ with this same reference configuration. By
	Proposition~\ref{prop:AC-nearest-point}, their minimizing parameters
	$\theta_u,\theta_v$ belong to $B_{r/2}(0)$, and
	\begin{equation*}
		|\theta_u-\theta_v|
		\le C\|u-v\|_{\Hs}.
	\end{equation*}
\end{Lemma}

\Proof
For $w\in\{u,v\}$, let $F_w(\theta)=\frac12\|w-\Phi(\theta)\|_{\hs}^2$. The minimizing parameters satisfy
\[
\nabla F_u(\theta_u)=0,
\qquad
\nabla F_v(\theta_v)=0.
\]

Put $d\coloneqq\theta_u-\theta_v$. The segment $\theta_v+q d$, $0\le q\le1$, is contained in $B_{r/2}(0)$, which is convex. The fundamental theorem of calculus and~\eqref{AC-strict-convexity-estimate} give
\begin{align*}
\left\langle
\nabla F_u(\theta_u)-\nabla F_u(\theta_v),d
\right\rangle
=\int_0^1D^2F_u(\theta_v+qd)[d,d]\,\mathrm dq\ge\tfrac{c_0^2}{8}|d|^2.
\end{align*}

Since $\nabla F_u(\theta_u)=0$, the left-hand side equals $-\langle\nabla F_u(\theta_v),d\rangle$. Since $\nabla F_v(\theta_v)=0$ and $\nabla F_u(\theta_v)-\nabla F_v(\theta_v)=D\Phi(\theta_v)^*(v-u)$, we obtain, using~\eqref{AC-chart-first-derivative},\[\tfrac{c_0^2}{8}|\theta_u-\theta_v|^2\le K\|u-v\|_{\hs}|\theta_u-\theta_v|.\]

The desired estimate follows by division when the two parameter vectors are different otherwise it is immediate.\qed

\section*{Acknowledgments} The first author thanks TIFR-CAM, Bengaluru, for postdoctoral support and an excellent research environment, and Debabrata Karmakar for insightful discussions and helpful suggestions. The second author thanks the Department of Mathematics, IIT Bombay, for providing ideal working conditions, and expresses his gratitude to Chandan Biswas for illuminating discussions related to the commutator estimate in Appendix~\ref{AA}.

\bibliographystyle{acm}
\bibliography{bibliography}
\end{document}